\documentclass[11pt,a4paper]{scrartcl}
\usepackage[round]{natbib}

\usepackage[top=3cm, bottom=3cm, left=3cm, right=3cm]{geometry}

\usepackage{booktabs}

\usepackage{color}
\usepackage{latexsym}              
\usepackage{amsmath}               
\usepackage{amssymb}               
\usepackage{amsfonts}              
\usepackage{amsthm}                
\usepackage{multirow}
\usepackage{tikz}                  
\usetikzlibrary{arrows,positioning,shapes}
\usepackage{tcolorbox}
\usepackage{placeins}
\usepackage{subcaption}
\usepackage{algorithm}
\usepackage{bm}

\usepackage[english]{babel} 

\RequirePackage[%
  pdfstartview=FitH,%
  breaklinks=true,%
  bookmarks=true,%
  colorlinks=true,%
  linkcolor= blue,
  anchorcolor=blue,%
  citecolor=orange,
  filecolor=blue,%
  menucolor=blue,%
  urlcolor=blue%
  ]{hyperref}
  \usepackage[capitalize,nameinlink,noabbrev]{cleveref}

\renewenvironment{proof}[1][\proofname]{{\bfseries #1.}}{\qed \\ }

\makeatother

\theoremstyle{plain}  
\newtheorem{theorem}{Theorem}[section]
\newtheorem{definition}[theorem]{Definition}

\newtheorem{lemma}[theorem]{Lemma}
\newtheorem{proposition}[theorem]{Proposition}

\newtheorem{corollary}[theorem]{Corollary}

\newcommand{\bbE}{\mathbb{E}}

\newcommand{\bbR}{\mathbb{R}}

\newcommand{\bbV}{\mathbb{V}}

\newcommand{\Real}{{\mathbb R}}

\DeclareMathOperator{\gammadist}{Gamma}

\DeclareMathOperator{\CRM}{CRM}

\DeclareMathOperator{\rhoca}{\rho_{mSt}}
\newcommand{\psica}{\psi_{\mathrm{mSt}}}
\newcommand{\psicainv}{\psi^{-1}_{\mathrm{mSt}}}
\DeclareMathOperator{\kappaca}{\kappa_{mSt}}
\DeclareMathOperator{\rhomix}{\rho_{mGG}}
\DeclareMathOperator{\psimix}{\psi_{mGG}}

\DeclareMathOperator{\kappamix}{\kappa_{mGG}}

\DeclareMathOperator{\rhoSt}{\rho_{St}}
\DeclareMathOperator{\psiSt}{\psi_{St}}
\DeclareMathOperator{\rhoGG}{\rho_{GG}}
\DeclareMathOperator{\psiGG}{\psi_{GG}}

\DeclareMathOperator{\LambertW}{{LambertW}}

\DeclareMathOperator{\logit}{logit}

\newcommand{\1}[1]{\mathbf{1}_{\{#1\}}}

\def\[#1\]{\begin{align}#1\end{align}}

\begin{document}
\title{ Rapidly Varying Completely Random Measures for Modeling Extremely Sparse Networks}

\author{
  Valentin Kilian\textsuperscript{1} \quad
  Benjamin Guedj\textsuperscript{2,3} \quad
  François Caron\textsuperscript{1} \\
  \small
  \textsuperscript{1}Department of Statistics, University of Oxford \\
    \small
  \textsuperscript{2}Department of Computer Science, University College London \\
    \small
  \textsuperscript{3}Inria, France \\
    \small
  \texttt{kilian@stats.ox.ac.uk}, \quad
  \texttt{b.guedj@ucl.ac.uk}, \quad
  \texttt{caron@stats.ox.ac.uk}
}
\date{}

\maketitle

\begin{abstract}
Completely random measures (CRMs) are fundamental to Bayesian nonparametric models, with applications in clustering, feature allocation, and network analysis. A key quantity of interest is the Laplace exponent, whose asymptotic behavior determines how the random structures scale. When the Laplace exponent grows nearly linearly — known as rapid variation — the induced models exhibit approximately linear growth in the number of clusters, features, or edges with sample size or network nodes. This regime is especially relevant for modeling sparse networks, yet existing CRM constructions lack tractability under rapid variation. We address this by introducing a new class of CRMs with index of variation $\alpha\in(0,1]$, defined as mixtures of stable or generalized gamma processes. These models offer interpretable parameters, include well-known CRMs as limiting cases, and retain analytical tractability through a tractable Laplace exponent and simple size-biased representation. We analyze the asymptotic properties of this CRM class and apply it to the Caron–Fox framework for sparse graphs. The resulting models produce networks with near-linear edge growth, aligning with empirical evidence from large-scale networks. Additionally, we present efficient algorithms for simulation and posterior inference, demonstrating practical advantages through experiments on real-world sparse network datasets.
\end{abstract}



\section{Introduction} 
\label{sec:intro}
Over the past twenty-five years, completely random measures (CRMs) have emerged as fundamental components of modern Bayesian methodology, with applications in clustering, feature modeling, species discovery, survival analysis, finance, and network analysis~\citep{Brix1999,Regazzini2003,James2002,Nieto-Barajas2004,Cont2004,Thibaux2007,Lijoi2007,Teh2009,James2009,Lijoi2010,Broderick2012,Favaro2013,Cai2016,Caron2017,Camerlenghi2019}. A key quantity associated with a CRM $G$ defined on a space $\Theta$ is its Laplace exponent, the Bernstein function $\psi : [0, \infty) \to [0, \infty)$ given by
\begin{equation}
\psi(t)=-\log\mathbb{E}\left[  e^{-tG(\Theta)}\right]  ,\text{ }t\geq0.
\label{eq:LaplaceExponentGeneric}
\end{equation}
The asymptotic behavior of $\psi$ as $t$ grows to infinity determines critical properties of the CRM and of the random structures built upon it. In many cases of practical interest, $\psi$ exhibits regular variation at infinity: \begin{equation} \psi(t) \underset{t \rightarrow \infty}{\sim} t^{\alpha} \ell(t), \label{eq:LaplaceExponentAsympGeneric} \end{equation} for some $\alpha \in [0,1]$ and a slowly varying function $\ell(t)$, which is typically either a positive constant or a power of $\log(t)$.

This asymptotic form is central in determining the growth rate of observable quantities in CRM-based models. For example, in normalized CRM mixture models, latent feature models, or species sampling schemes, the number of clusters, features, or distinct species $K_n$ scales as $n^{\alpha} \ell(n)$ as the sample size $n \to \infty$~\citep{Pitman2003,Gnedin2007,Broderick2012}. Similarly, in the Caron–Fox family of sparse graph models, which are based on CRMs~\citep{Caron2017}, the number of nodes $N$ scales as $E^{(1+\alpha)/2} \ell(\sqrt{E})$, where $E$ denotes the number of edges~\citep{Caron2023}. 

The parameter $\alpha \in [0,1]$ is referred to as the \emph{index of variation}. The limiting cases $\alpha = 0$ and $\alpha = 1$ correspond to \emph{slow} and \emph{rapid} variation respectively~\citep{Gnedin2007}. Classical CRMs with favorable analytical properties -- such as conjugacy -- typically fall in the range $\alpha \in [0,1)$, and include the gamma, beta, stable, generalized gamma, and stable beta processes~\citep{Hougaard1986,Hjort1990,Brix1999,Lijoi2007,Teh2009}; we collect a summary in \cref{tab:examplesCRMs}.

Despite its importance, the extreme case $\alpha = 1$, corresponding to rapid variation, lacks a tractable CRM construction in the literature. In clustering and feature allocation models, it yields regimes where the number of clusters or features grows approximately linearly in $n$, modulo a logarithmic correction. Within the Caron–Fox framework for random graphs, the case $\alpha = 1$ produces networks in which the number of edges grows nearly linearly with the number of nodes—a property widely regarded as more consistent with empirical observations of real-world sparse networks~\cite[Figure 1.1]{VanDerHofstad2024}.

In this article, we introduce a novel class of completely random measures with index of variation $\alpha \in (0,1]$, thereby filling the gap in the existing literature for the rapid variation regime. The proposed CRM is constructed as a mixture of stable, or more generally, generalized gamma CRMs. It features interpretable parameters, includes the stable and generalized gamma processes as limiting cases, and admits a simple and tractable Laplace exponent. Moreover, it supports a size-biased construction that facilitates inference and simulation.

We apply this model to the construction of sparse random graphs. Within the Caron–Fox framework, the resulting graphs exhibit edge growth that is approximately linear in the number of nodes. We develop algorithms for both simulation and posterior inference, leveraging the mixture structure of the model, and demonstrate its empirical validity on two real-world sparse networks.

The remainder of the paper is structured as follows. In \cref{sec:model}, we introduce our novel CRM family and detail its construction. In \cref{sec:properties}, we study its theoretical properties and describe appropriate simulation methods. In \cref{sec:graphs}, we apply our model to network data, describe graph properties, and present an efficient approach to posterior inference. In \cref{sec:experiments}, we present numerical experiments on real-world networks, and we gather concluding remarks in \cref{sec:discussion}. We defer to the supplementary material proofs, implementation details, additional experiments, and some background on regular variation.

\begin{table}[t]
\begin{center}
\caption{Index of variation and slowly varying function associated to some CRMs.}
\label{tab:examplesCRMs}
\begin{tabular}
[c]{lll}%
CRM & $\alpha$ & $\ell_{\alpha}$\\
\hline
Gamma & 0 & $\log$\\
Generalized gamma & (0,1) & constant\\
Beta & 0 & $\log$\\
Stable beta & (0,1) & constant\\
Stable & (0,1) & constant\\
Mixed Stable (ours) & (0,1] & $\log^{-1}$\\
Mixed Generalized Gamma (ours) & (0,1] & $\log^{-1}$%
\end{tabular}
\end{center}
\end{table}

\paragraph{Asymptotic notation.} Throughout this article, we write $f(x) \underset{x \to a}{\sim} g(x)$ to denote that $\lim_{x \to a} f(x)/g(x) = 1$, where $a \in \mathbb{R} \cup \{-\infty, +\infty\}$. If $f(x)$ and $g(x)$ are random variables, the notation indicates that the equivalence holds almost surely. This usage of $\sim$ should not be confused with the standard distributional notation $X \sim \mathcal{D}$, which means that the random variable $X$ has distribution $\mathcal{D}$.

\section{Rapidly Varying Completely Random Measures} 
\label{sec:model}
\subsection{Background on homogeneous CRMs} 

An homogeneous CRM on some Polish space $\Theta$ with no deterministic
component nor fixed components takes the form \citep{Kingman1967,Kingman1993,Lijoi2010}
\begin{equation}
G=\sum_{j\geq1}W_{j}\delta_{\theta_{j}},%
\end{equation}
where $\{(W_{j},\theta_{j})\}_{j\geq1}$ are realisations of a Poisson point
process on $(0,\infty)\times\Theta$ with mean measure $\nu(\mathrm{d}w)H(\mathrm{d}\theta)$
where $H$ is some probability distribution on $\Theta$, and $\nu$ is a L\'{e}vy measure on $(0,\infty)$  that satisfies
\begin{equation}
\int_{0}^{\infty}(1-e^{-w})\nu(\mathrm{d}w)<\infty.
\end{equation}
We assume here that $\nu$ admits a density $\rho$ with respect to the Lebesgue measure,
with $\nu(\mathrm{d}w)=\rho(w)\mathrm{d}w$.  We write $G\sim \CRM(\rho,H)$. Note that, if $\int_0^\infty \rho(w)\mathrm{d}w=\infty$, then $0<G(\Theta)=\sum_{j\geq 1} W_j<\infty$ almost surely. The Laplace exponent is related to the L\'{e}vy
intensity $\rho$ through
\begin{equation}
\psi(t)=\int_{0}^{\infty}(1-e^{-wt})\rho(w)\mathrm{d}w.
\label{eq:psifunction}
\end{equation}
Two standard CRMs are the stable CRM, with L\'{e}vy intensity and Laplace exponent
\begin{align*}
\rhoSt(w;\alpha)   =\frac{\alpha}{\Gamma(1-\alpha)}w^{-1-\alpha},~~
\psiSt(t;\alpha)  =t^{\alpha}%
\end{align*}
where $\alpha\in(0,1)$ is the index of variation, and the generalized gamma (GG) CRM (also known as exponentially tilted stable CRM), whose L\'{e}vy intensity and Laplace exponent are given by
\begin{align*}
\rhoGG(w;\alpha,\beta) &= \frac{\alpha}{\Gamma(1-\alpha
)}w^{-1-\alpha}e^{-\beta w}, \\
 \psiGG(t;\alpha,\beta) &= \psiSt(t+\beta
;\alpha)-\psiSt(\beta;\alpha) =(t+\beta)^{\alpha}-\beta^{\alpha}%
\end{align*}
where $\alpha\in(0,1)$ is the index of variation and $\beta\geq0$. For $\beta=0$, one recovers the stable CRM.

\subsection{A novel class CRMs with rapid variation} 

We now introduce a novel class of CRMs whose index of variation $\alpha\in(0,1]$, including the rapid variation case $\alpha=1$. This class of CRMs is defined as a mixture of generalized gamma CRMs, where the mixing is done with respect to the index of variation. We will refer to them as mixed Generalized Gamma (mGG). For ease of presentation, we first present the canonical model, based on a mixture of Stable CRMs (mSt) with parameters $\alpha$ and $\tau$. We then present the full model, with five parameters, whose properties can be derived from those of the canonical model.

\paragraph{Canonical model: Mixed Stable.} 

The L\'evy intensity of the canonical model is obtained by mixing over the index of variation in a stable model. It is given by
\begin{align}
\rhoca(w;\alpha,\tau) &  =\frac{1}{\alpha-\tau}\int_{\tau}^{\alpha
}\rhoSt(w;s)\mathrm{d}s =\frac{1}{\alpha-\tau}\int_{\tau}^{\alpha}\frac{s}{\Gamma(1-s)}w^{-1-s}\mathrm{d}s\label{eq:rhoca}
\end{align}
where $0\leq\tau<\alpha\leq1$. Although the intensity $\rhoca$ does not have a closed-form expression, the corresponding Laplace exponent has a remarkably simple expression:
\begin{align}
\psica(t;\alpha,\tau) &  =\frac{1}{\alpha-\tau}\int_{\tau}^{\alpha
}\psiSt(t;s)\mathrm{d}s   =\left\{
\begin{array}
[c]{ll}%
\frac{t^{\alpha}-t^{\tau}}{(\alpha-\tau)\log t} & t\in(0,1)\cup(1,\infty),\\
1 & t=1,\\
0 & t=0.
\end{array}
\right.
\label{eq:psica}
\end{align}

In this model, $\alpha\in(0,1]$ is the index of variation, as will be shown in \cref{thm:AsymptoticsLevyLaplace}. Note that the function $\psi$ is monotonically increasing and
continuous on $[0,\infty)$. Additionally, for any $\alpha\in(0,1)$,
\[
\lim_{\tau\rightarrow\alpha}\rhoca(w;\alpha,\tau)=\rho_{\text{St}%
}(w;\alpha),
\]
hence the stable model is recovered as a limiting case. Of the values for
$\tau$, $\alpha$, we will be particularly interested in the rapid variation
case $\alpha=1$, $\tau=0$. In this case, the Laplace exponent has a particularly tractable form, given by
\[
\psica(t;1,0)=\left\{
\begin{array}
[c]{ll}%
\frac{t-1}{\log t} & t>0\text{, }t\neq1,\\
1 & t=1,\\
0 & t=0.
\end{array}
\right. 
\]

\paragraph{Full model: Mixed Generalized Gamma.} 

The full five-parameters model is obtained by scaling and exponentially tilting the density of the mixed Stable CRM. Its L\'evy intensity can alternatively be expressed as a mixture of Generalized Gamma CRMs, hence the name mixed Generalized Gamma (mGG). Its L\'{e}vy intensity and Laplace exponent are given by
\begin{align}
\rhomix(w;\alpha,\tau,\beta,c,\eta) &  =\frac{\eta}{c}\rhoca\left (\frac
{w}{c};\alpha,\tau\right)e^{-\beta w/c} =\frac{\eta}{\alpha-\tau}\int_{\tau}^{\alpha}\frac{sc^{s}}{\Gamma
(1-s)}w^{-1-s}e^{-\beta w/c}\mathrm{d}s\label{eq:rhofull},\\
\psimix(t;\alpha,\tau,\beta,c,\eta) &  =\eta\left(  \psica%
(\beta+ct;\alpha,\tau)-\psica(\beta;\alpha,\tau)\right).
\label{eq:psifull}
\end{align}
The five parameters of the model have the following interpretation:
\begin{itemize}
\item $\alpha\in(0,1]$ is the index of variation;
\item $\tau\in[0,\alpha)$ is a power-law exponent that tunes the decay
of large weights, see \eqref{eq:asymptoticslargeweights};
\item $\beta$ is the exponential tilting parameter; if $\beta>0$, it tunes the
exponential decay of large weights, see \eqref{eq:asymptoticslargeweights};
\item $c$ is a scaling parameter: if $G\sim \CRM(\rhomix(\cdot;\alpha,\tau
,\beta,1,\eta),H)$, then 

$cG\sim \CRM(\rhomix(\cdot;\alpha,\tau,\beta,c,\eta),H)$;
\item $\eta$ is a rate  parameter: if $G_{1}\sim \CRM(\rhomix
(\cdot;\alpha,\tau,\beta,c,\eta_{1}),H)$ and

$G_{2}\sim \CRM(\rhomix(\cdot
;\alpha,\tau,\beta,c,\eta_{2}),H)$, then 
$G_{1}+G_{2}\sim \CRM(\rhomix
(\cdot;\alpha,\tau,\beta,c,\eta_{1}+\eta_{2}),H)$.
\end{itemize}


\section{Properties and Simulation}
\label{sec:properties}

In this section, we conduct a theoretical analysis of the properties of the mGG CRM. Specifically, we prove in \cref{thm:AsymptoticsLevyLaplace} that the parameter $\alpha\in(0,1]$ in this model is the index of variation, and the model may therefore exhibit rapid variation when $\alpha=1$. All these results will serve as essential tools for our analysis of the asymptotic properties of random graphs based on this model in \cref{sec:graphs}.

\subsection{Asymptotic properties of the L\'evy intensity and Laplace exponent} 

We derive here some asymptotic properties of the CRM. The behavior of the Laplace exponent, as explained in the introduction, tunes important asymptotic properties related to the behavior of small weights of the CRM. The behaviour of the L\'evy intensity at infinity is also of interest for some applications, as it relates to the behavior of large weights.

\begin{proposition}
\label{thm:AsymptoticsLevyLaplace}
We have, for the canonical model
\begin{align}
\psica(t;\alpha,\tau) &  \underset{t\rightarrow\infty}{\sim}\frac{t^{\alpha}}{(\alpha-\tau)\log t}\\
\rhoca(w;\alpha,\tau)&\underset{w\rightarrow0}{\sim}  \left \{\begin{array}{ll}
                                     \frac{1}{1-\tau}\frac{w^{-2}}{\log^2(1/w)} & \text{if }\alpha=1, \\
                                     \frac{\alpha}{(\alpha-\tau)\Gamma(1-\alpha)} \frac{w^{-1-\alpha}}{\log(1/w)}  & \text{if }\alpha<1,
                                   \end{array}\right .\\
\rhoca(w;\alpha,\tau)&\underset{w\rightarrow\infty}{\sim}  \left \{\begin{array}{ll}
                                     \frac{1}{\alpha-\tau}\frac{w^{-1-\tau}}{\log^2(w)} & \text{if }\tau=0, \\
                                     \frac{\tau}{(\alpha-\tau)\Gamma(1-\tau)} \frac{w^{-1-\tau}}{\log(w)}  & \text{if }\tau>0.
                                   \end{array}\right .
\end{align}
For the full model,%
\begin{align}
\psimix(t;\alpha,\tau,\beta,c,\eta) &  \underset{t\rightarrow\infty}{\sim} \eta c^\alpha \psica(t;\alpha,\tau),\\
\rhomix(w;\alpha,\tau,\beta,c,\eta)&\underset{w\rightarrow0}{\sim} \eta c^\alpha \rhoca(w;\alpha,\tau),\\
\rhomix(w;\alpha,\tau,\beta,c,\eta)&\underset{w\rightarrow\infty}{\sim} \eta c^\tau e^{-\beta w/c}  \rhoca(w;\alpha,\tau).\label{eq:asymptoticslargeweights}
\end{align}
\end{proposition}

\subsection{Moments} 

Derivatives of the Laplace exponent $\psi$ of a CRM are of interest to derive moments or other properties. For $n\geq 1$ and $t> 0$, let
\begin{equation}
\kappa(n,t)=(-1)^{n+1} \psi^{(n)}(t)=\int_0^\infty w^n e^{-wt}\rho(w)\mathrm{d}w.
\label{eq:kappafunction}
\end{equation}
If $G\sim\CRM(\rho, H)$, then $\bbE[G(\Theta)]=\kappa(1,0)$ and $\bbV[G(\Theta)]=\kappa(2,0)$. We provide below the expression of the $\kappa$ function for the mGG CRM, and analytic expressions for the expectation and variance of $G(\Theta)$ under this model.

\begin{proposition}
\label{thm:moment}
Let $m\geq1$ and $z\geq0$. Define $$\kappamix(m,z;\alpha,\tau,\beta,c,\eta)=\int
_{0}^{\infty}w^{m}e^{-zw}\rho(w;\alpha,\tau,\beta,c,\eta)\mathrm{d}w$$ and $\kappaca(m,z;\alpha
,\tau)=\int_{0}^{\infty}w^{m}e^{-zw}\rhoca(w;\alpha,\tau)\mathrm{d}w$. Then
\begin{align*}
\kappaca\left(  m,z;\alpha,\tau\right)   &  =\frac{z^{-m}}{\alpha-\tau}\int_{\tau}^{\alpha}sz^{s}\frac{\Gamma
(m-s)}{\Gamma(1-s)}\mathrm{d}s,\\
\kappamix(m,z;\alpha,\tau,\beta,c,\eta)&=\eta c^{m}\kappaca\left(
m,\beta+cz;\alpha,\tau\right).
\end{align*}
For $m=1$, $z>0$, we have
$$
\kappaca\left(  1,z;\alpha,\tau\right)=\left \{\begin{array}{cc}
                                        \frac{z^{-1}}{\alpha-\tau}\frac{z^{\tau}-z^{\alpha}+(\alpha z^{\alpha
}-\tau z^{\tau})\log z}{\left(  \log z\right)  ^{2}} & z\neq 1,\ \\
                                         \frac{\alpha+\tau}{2} & z=1.
                                       \end{array}\right . 
$$
If G is a CRM with parameters $\alpha,\tau,\beta>0,c,\eta$, then the total mass $G(\Theta)=\sum_{j\geq 1} W_j$ has expectation $\bbE [G(\Theta)]=\eta c \kappaca(1,\beta;\alpha,\tau)$ and variance
$\bbV[G(\Theta)]=\eta c^2 \kappaca(2,\beta;\alpha,\tau)$. In particular, if $\alpha=1$ and $\tau=0$, we obtain
\begin{align*}
\bbE [G(\Theta)]  &  =
\left \{\begin{array}{cc}
                                       \eta c\frac{1-\beta+\beta\log(\beta)}{\beta\log^{2}(\beta)}  & \beta\neq 1,\ \\
                                    \frac{\eta c}{2} & \beta=1,
                                       \end{array}\right.\\
\bbV[G(\Theta)] &  = \left \{\begin{array}{cc}
                                        \eta c^{2} \frac{(\beta+1)\log\beta+2(1-\beta)}{\beta^{2}%
\log^{3}\beta} & \beta\neq 1,\ \\
                                    \frac{\eta c^2}{6} & \beta=1.
                                       \end{array}\right.
\end{align*}
\end{proposition}

\subsection{Simulation of the weights and total mass}
\label{sampling}

Among the key properties expected from a tractable statistical model is the ability to sample from it. However, since CRMs are infinite objects that cannot be fully represented on a finite computer, only approximate sampling is possible. Several approaches for approximating CRM sampling exist, and here we propose a size-biased method, which is advantageous due to its ease of implementation. Additionally, we provide an estimate for the total mass lost in the approximation, which is proportional to $1/\log(n)$, where $n$ is the number of size-biased samples. In certain applications, only the total mass of the CRM is required. For such cases, we introduce an approximate method for sampling the total mass based on a Riemann sum, which converges at a rate of $1/n$. This convergence rate outperforms the alternative approach of computing individual weights and summing them.

\subsubsection{Size-biased sampling of the CRM} 
\label{sec:size-biased-sampling}

We follow here the strategy of \cite{Perman1992} to sample the size-biased weights, see \cref{app:sizebiased} in the supplementary material for details.

\begin{proposition}
\label{thm:sizebiasedmix}
Let $\xi_{1},\xi_{2},\ldots$ be the ordered points of a unit-rate Poisson
process on $(0,\infty)$; that is, $\xi_{1}$, $\xi_{2}-\xi_{1}$, $\xi_{3}%
-\xi_{2}$, \ldots\ are iid unit-rate exponential random variable. We have the following size-biased construction for $G\sim
\CRM(\rhomix(\cdot;\alpha,\tau,\beta,c,\eta),H)$. For $j\geq1,$%
\begin{align*}
T_{j}  &  =\psicainv\left(  \frac{\xi_{j}}{\eta}+\psica(\beta;\alpha,\tau)~;~\alpha,\tau\right)  -\beta,\\
S_j\mid \{T_j=t\} & \sim p(s|t)\propto s(t+\beta)^s \1{s\in(\tau,\alpha)},\\
W'_{j}|\{T_{j}=t,S_{j}=s\}&\sim\gammadist(1-s,t+\beta),\\
W_{j}  &  =cW_{j}^{\prime},%
\end{align*}
where $\psicainv$ is the inverse of $\psica$.
\end{proposition}

One can sample exactly from $p(s|t)$ using the inverse CDF, which has a tractable expression in terms of standard functions, see \cref{prop:inverseCDF} in the supplementary material. When $\tau=0$, the inverse $\psica^{-1}$ can be expressed using standard functions, see \cref{prop:inverselaplace}. The following proposition describes the asymptotic L1 error of the size-biased approximation. Its proof, similar to that in \cite{Lee2023}, is given in \cref{sup:omitedproof}.

 \begin{proposition}
 \label{thm:totalmasslost}
Let $\alpha=1$ and $\tau=0$. Let $(W_j)_{j \geq 1}$ be the sequence of weights sampled via the size-biased algorithm in \cref{thm:sizebiasedmix}. Then
 $$R_n=\sum_{j>n}^\infty W_j \underset{n\to\infty}{\sim}\frac{c\eta}{\log(n)} \quad\textrm{almost surely.}$$
 \end{proposition}

\subsubsection{Approximate sampling of the total mass}

One can use the mixture representation of the CRM in order to simulate random variables that are approximately distributed as $G(\Theta)=\sum_{j\geq 1}W_j$. The approximation relies on a weighted sum of exponentially tilted positive stable random variables, for which efficient simulation algorithms exist~\citep{Devroye2009,Hofert2011}.
\begin{proposition}
\label{thm:mathapprox}
Let $n\geq 1$. For $i=1,\ldots,n$ let $s_{i}=\tau + \frac{(\alpha-\tau)(i-1)}{n}$ and let $X_i$ be an exponentially tilted stable random variables, with stable parameter $s_i$ and exponential tilting parameter $\beta_i=\left(\frac{\eta}{n}\right)^{\frac{1}{s_i}}\beta$. That is, for any $t\geq 0$,
$$
\bbE [e^{-tX_i}]=\exp\left(-( (t+\beta_i)^{s_i}-\beta_i^{s_i})\right).
$$
Then
$$
S_n=c\sum_{i=1}^n \left(\frac{\eta}{n}\right)^{\frac{1}{s_i}} X_i \to G(\Theta)
$$
in distribution as $n\to\infty$.
\end{proposition}
\FloatBarrier

\section{Sparse network modelling} 
\label{sec:graphs}

We consider in this section the use of the mGG CRM on $\Theta=\bbR_+$, with $\alpha=1$, $\tau=0$, $\beta> 0$, $c>0$, $\eta>0$, within the Caron-Fox graph construction in order to obtain extremely sparse sequences of graph. We describe the model, its asymptotic properties, and a Markov chain Monte Carlo method for posterior inference of the model parameters.

\subsection{mGG graph model} 
\label{sec:CaronFox}

\cite{Caron2017} introduced a novel class of graph models based on CRMs which falls within the broader framework of graphex processes  -- a generalisation of dense graphons to the sparse setting~\citep{Caron2017,Veitch2015,Borgs2018}. For a CRM $$G=\sum_{j\geq 1}W_j\delta_{\theta_j}\sim  \CRM(\rhomix(\cdot;\alpha,\tau,\beta,c,\eta),\lambda),$$ where $\lambda$ is the Lebesgue measure, we consider the atomic measure on $\mathbb{R}_+^2$ 
\begin{equation}
\label{eq:undirectedgraph}
U=\sum_{i=1}^\infty\sum_{j=1}^\infty Z_{i,j}\delta_{(\theta_i,\theta_j)},
\end{equation}
where the binary random variables $Z_{ij}\in\{0,1\}$ indicate if there is an edge between node $i$ and node $j$. They are defined, for $i< j$, by
$$
Z_{ij}\mid (W_k)_{k\geq 1}\sim \text{Bernoulli}\left(1-e^{-2W_i W_j}\right),
$$
$Z_{ji}=Z_{ij}$, and $Z_{ii}\mid (W_k)_{k\geq 1}\sim\text{Bernoulli}\left(1-e^{-W_i^2}\right)$. The weight $W_i$  may be interpreted as a sociability parameter of node $i$; the larger $W_i$, the more likely node $i$ is to connect to other nodes. From $Z$ we derive a  growing family of graphs $(\mathcal G_t)_{t\geq 0}$, where  $\mathcal G_t=(\mathcal V_t,\mathcal E_t) $ is the graph of size $t$, with vertex set $\mathcal V_t$ and edge set $\mathcal E_t$, defined as
$$
\begin{aligned}
& \mathcal{V}_t=\left\{\theta_i \mid \theta_i \leq t \text { and } \exists \theta_k \leq t \text { s.t. } Z_{i k}=1\right\}, \\
& \mathcal{E}_t=\left\{\{\theta_i, \theta_j\} \mid t, \theta_j \leq t \text { and } Z_{i j}=1\right\}.
\end{aligned}
$$

\begin{figure}[t]
\begin{subfigure}{0.45\textwidth}
        \centering
\begin{tikzpicture}[node distance=1.4cm, auto,>=latex',scale=.8]
\draw[->,ultra thick] node[left] {0} (0,0) -- (6,0) node[right] {};
\draw[->,ultra thick] (0,0) -- (0,-6) node[right] {};
\draw[fill=blue,opacity=0.1] (0,0) -- (0,-6) -- (6,-6) -- (6,0)  ;
\draw[dashed] (1,0) -- (1,-6) node[right] {};
\draw[dashed] (3,0) -- (3,-6) node[right] {};
\draw[dashed] (3.5,0) -- (3.5,-6) node[right] {};
\draw[dashed] (5,0) -- (5,-6) node[right] {};
\draw[dashed] (0,-1) -- (6,-1) node[right] {};
\draw[dashed] (0,-3) -- (6,-3) node[right] {};
\draw[dashed] (0,-3.5) -- (6,-3.5) node[right] {};
\draw[dashed] (0,-5) -- (6,-5) node[right] {};
\node[draw,circle,inner sep=2.5pt,fill=blue,shading=ball] at (1,-1) {};
\node[draw,circle,inner sep=2.5pt,fill=blue,shading=ball] at (1,-5) {};
\node[draw,circle,inner sep=2.5pt,fill=blue,shading=ball] at (3,-1) {};
\node[draw,circle,inner sep=2.5pt,fill=blue,shading=ball] at (3.5,-1) {};
\node[draw,circle,inner sep=2.5pt,fill=blue,shading=ball] at (5,-1) {};
\node[draw,circle,inner sep=2.5pt,fill=blue,shading=ball] at (1,-3) {};
\node[draw,circle,inner sep=2.5pt,fill=blue,shading=ball] at (1,-3.5) {};
\node[draw,circle,inner sep=2.5pt,fill=blue,shading=ball] at (5,-3.5) {};
\node[draw,circle,inner sep=2.5pt,fill=blue,shading=ball] at (3.5,-5) {};
\node[] at (5.4,-3.8) {$\color{blue}Z_{ij}$};
\node[red] at (5,-.3) {$\color{red}\theta_{i}$};
\node[red] at (.3,-3.5) {$\color{red}\theta_{j}$};
\draw[-,ultra thick,red] (1,0) -- (1,1) node[right] {};
\draw[-,ultra thick,red] (3,0) -- (3,.3) node[right] {};
\draw[-,ultra thick,red] (3.5,0) -- (3.5,.35) node[right] {};
\draw[-,ultra thick,red] (5,0) -- (5,.5) node[above] {$\color{red}W_i$};
\draw[-,ultra thick,red] (.6,0) -- (.6,.2) node[right] {};
\draw[-,ultra thick,red] (1.2,0) -- (1.2,.1) node[right] {};
\draw[-,ultra thick,red] (4.1,0) -- (4.1,.2) node[above] {};
\draw[-,ultra thick,red] (0, -1) -- (-1,-1) node[right] {};
\draw[-,ultra thick,red] (0,-3) -- (-.3,-3) node[right] {};
\draw[-,ultra thick,red] (0,-3.5) -- (-.35,-3.5) node[left] {$\color{red}W_j$};
\draw[-,ultra thick,red] (0,-5) -- (-.5,-5) node[right] {};
\draw[-,ultra thick,red] (0,-.6) -- (-.2,-.6) node[right] {};
\draw[-,ultra thick,red] (0,-1.2) -- (-.1,-1.2) node[right] {};
\draw[-,ultra thick,red] (0,-4.1) -- (-.2,-4.1) node[right] {};
\end{tikzpicture}
                \caption{}

    \end{subfigure}
    \begin{subfigure}{0.45\textwidth}
        \centering
               \begin{tikzpicture}[-, thick, node distance=1.5cm]

\node[circle, draw, fill=blue!10] (1) at (0,4) {1};
\node[circle, draw, fill=blue!10] (2) at (2,2) {2};
\node[circle, draw, fill=blue!10] (3) at (0,0) {3};
\node[circle, draw, fill=blue!10] (4) at (-2,2) {4};

\draw (1) edge[loop above, every loop/.style={-}] (1);      
\draw (1) -- (2);                    
\draw (1) -- (3);                    
\draw (1) -- (4);                    
\draw (4) -- (3);                    

\end{tikzpicture}                \caption{}
    \end{subfigure}
   \caption{Point process representation of a random graph. Each node $i$ is
embedded in $\mathbb{R}_{+}$ at some location $\theta_{i}$ and is associated with a
sociability parameter $W_{i}$. An edge between nodes $\theta_{i}$ and
$\theta_{j}$ is represented by a point at locations $(\theta_{i},\theta_{j})$
and $(\theta_{j},\theta_{i})$ in $\mathbb{R}_{+}^{2}$ (Figure 1, \citealp{Caron2017}). The corresponding graph is plotted in (b).}
\label{fig:pointprocess}%
\end{figure}
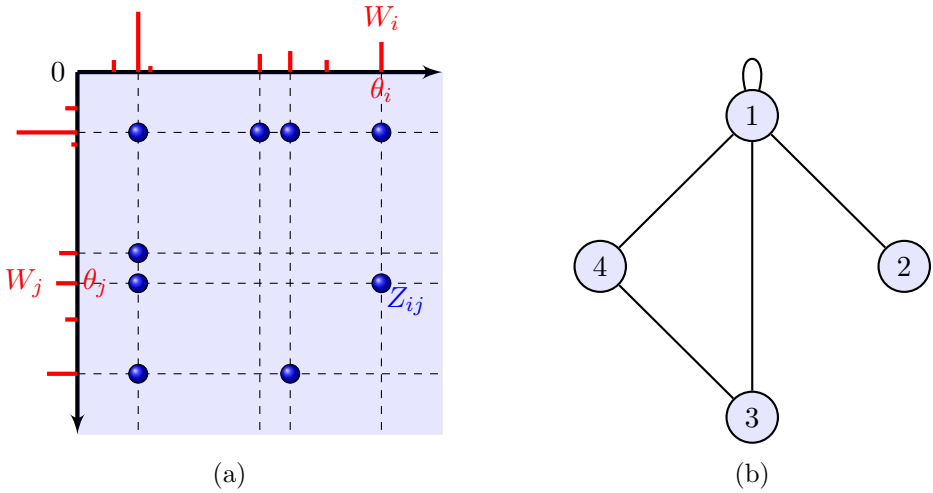

\subsection{Properties}
\label{sec:graphproperties}

For a graph $\mathcal{G}_t$, let $N_t=|\mathcal{V}_t|$ be the number of nodes, $N_t^{(e)}=|\mathcal{E}_t|$ the number of edges and $N_{t,j}$ the number of nodes of degree $j, j\geq1$. 

The following sparsity results follow from  \citet[Proposition 11]{Caron2023}. It states that the graphs are extremely sparse, with the number of edges scaling approximately linearly (up to a log factor) with respect to the number of nodes. \cref{fig:properties} shows a random graph drawn from the model compared to a (dense) Erd\H{o}s--R\'enyi graph.

\begin{proposition}[Asymptotic number of nodes and edges] We have almost surely that
\begin{align}
 N_t&\underset{t\rightarrow\infty}{\sim}  t^2 \frac{C}{\log(t)},~~~N_t^{(e)}\underset{t\rightarrow\infty}{\sim} \frac{t^2}{2}\overline{W}
\end{align}
 where $C=(\eta c)^2$ if $\beta=1$ and $C=2(\eta c)^2\frac{\beta^{-1}-1+\log \beta}{\left(  \log \beta\right)  ^{2}}$ otherwise, and $$\overline W=\int_0^\infty \psimix(2w;1,0,\beta,c,\eta) \rhomix(w;1,0, \beta,c,\eta)\mathrm{d}w.$$
\end{proposition}

\begin{corollary}[Extreme sparsity]
\label{thm:extremesparse}
We have almost surely:
$$N_t^{(e)}\underset{t\rightarrow\infty}{\sim} \frac{\overline{W}}{4C}N_t\log(N_t).$$
Hence $\frac{N_t^{(e)}}{N_t}\rightarrow \infty$ and for any $\epsilon>0$ we have  $\frac{N_t^{(e)}}{N^{1+\epsilon}_t}\rightarrow 0$. In other words the graph is almost extremely sparse.
\end{corollary}

\begin{figure}[!ht]
    \centering
    \begin{subfigure}{0.45\textwidth}
        \centering
        \includegraphics[width=\textwidth]{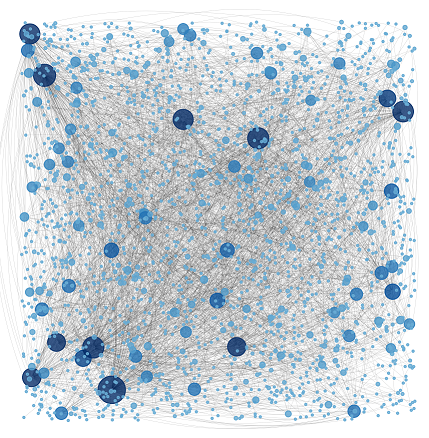}
        \caption{}
    \end{subfigure}
    \begin{subfigure}{0.45\textwidth}
        \centering
        \includegraphics[width=\textwidth]{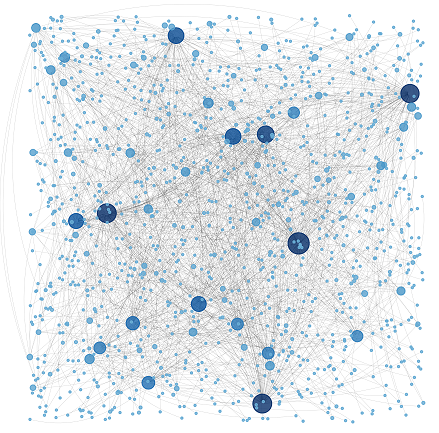}
        \caption{}
    \end{subfigure}
       \caption{Sampled mGG graphs (a)  with parameters $\alpha=1$, $\tau=0$, $\beta=1$, $c=1$, and $\eta=100$, 2122 nodes, and 2583 edges; (b)  with parameters $\alpha=1$, $\tau=0$, $\beta=1.5$, $c=5$, and $\eta=20$, 1249 nodes, and 1366 edges. The node size is proportional to its degree. The graphs were generated using our size-biased method and displayed using NetworkX.}
    \label{fig:graph}
\end{figure}
Another consequence of \cite[Proposition 11]{Caron2023} is the asymptotic degree distribution.
\begin{proposition}[Asymptotic degree distribution]
\label{thm:degdistrimix}
We have that, almost surely,
 $$N_{t,j}\underset{t\rightarrow\infty}{\sim} \left \{ \begin{array}{ll}
t^2 \frac{C}{\log(t)} & j=1,\\
\frac{t^2}{j(j-1)}\frac{C}{\log^2(t)} & j\geq 2.
\end{array}\right.$$
Consequently, the nodes of degree 1 dominate in the graphs as $N_{t,j}=o(N_{t,1})$ for all $j\geq 2$. Let $\tilde{N}_{t,2}=\sum_{k\geq2}N_{t,k}$ be the number of nodes of degree at least 2. Then, almost surely for all $j\geq 2$
\begin{equation}
\label{eq:powerlawdistribb}
\frac{N_{t,j}}{\tilde{N}_{t,2}}\underset{t\rightarrow\infty}{\longrightarrow} \frac{1}{j(j-1)}.
\end{equation}
This corresponds, for nodes of degree larger than 2, to a degree distribution with power-law of exponent 2.
\end{proposition}

In \cref{fig:properties} we compare our methods to the Generalized Gamma without mixture of \citet{Caron2017}  and with the Barabási–Albert model \citep{Barabasi1999a}. All these models are known to lead to sparse sequences of graphs presenting power law degree distribution. The expected asymptotic behaviour of our model is given in \cref{thm:extremesparse} and \eqref{eq:powerlawdistribb}. For the Generalized Gamma, here with $\alpha=0.5$, we have asymptotically $N_t^{(e)}\underset{t\rightarrow\infty}{\sim} N_t^{\frac{4}{3}}$ and a degree distribution following a power law of exponent $1.5$. For the Barabási–Albert model, the number of edges is linear in the number of nodes and presents an asymptotic degree distribution following a power law of exponent $3$.

\begin{figure}[ht]
    \centering
    \begin{subfigure}{0.45\textwidth}
        \centering
        \includegraphics[width=\textwidth]{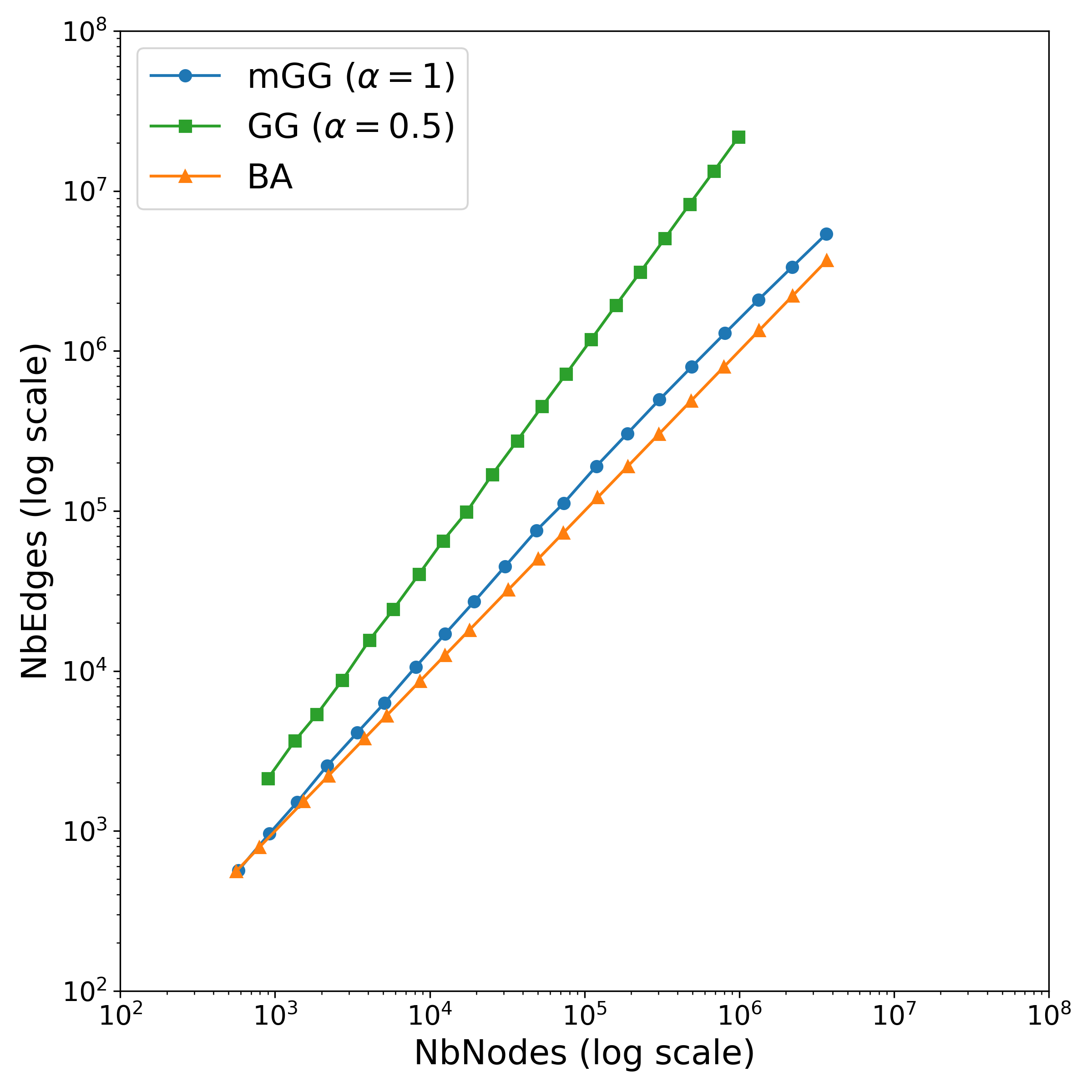}
    \end{subfigure}
    \begin{subfigure}{0.45\textwidth}
        \centering
        \includegraphics[width=\textwidth]{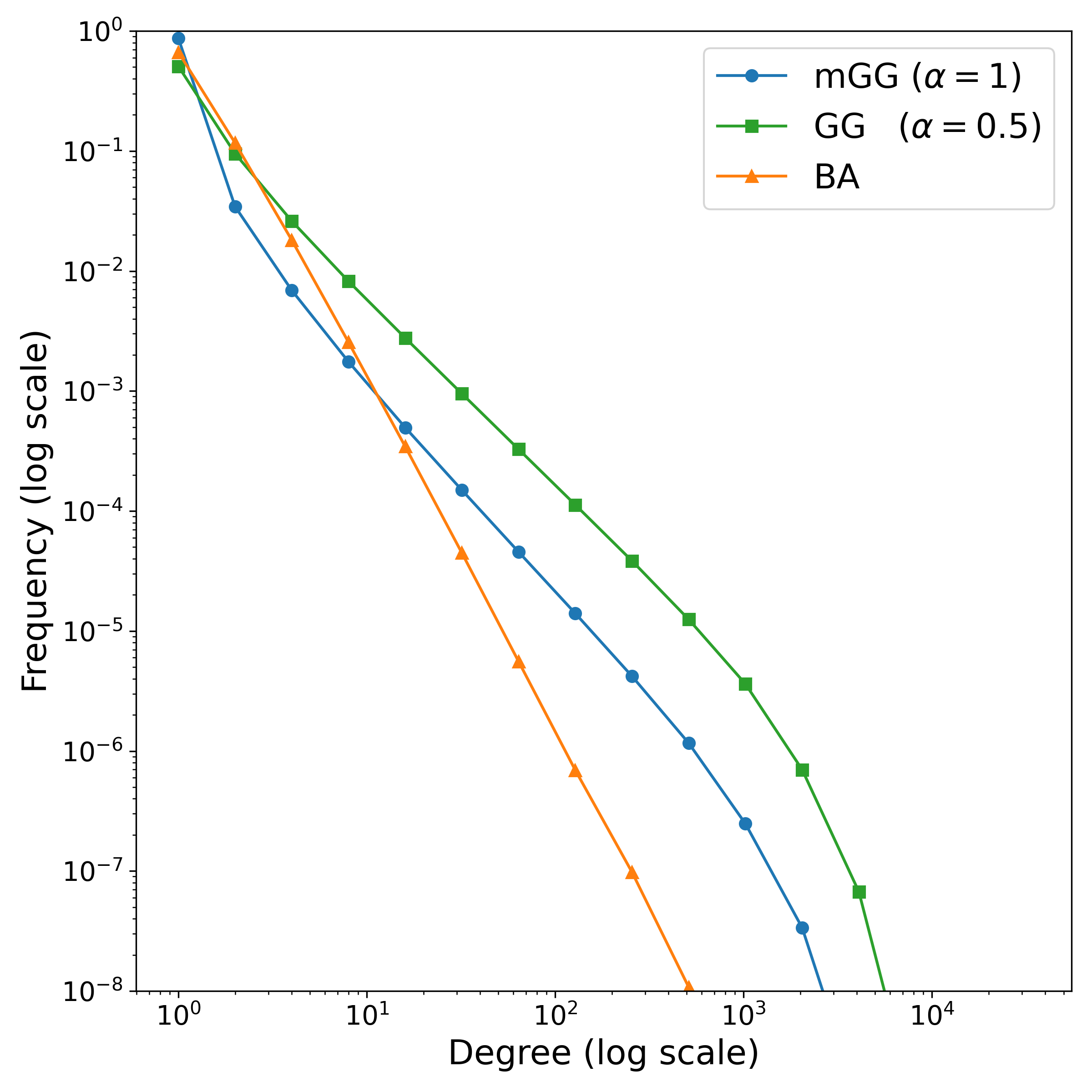}
    \end{subfigure}
   \caption{Examination of the mGG graph properties ({\Large$\bullet$}) with parameters $\alpha=1$, $\tau=0$, $c=1$, and $\beta=1$ for various values of $\eta$ ranging from 50 to 6000, resulting in graphs of different sizes. Comparison with the Generalized Gamma CRM ({\footnotesize $\blacksquare$}) with parameters $\tau=1$ and $\sigma=0.5$, and with the Barabási–Albert model ($\blacktriangle$). For every configuration we simulate 20 graph samples and plot the mean of the quantity of interest.}
   \label{fig:properties}
\end{figure}

\subsection{Posterior inference}
\label{sec:graphinference}

Here we briefly describe the MCMC algorithm to approximate the posterior density over the model parameters. Assume that we have observed a graph $\mathcal{G}$, which we assume was sampled from a mGG graph model. We aim to infer the sociability parameter $W_i$ of each nodes and the parameters $\phi=(\beta,c,\eta)$ of the mGG model. We consider the following improper priors:
\begin{align}
p(\beta) &\propto \frac{1}{\beta},~~p(c)\propto \frac{1}{c},~~p(\eta) \propto \frac{1}{\eta}.\label{eq:priorshyper}
\end{align}
We introduced two sets of auxiliary variables. First, as in \cite{Caron2017}, we introduce latent count variables $\tilde{q}_{i j}$ (an interpretation of this $\tilde{q}_{i j}$ in term of directed multigraph is given in \cref{sup:multigraph}) , with conditional distribution, for $i<j$,
\begin{equation}
\label{latentcountdistrib}
\tilde{q}_{i j} \mid Z, W \sim \begin{cases}\delta_0 & \text { if } Z_{i j}=0, \\ \operatorname{tPoisson}\left(2 W_i W_j\right) & \text { if } Z_{i j}=1, i \neq j, \\  \operatorname{tPoisson}\left( W_i^2\right) & \text { if } Z_{i j}=1, i = j,
\end{cases}
\end{equation}
where $\operatorname{tPoisson}(\lambda)$ is the zero-truncated Poisson distribution, and with $\tilde{q}_{ji}= \tilde{q}_{ij}$. Second, to take advantage of the mixture representation of the L\'evy intensity, we introduce, as in the size-biased algorithm, local indices of variations $S_i\in(0,1)$ for each node $i$. In the following, we also write $w_*$ for the sum of the sociabilities of nodes with no connection in $\mathcal G$. The algorithm is presented in \cref{alg:cap}.

\begin{algorithm}[t]
\caption{Posterior inference}
\label{alg:cap}
\textbf{Step 1:} update the weights $w_i$ and the latent space variables $s_i$ given the rest, using a Hamiltonian Monte Carlo step.

\textbf{Step 2:} update $(\phi,w_*)$ given the rest, using a Gibbs sampler update.

\textbf{Step 3:} update the latent count $(\tilde{q}_{ij})$ given the rest using \eqref{latentcountdistrib}.
\end{algorithm}
\FloatBarrier

\section{Numerical experiments on networks} 
\label{sec:experiments}

\subsection{Synthetic data}
\label{art:synthetic}

We first examine the convergence of the MCMC algorithm on simulated data, where the graph is generated from the mGG graph model described in \cref{sec:CaronFox}. We simulate a mGG graph with parameters $\alpha=1$,  $\tau=0$, $\beta=1$, $c=2$ and $\eta=130$ using size-biased sampling with $N=10^5$ weights.  The resulting graph contains 11,613 nodes and 17,427 edges. We run three MCMC chains, each consisting of 6 million iterations and initialized from different starting values. Half of the iterations are discarded as burn-in. The hyperparameters $\alpha$ and $\tau$ are fixed to their true values, while the remaining parameters are estimated. Details of the MCMC parameterization are provided in \cref{sup:inference}.

\begin{figure}[t]
    \centering
    \begin{subfigure}{0.45\textwidth}
        \centering
        \includegraphics[width=\textwidth]{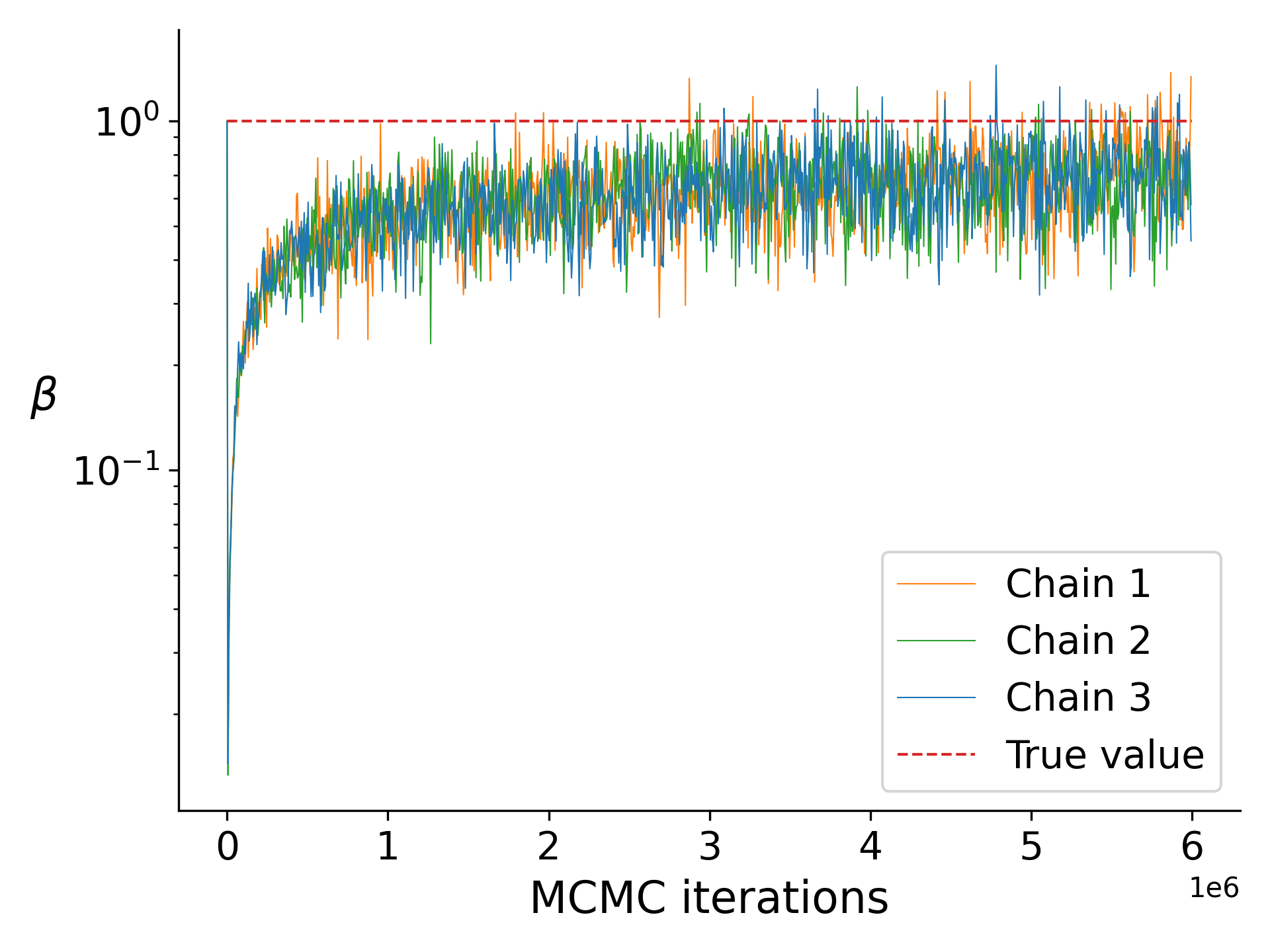}
        \caption{}
    \end{subfigure}
    \begin{subfigure}{0.45\textwidth}
        \centering
        \includegraphics[width=\textwidth]{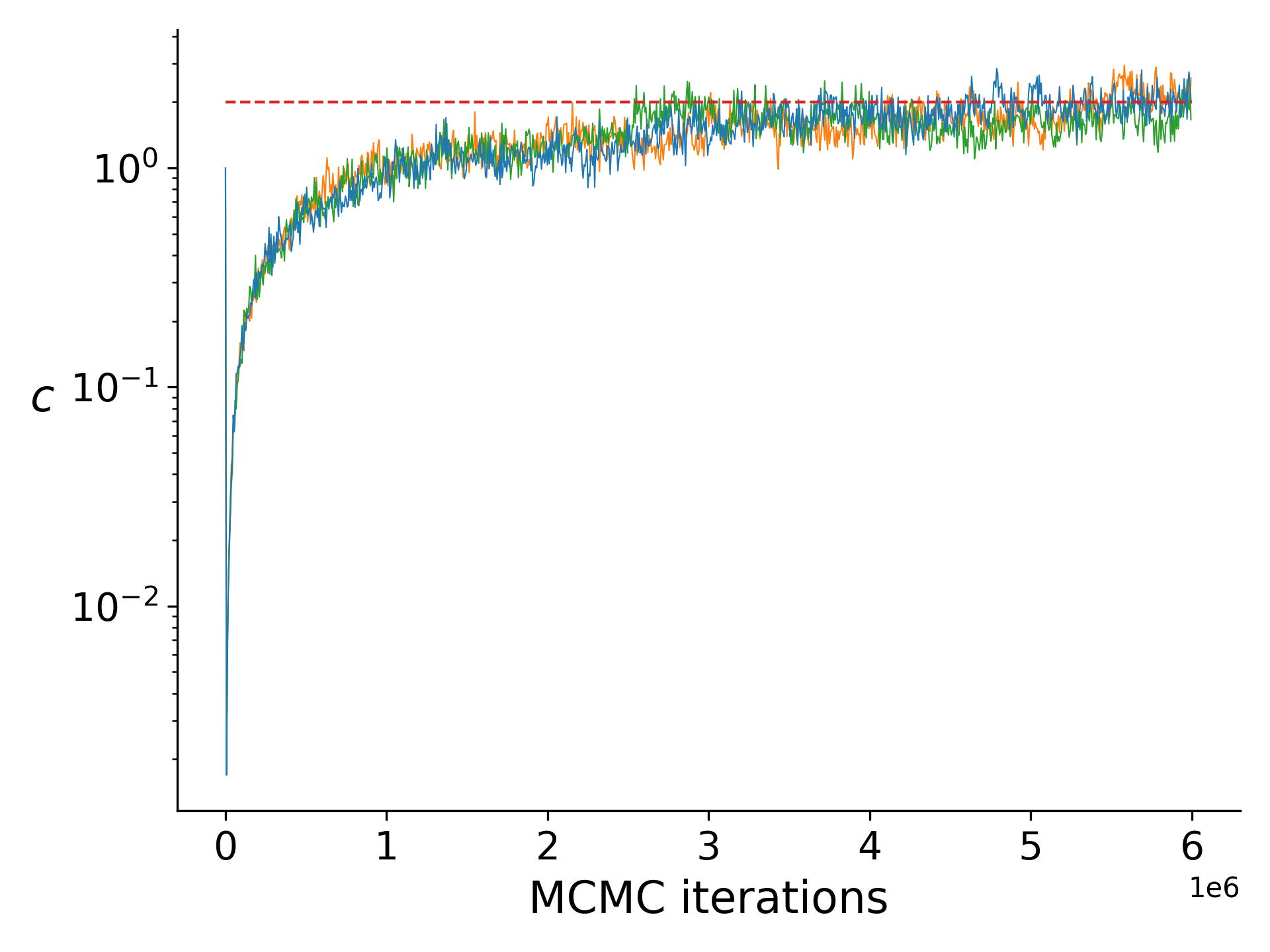}
        \caption{}
    \end{subfigure}
    \begin{subfigure}{0.45\textwidth}
        \centering
        \includegraphics[width=\textwidth]{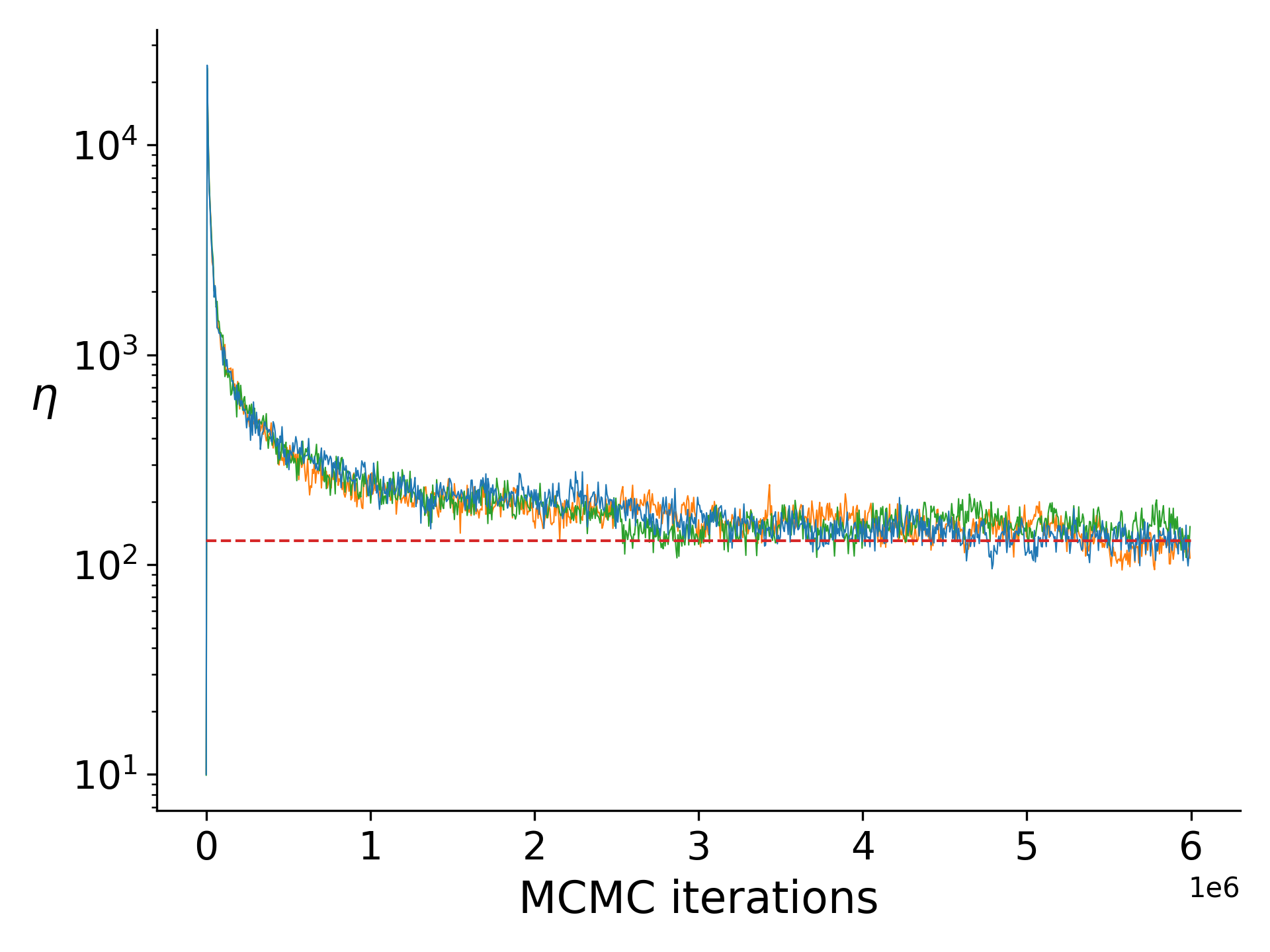}
        \caption{}
    \end{subfigure}
        \begin{subfigure}{0.45\textwidth}
        \centering
        \includegraphics[width=\textwidth]{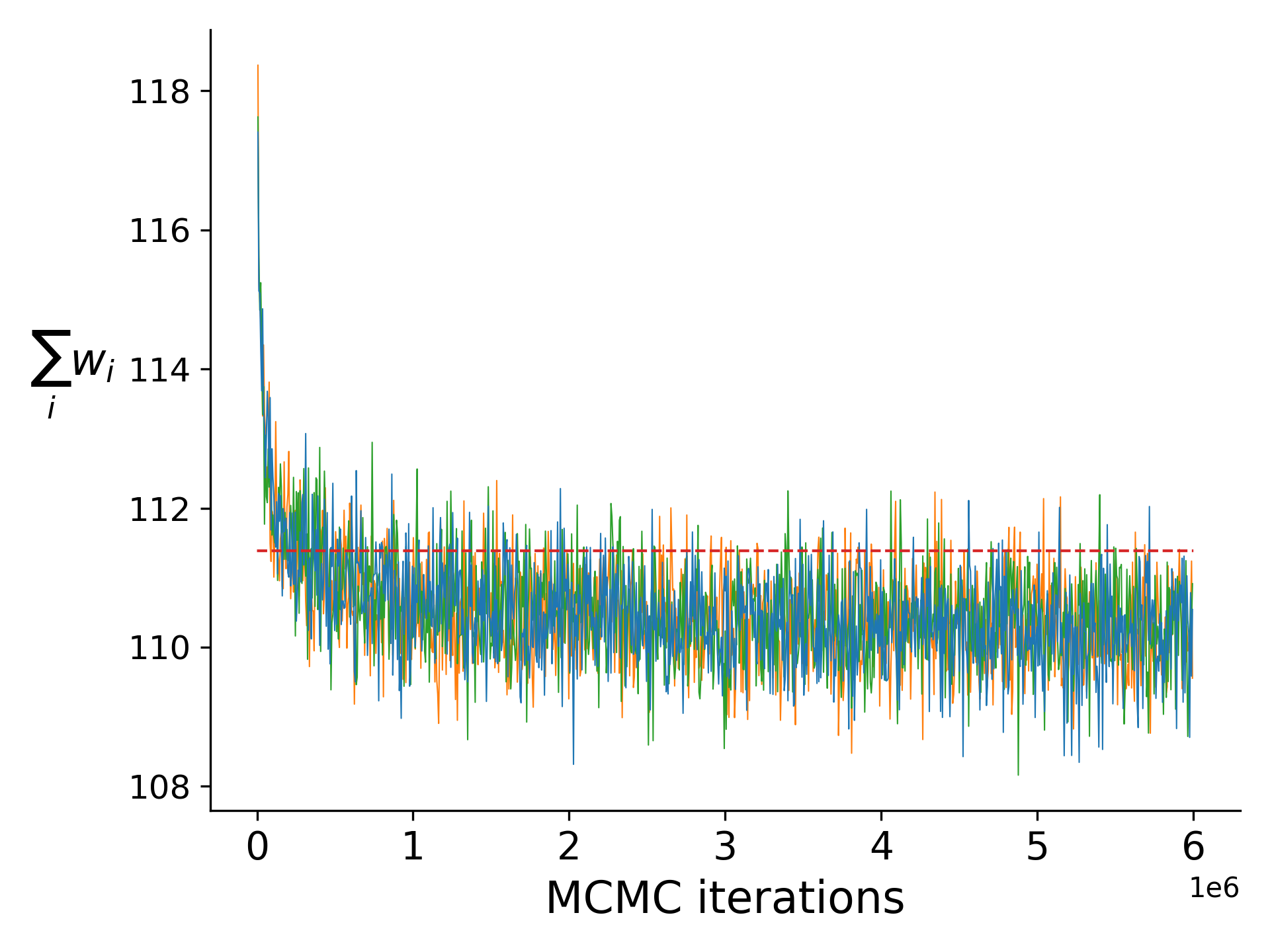}
        \caption{}
    \end{subfigure}
       \caption{MCMC traceplot of parameters (a) $\beta$, (b), c, (c), $\eta$ and (d) the total sum of the weights $w_i$ for a mGG graph with model parameters $\alpha=1, \tau=0, \beta=1, c=2$ and $\eta=130$. Three chains are displayed and the true value is in red dashed line.}
    \label{fig:traceplothyperparam}
\end{figure}

\begin{figure}[t]
    \centering
    \begin{subfigure}{0.45\textwidth}
        \centering
        \includegraphics[width=\textwidth]{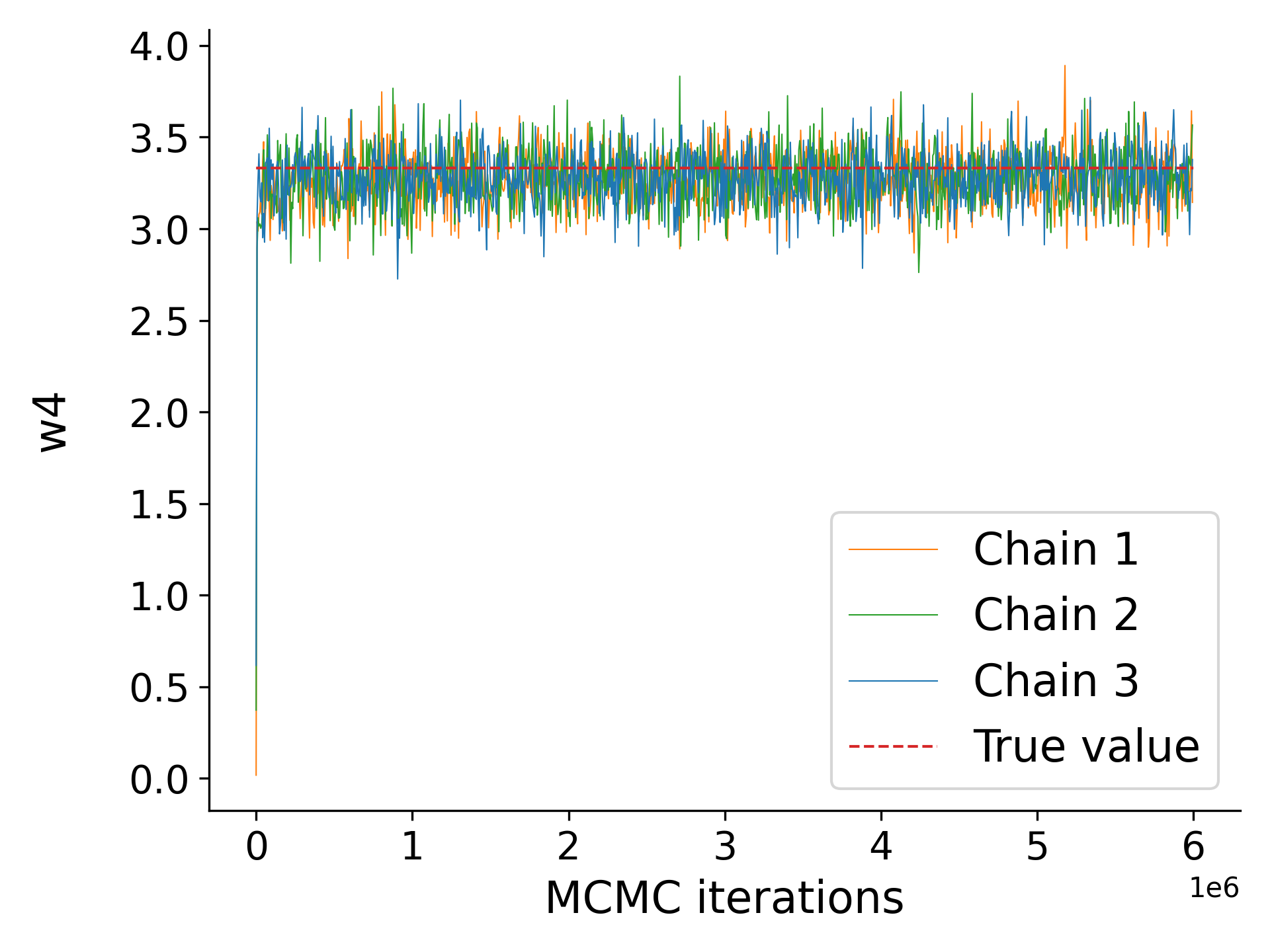}
        \caption{}
    \end{subfigure}
     \begin{subfigure}{0.45\textwidth}
        \centering
        \includegraphics[width=\textwidth]{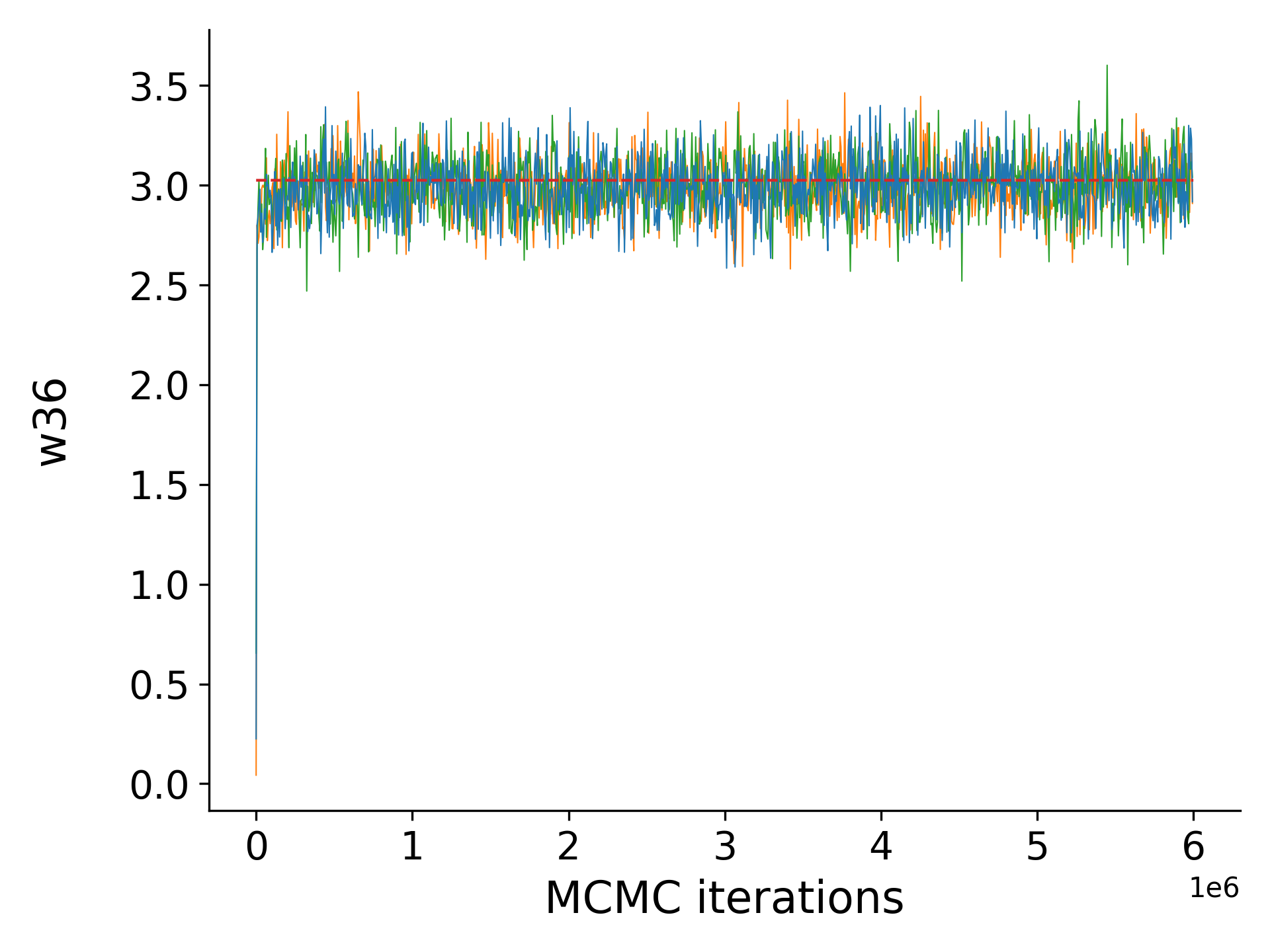}
        \caption{}
    \end{subfigure}

    \begin{subfigure}{0.45\textwidth}
        \centering
        \includegraphics[width=\textwidth]{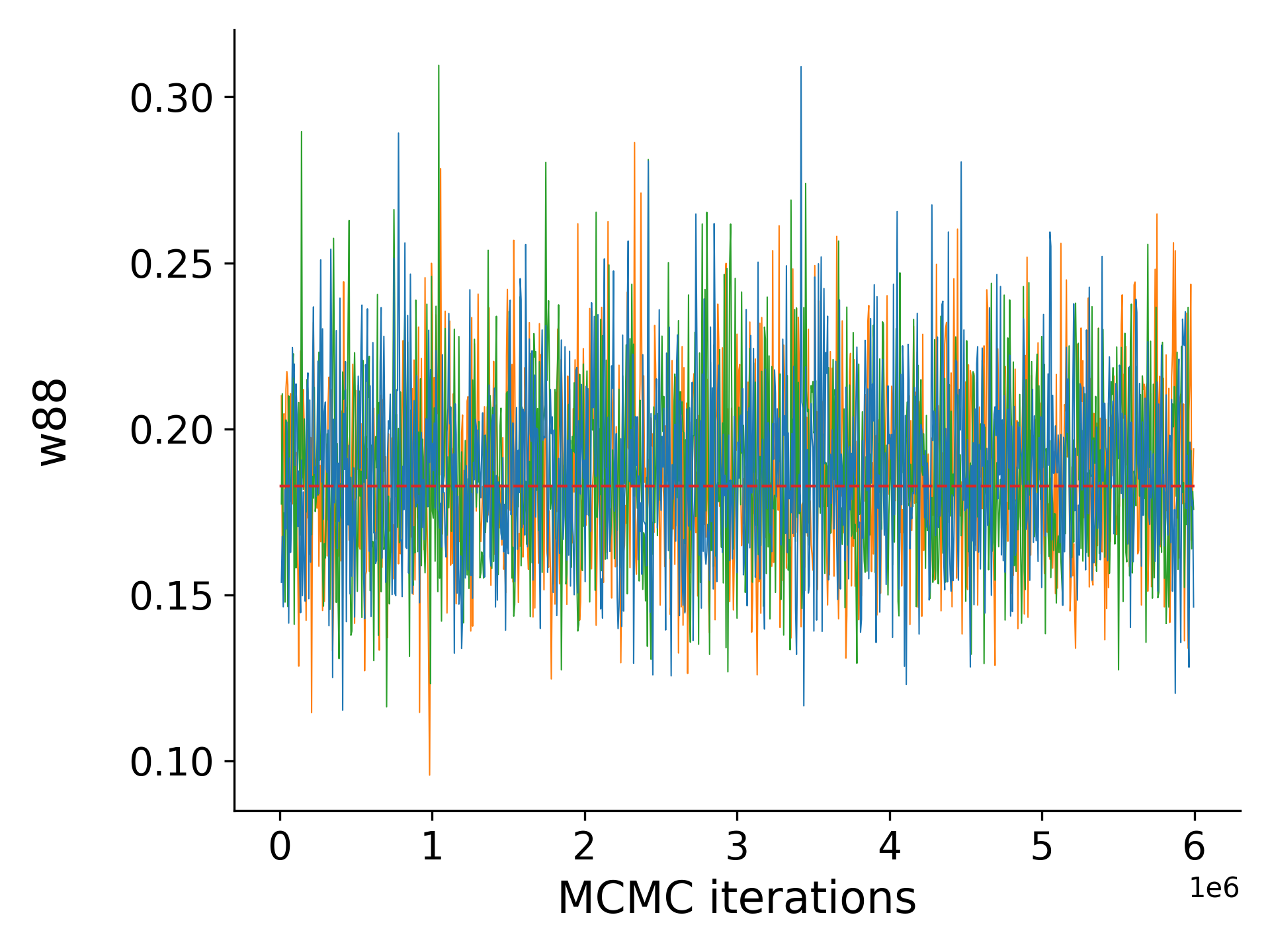}
        \caption{}
    \end{subfigure}
    \begin{subfigure}{0.45\textwidth}
        \centering
        \includegraphics[width=\textwidth]{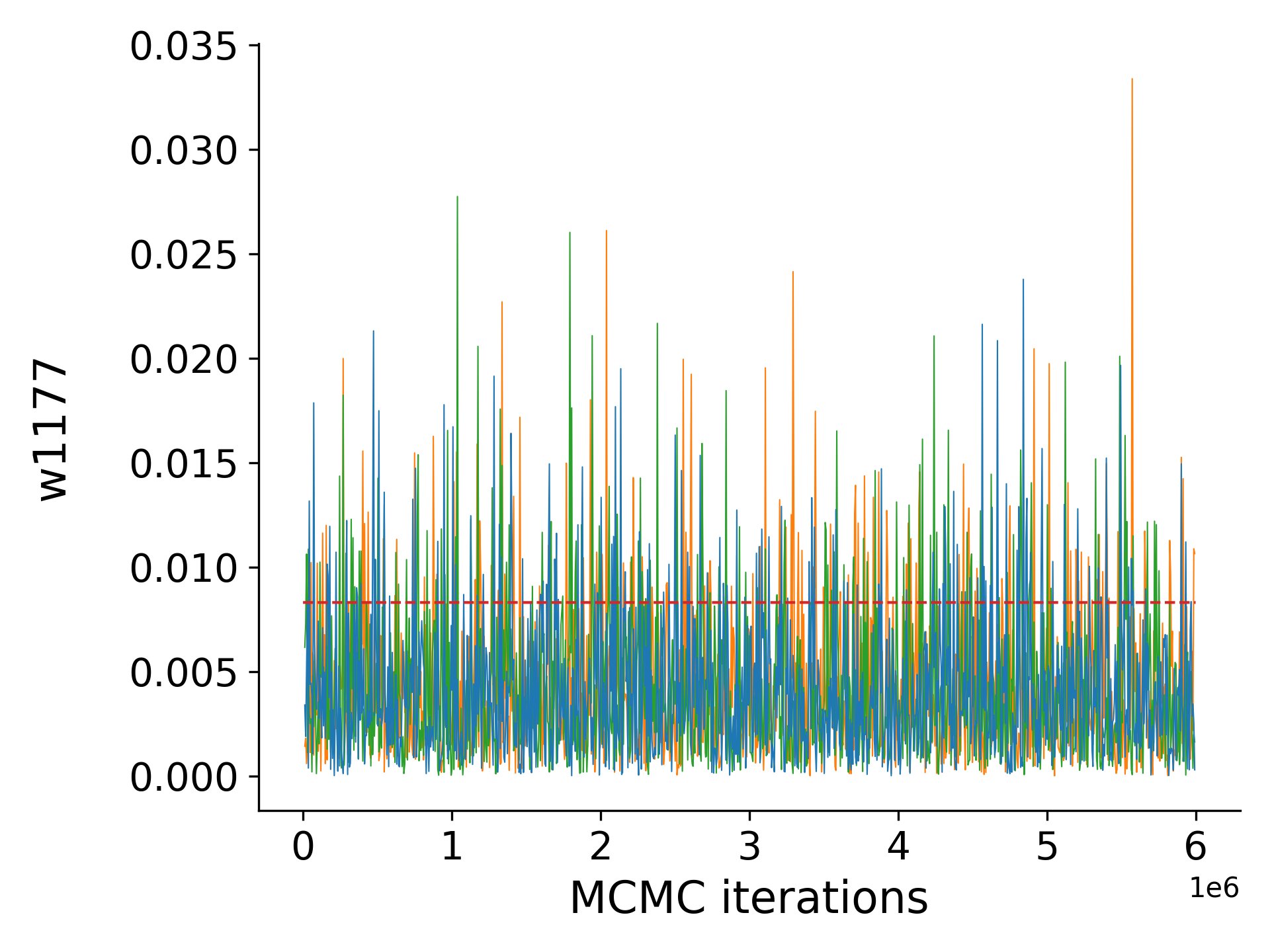}
        \caption{}
    \end{subfigure}
       \caption{MCMC traceplot of four weights for a graph generated from a mGG graph model with parameters $\alpha=1$, $\tau=0$, $\beta=1$, $c=2$ and $\eta=130$. Three chains are displayed and the true value is in red dashed line. The degree of the corresponding node is (a) 536, (b) 496, (c) 49 and (d) 2.}
    \label{fig:traceploweights}
\end{figure}

\begin{figure}[t]
    \centering
    \begin{subfigure}{0.45\textwidth}
        \centering
        \includegraphics[width=\textwidth]{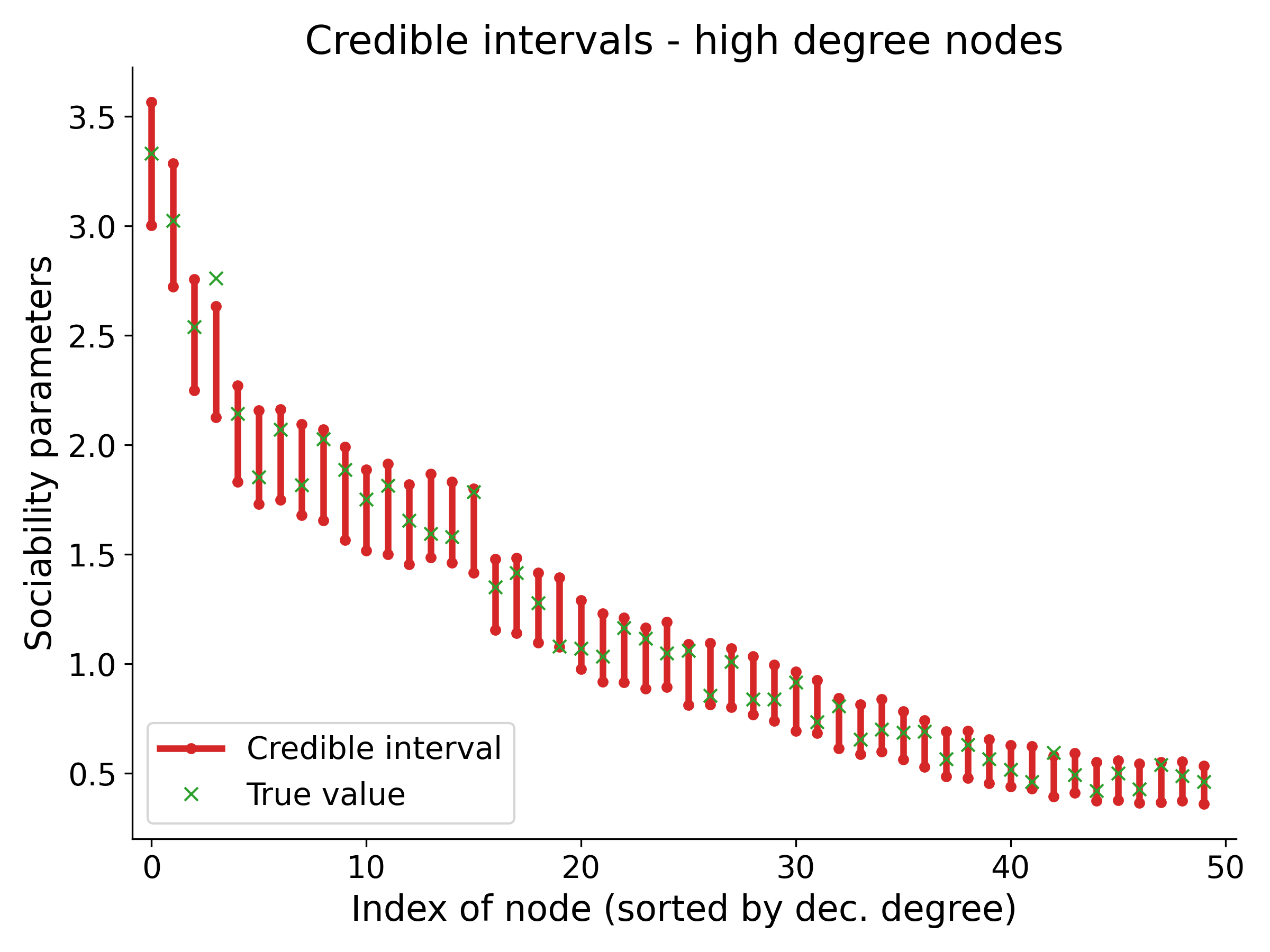}
        \caption{}
    \end{subfigure}
    \begin{subfigure}{0.45\textwidth}
        \centering
        \includegraphics[width=\textwidth]{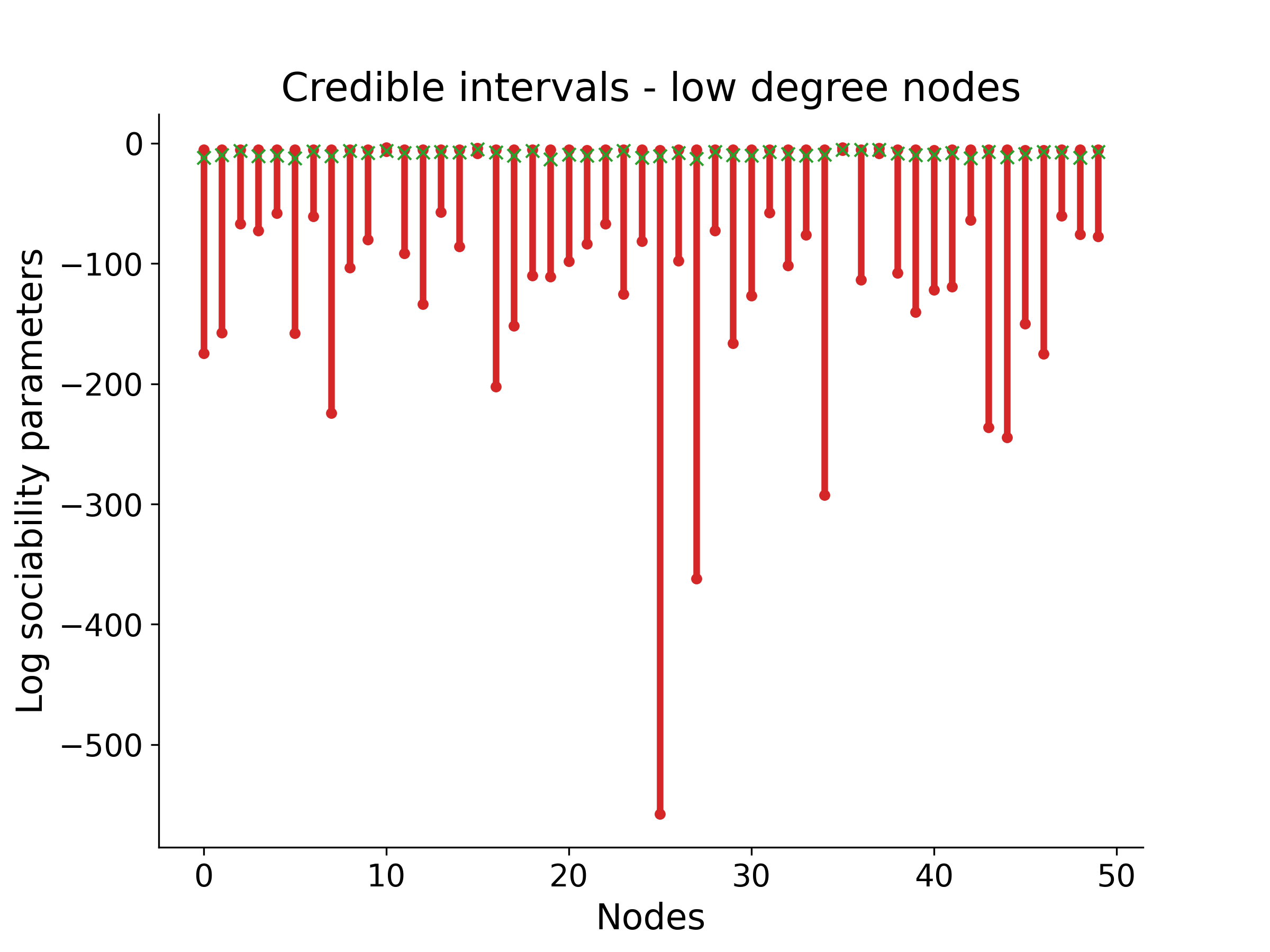}
        \caption{}
    \end{subfigure}
        \caption{95\% posterior intervals of (a) the sociability parameters $w_i$ of the 50 nodes with highest degree and (b) the log-sociability parameters  $\log(w_i)$ of the 50 nodes with lowest degree, for a sample from a mGG graph model with parameters $\alpha=1, \tau=0, \beta=1, c=2$ and $\eta=130$. The true values are in green.}
            \label{fig:credibleinterval}
\end{figure}

Trace plots for the parameters $\beta$, $c$, and $\eta$ are shown in \cref{fig:traceplothyperparam}, while trace plots for four of the 11,613 weights are presented in \cref{fig:traceploweights}. To further illustrate convergence of the weights (or sociability parameters) $w_i$, we also include in \cref{fig:traceplothyperparam}(d) the trace plot of the total sum of the weights. In addition, \cref{fig:credibleinterval} displays posterior credible intervals for the sociability parameters $w_i$ of the 50 nodes with the highest degrees, as well as for the log-sociability parameters $\log(w_i)$ of the 50 nodes with the lowest degrees. These results highlight the ability of our method to accurately recover sociability parameters for both high- and low-degree nodes.

\begin{figure}[t]
           \centering
    \begin{subfigure}{0.45\textwidth}
        \centering
        \includegraphics[width=\textwidth]{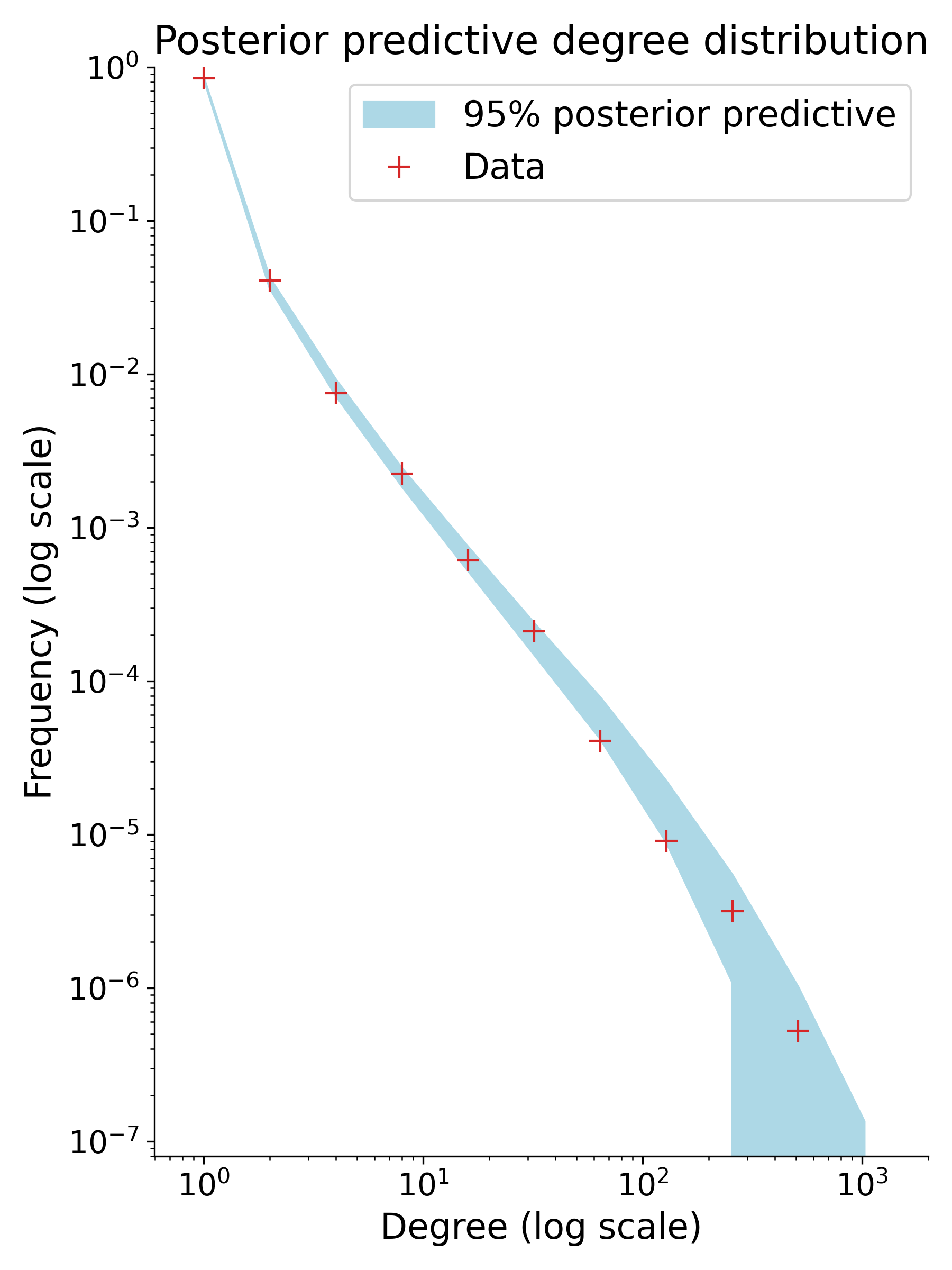}
        \caption{}
    \end{subfigure}
    \begin{subfigure}{0.45\textwidth}
        \centering
        \includegraphics[width=\textwidth]{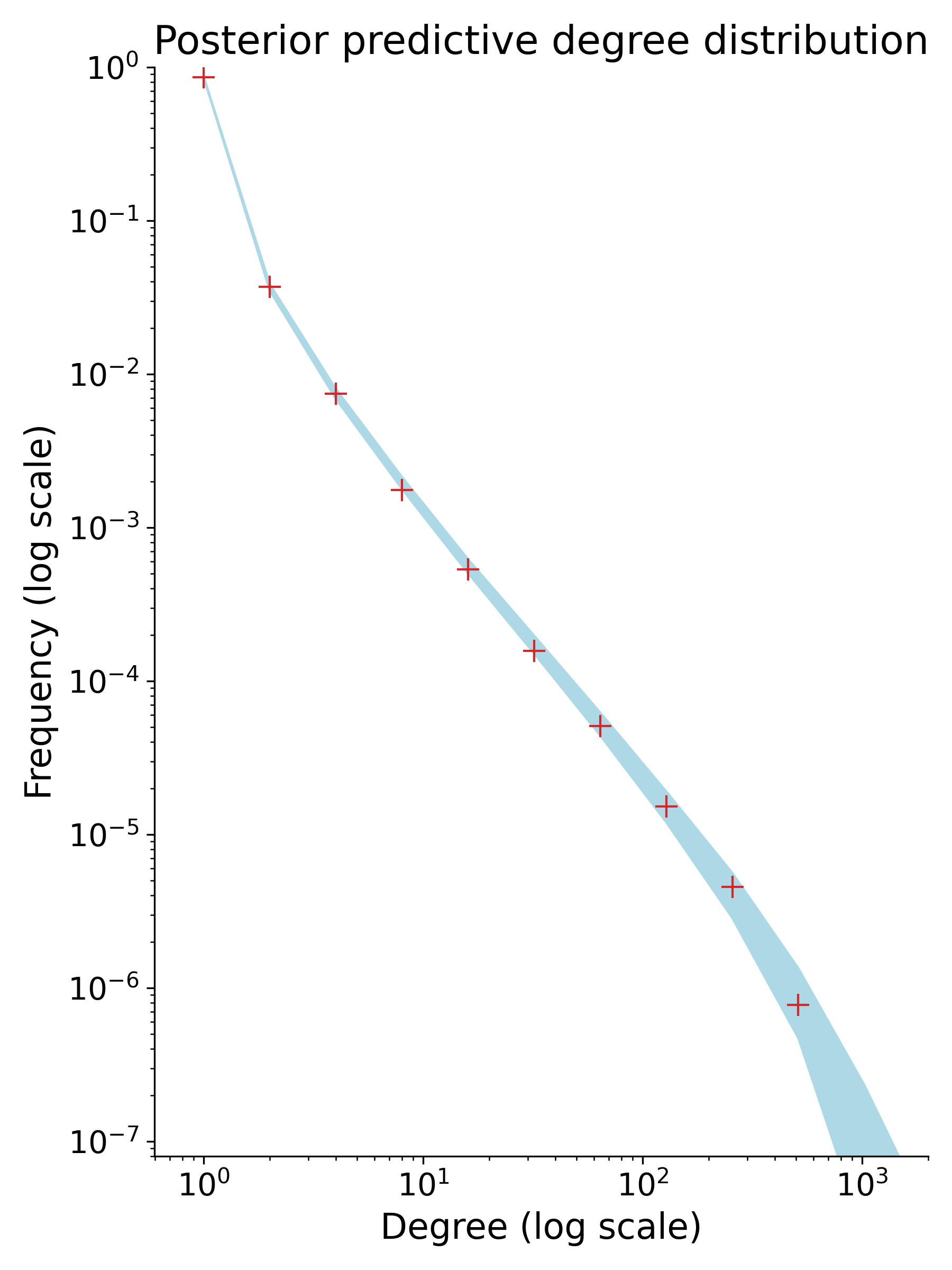}
        \caption{}
    \end{subfigure}
      \caption{(a) Posterior predictive degree distributions from our MCMC run, compared to those of a graph—distinct from the training graph—generated using the true parameters. (b) Same comparison, but with a reference graph that is twice as large; the corresponding $\hat\eta$ values from the MCMC run have been rescaled accordingly.}
            \label{fig:postdegree}
            \end{figure}

We use our sampling method to simulate 200 graphs from the posterior predictive distribution. \cref{fig:postdegree}(a) compares the degree distribution of a graph generated using the true parameter values with those of the simulated graphs. We then repeat the procedure with $\eta$ rescaled by a factor of 2: we generate a new reference graph using the true parameters but with $\eta' = 2\eta$, and simulate 200 posterior predictive graphs using the same MCMC run, but with the rescaled $\hat\eta'$. The resulting degree distributions are again compared in \cref{fig:postdegree}(b).

We compute the Gelman–Rubin convergence diagnostic \citep{Gelman1992} for the four key parameters, all of which yield values below 1.05 (see \cref{tab:GR} for exact figures). To assess convergence across all parameters -- $\beta$, $c$, $\eta$, $(w_i)$, and $(s_i)$ -- we use the multivariate Gelman–Rubin diagnostic introduced by \citet{Vats2021}, see \eqref{equa:Rmulti}, obtaining a value of 1.001, which strongly suggests convergence.

\subsection{Real World Data}
\label{art:real}
We now evaluate our method on real-world networks. We consider three datasets, all consisting of undirected graphs, downloaded from \citet{konect} and \citet{nr}.

\begin{itemize}
\item \textbf{Flickr:} A crawl of Flickr, a photo and video sharing platform. The dataset contains all links between users.  (\url{https://socialnetworks.mpi-sws.org/data-imc2007.html}, \citealp{Mislove2007a}).
\item \textbf{Douban:} A crawl of Douban, an online social network that provides user reviews and recommendations for movies, books, and music. The dataset contains all links between users (\url{https://networkrepository.com/soc-douban.php}, \citealp{Zafarani2014}).
\item \textbf{TwitterCrawl:} A crawl of Twitter, where nodes represent users and edges correspond to retweets collected from various social and political hashtags (\url{https://networkrepository.com/rt-retweet-crawl.php}, \citealp{Rossi2014}).

\end{itemize}

The sizes of the datasets are shown in \cref{table:size}. To train our model we extracted a subgraph containing approximately 5\% of the nodes from each network using $p$-sampling~\citep{Veitch2019}. The corresponding values of $p$ and subgraph statistics are reported in \cref{tab:sizepsample}. Further details and additional plots are provided in \cref{sup:realworld}.

\begin{table}[htp]
\caption{Sizes of the real-world datasets}
\begin{center}
\begin{tabular}{lcccc}
\hline Dataset & Nodes & Edges & Max degree & Mean degree \\
\hline Flickr & 1,861,232 & 155,55,041 & 54,472 & 33.4 \\
Douban & 154,908 & 327,162 & 574 & 8.4  \\
TwitterCrawl & 1,112,702 & 2,278,852 & 10,140 & 8.2  \\
\hline
\end{tabular}
\end{center}
\label{table:size}
\end{table}%

\begin{table}[ht]
\caption{Sizes of the subgraphs used for inference}
\begin{center}
\begin{tabular}{lccccc}
\hline Dataset & Nodes
 & 
Edges
 & Max degree &  Mean degree & p \\

\hline Flickr & 85,613 & 299,358 & 2890 & 13.9 & 0.14\\
Douban & 7775 & 9586 & 82 & 4.9 & 0.17 \\
TwitterCrawl & 70,944 & 87,914 & 618 & 5.0 & 0.2 \\
\hline
\end{tabular}
\end{center}
\label{tab:sizepsample}
\end{table}%

To assess model fit, we apply $p$-sampling to split each dataset into a training set (approximately 5\% of the nodes) and a test set (the remaining nodes). We then compare the empirical degree distributions of the test sets with the predictive degree distributions generated by our model, using hyperparameters inferred from the training data. Specifically, we follow the procedure described in \cref{sampling} to simulate 200 graphs from the posterior predictive distribution for each dataset. Parameters are drawn from the post-burn-in MCMC samples, and we scale $\hat\eta$ by $(1-p)/p$ to match the expected size of the target graph.

\begin{figure}[ht]
    \centering
    \begin{subfigure}{0.31\textwidth}
        \centering
        \includegraphics[width=\textwidth]{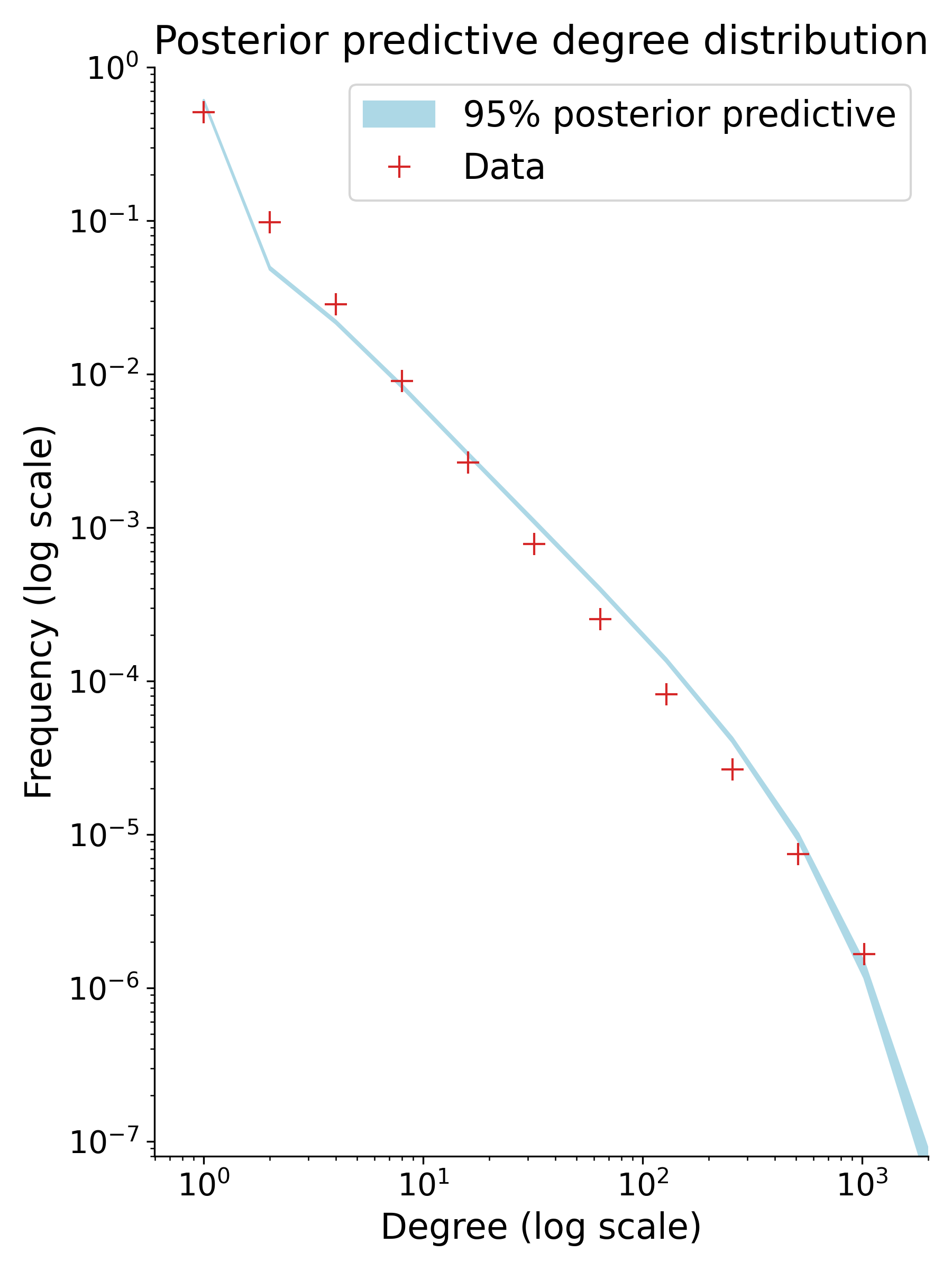}
        \caption{Flickr}
    \end{subfigure}
    \begin{subfigure}{0.31\textwidth}
        \centering
        \includegraphics[width=\textwidth]{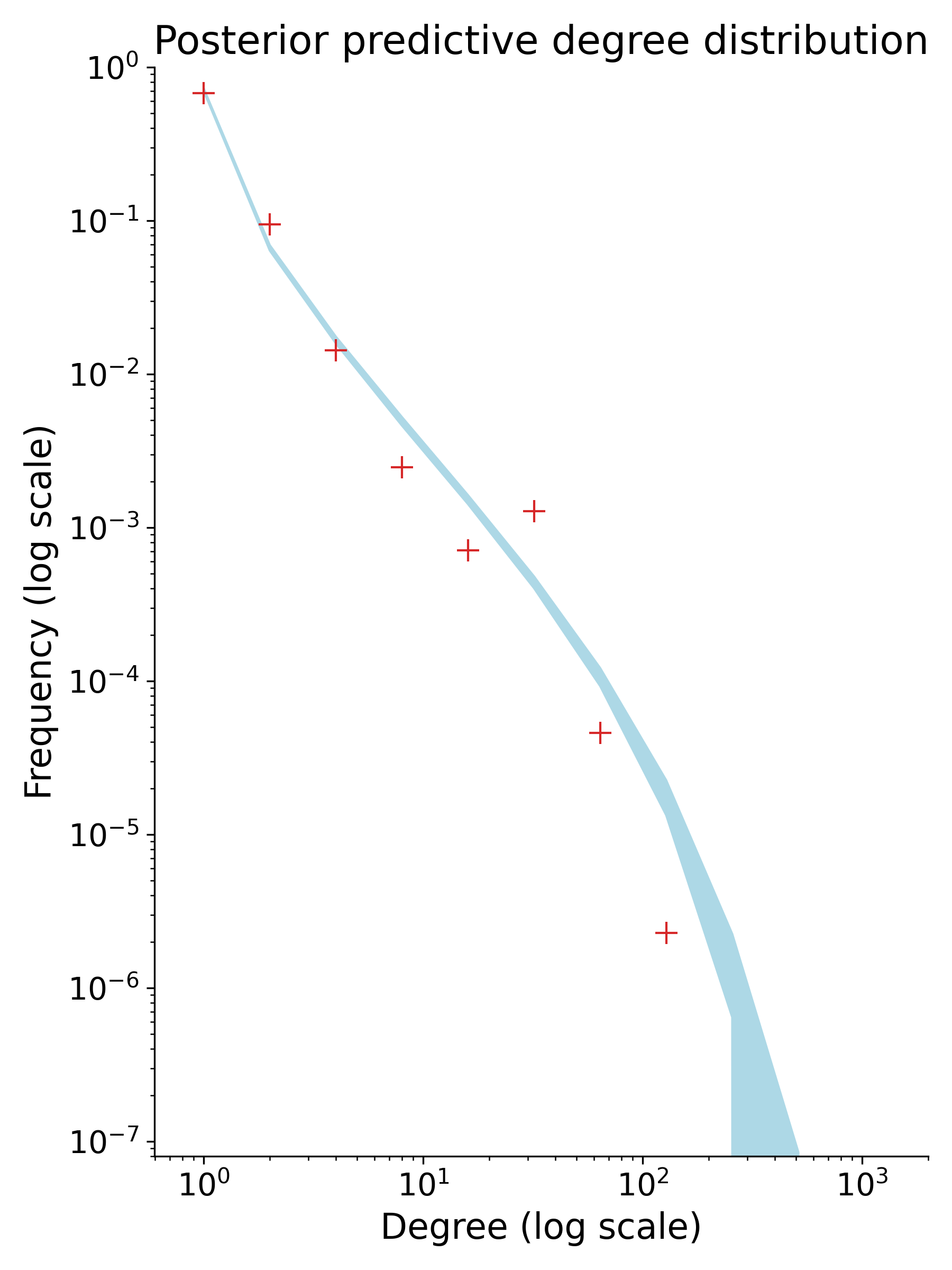}
        \caption{Douban}
    \end{subfigure}
    \begin{subfigure}{0.31\textwidth}
        \centering
        \includegraphics[width=\textwidth]{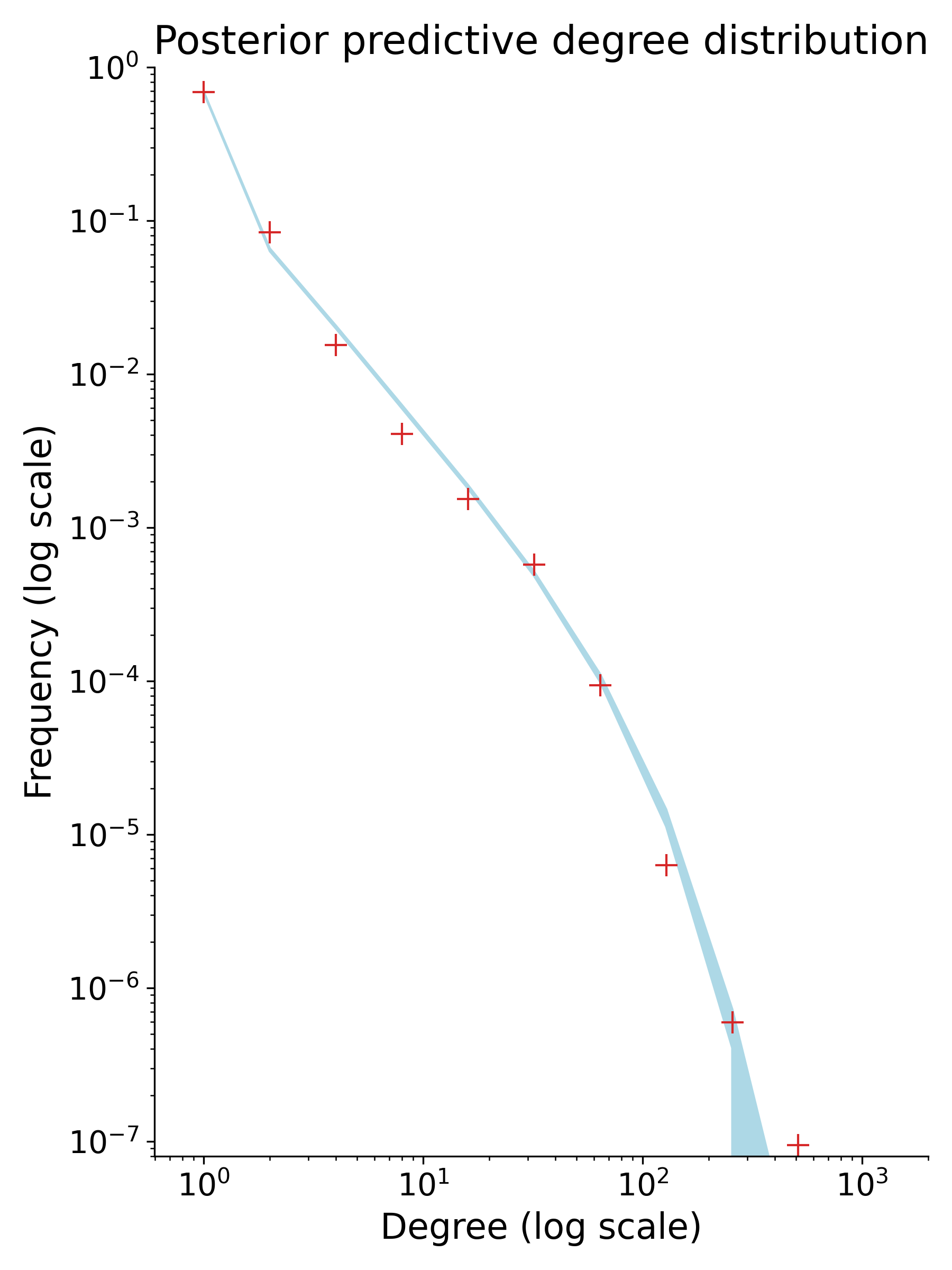}
        \caption{TwitterCrawl}
    \end{subfigure}
       \caption{Posterior predictive degree distributions for the real-world datasets. The predictive distribution is shown in blue, and the empirical distribution of the test set is shown in red.}
    \label{fig:postdegreepredictive}
\end{figure}

\Cref{fig:postdegreepredictive} shows the resulting comparisons between the empirical and predictive degree distributions. Overall, the model provides a good fit across the datasets, successfully capturing the heavy-tailed behavior characteristic of real-world networks. 

\FloatBarrier

\section{Discussion} 
\label{sec:discussion}

In this work, we propose a novel class of CRMs by mixing the L\'evy intensity of a stable or GG CRM over the stability parameter. The use of mixtures of existing L\'evy measures to obtain novel univariate or multivariate CRMs has a long history \citep{Cont2004,Griffin2017,Griffin2018,Ayed2019,Todeschini2020,Naik2021,Ayed2024}. A key novelty here is that the mixing is performed with respect to the index of variation of the base model, in order to obtain the desired rapid varying properties. The L\'evy intensities we consider have the form
\begin{align}
\int_0^1 g(s)h(w)w^{-1-s}\mathrm{d}s,
\label{eq:mixtureform}
\end{align} where $g(s)=\frac{s}{\Gamma(1-s)}$, $h(w)=1$ for mST \eqref{eq:rhoca} and $g(s)=\frac{\eta s c^s}{\Gamma(1-s)}$, $h(w)=e^{-\beta w/c}$ for mGG \eqref{eq:rhofull}. In both cases,
\begin{align}
g(s)&\underset{s\to 1}{\sim} C_1(1-s),~~~h(w)\underset{w\to 0}{\sim} C_2\label{eq:asymptoticgh}
\end{align}
 for some constants $C_1$ and $C_2$. These asymptotics imply, via Karamata's theorem, the behavior of the Laplace exponent at $\infty$ and L\'evy intensity a $0$ described in \cref{thm:AsymptoticsLevyLaplace}, and lead to extreme sparsity in the associated random graph. The specific forms of $g$ and $h$ were chosen to balance flexibility with analytical tractability: the Laplace exponents \eqref{eq:psica} and \eqref{eq:psifull} and associated moments admit simple expressions. This tractability also enables implementation of the size-biased sampling algorithm presented in \cref{sec:size-biased-sampling} without requiring numerical integration. More broadly, any mixture of the form in \eqref{eq:mixtureform} where $g$ and $h$ satisfy the asymptotics in \eqref{eq:asymptoticgh} could be used. For example, one can define a rapidly varying version of the stable beta process~\citep{Teh2009} by taking $g(s)=\frac{\eta s}{\Gamma(1-s)}$ and $h(w)=(1-w)^{\xi-1}1_{w\in(0,1)}$ for some $\eta,\xi>0$.

While we have focused on the application to sparse graph modeling, other domains are of significant interest. For instance, in clustering or partition models, such CRMs would yield random partitions with the number of clusters $K_n$ scaling as $n/\log(n)$. This extreme behavior may be desirable in problems such as entity resolution and deduplication~\citep{Zanella2016,DiBenedetto2021,Betancourt2022}, where most records correspond to unique individuals and thus $K_n$ is expected to grow nearly linearly with $n$.

In \cref{sec:graphs}, we introduced a class of extremely sparse random graph models based on the graphex framework. Another family of models -- edge-exchangeable networks~\citep{Crane2018,Cai2016,Janson2018} -- also employs CRMs as key components. Incorporating the CRMs developed in this paper within such models yields edge-exchangeable multigraphs whose number of edges also scales approximately linearly (up to a logarithmic factor) with the number of nodes.

\newpage
\appendix

\renewcommand{\thesection}{S\arabic{section}}
\renewcommand{\thetheorem}{S\arabic{theorem}}
\renewcommand{\theproposition}{S\arabic{proposition}}
\renewcommand{\thelemma}{S\arabic{lemma}}
\renewcommand{\thefigure}{S\arabic{figure}}
\renewcommand{\thetable}{S\arabic{table}}
\renewcommand{\theequation}{S\arabic{equation}}

  \begin{center}
    {\LARGE \textbf{ Rapidly Varying Completely Random Measures for Modeling Extremely Sparse Networks \\ Supplementary material}}
\end{center}
  \medskip

The supplementary material is organized as follows: \cref{sup:omitedproof} contains the proofs of the propositions and corollaries from \cref{sec:properties} and \cref{sec:graphs} of the main text. \cref{sup:mcmc} provides full details of the posterior inference algorithm. \cref{sup:experim} presents additional experiments. \cref{sup:identity,app:regularvariation,app:sizebiased} provide background on the Lambert W function, regularly varying functions, and size-biased sampling for CRMs, respectively. For clarity, all sections, theorems, propositions, lemmas in the supplementary material are prefixed with ``S" to distinguish them from those in the main text.

All the code used in this work is available at: \url{https://github.com/ValentinKil/rapidly_varying_crm}.



\section{Proofs}
\label{sup:omitedproof}
\subsection{Proofs of \cref{sec:properties}}

\subsubsection{Proof of \cref{thm:AsymptoticsLevyLaplace}}

The asymptotics for the Laplace exponent are obtained trivially. For the Lévy intensity, consider
\begin{align*}
(\alpha-\tau) x^{\tau+1}\rhoca(x;\alpha,\tau)  & =\int_{\tau}^{\alpha}\frac{s}{\Gamma
(1-s)}x^{\tau-s}\mathrm{d}s\\
& =\int_{0}^{\alpha-\tau}\frac{u+\tau}{\Gamma(1-u-\tau)}e^{-u\log
x}\mathrm{d}u
\end{align*}
with the change of variables $u=s-\tau$. Let $f(z)=\int_{0}^{\alpha-\tau}\frac{u+\tau}{\Gamma(1-u-\tau)}e^{-u z}du$. We have, as $u\to 0$,
\begin{align*}
\frac{u+\tau}{\Gamma(1-u-\tau)}\sim\left \{\begin{array}{ll}
                                     u & \text{if }\tau=0, \\
                                     \frac{\tau}{\Gamma(1-\tau)} & \text{if }\tau>0 .
                                   \end{array}\right .
\end{align*}
Using the Tauberian \cref{thm:RVTauberian}, we obtain, as $z\to\infty$
$$
f(z)\sim\left \{\begin{array}{ll}
                                     z^{-2} & \text{if }\tau=0, \\
                                     \frac{z^{-1}\tau}{\Gamma(1-\tau)} & \text{if }\tau>0 .
                                   \end{array}\right .
$$
Hence, as $x\to\infty$,
$$
(\alpha-\tau) x^{\tau+1}\rhoca(x;\alpha,\tau)=f(\log x)\sim \left \{\begin{array}{ll}
                                     \log^{-2}(x) & \text{if }\tau=0, \\
                                     \frac{\log^{-1}(x) \tau}{\Gamma(1-\tau)} & \text{if }\tau>0 .
                                   \end{array}\right .
$$
The asymptotics of the full model follow. For the behaviour at 0, we proceed similarly. Consider
\begin{align*}
(\alpha-\tau)x^{\alpha+1}\rhoca(x;\alpha,\tau)  & =\int_{\tau}^{\alpha}\frac{s}{\Gamma
(1-s)}x^{\alpha-s}\mathrm{d}s\\
& =\int_{0}^{\alpha-\tau}\frac{\alpha-u}{\Gamma(1+u-\alpha)}e^{-u\log(1/x)}\mathrm{d}u
\end{align*}
using the change of variables $u=\alpha-s$. Let $g(z)=\int_{0}^{\alpha-\tau}\frac{\alpha-u}{\Gamma(1+u-\alpha)}e^{-u z}\mathrm{d}u$. We have, as $u\to 0$,
\begin{align*}
\frac{\alpha-u}{\Gamma(1+u-\alpha)}\sim\left \{\begin{array}{ll}
                                     u & \text{if }\alpha=1, \\
                                     \frac{\alpha}{\Gamma(1-\alpha)} & \text{if }\alpha<1 .
                                   \end{array}\right .
\end{align*}
Using the Tauberian \cref{thm:RVTauberian}, we obtain, as $z\to\infty$
$$
g(z)\sim\left \{\begin{array}{ll}
                                     z^{-2} & \text{if }\alpha=1, \\
                                     \frac{z^{-1}\alpha}{\Gamma(1-\alpha)} & \text{if }\alpha<1.
                                   \end{array}\right .
$$
It follows that, as $x\to 0$,
$$
(\alpha-\tau)x^{\alpha+1}\rhoca(x;\alpha,\tau)=g(\log 1/x)\sim\left \{\begin{array}{ll}
                                     \log^{-2}(1/x) & \text{if }\alpha=1, \\
                                     \frac{\log^{-1}(1/x) \alpha}{\Gamma(1-\alpha)} & \text{if }\alpha<1.
                                   \end{array}\right .
$$
Hence the result for the full model.

\subsubsection{Proof of \cref{thm:moment}}

For $m\geq1$ and $z\geq0$, recall that 
\begin{align*}
\kappamix(m,z;\alpha,\tau,\beta,c,\eta)&=\int
_{0}^{\infty}w^{m}e^{-zw}\rhomix(w;\alpha,\tau,\beta,c,\eta)\mathrm{d}w,\\
\kappaca(m,z;\alpha,\tau)&=\int_{0}^{\infty}w^{m}e^{-zw}\rhoca(w;\alpha,\tau)\mathrm{d}w.
\end{align*}
 Note that
\[
\kappamix(m,z;\alpha,\tau,\beta,c,\eta)=\eta c^{m}\kappaca\left(
m,\beta+cz;\alpha,\tau\right).
\]
We have, for $z>0$
\begin{align*}
\kappaca\left(  m,z;\alpha,\tau\right)   &  =\int_{0}^{\infty}%
w^{m}e^{-zw}\frac{1}{\alpha-\tau}\int_{\tau}^{\alpha}\frac{s}{\Gamma
(1-s)}w^{-1-s}\mathrm{d}s\mathrm{d}w\\
&  =\frac{1}{\alpha-\tau}\int_{\tau}^{\alpha}\frac{s}{\Gamma(1-s)}\frac
{\Gamma(m-s)}{z^{m-s}}\mathrm{d}s\\
&  =\frac{z^{-m}}{\alpha-\tau}\int_{\tau}^{\alpha}sz^{s}\frac{\Gamma
(m-s)}{\Gamma(1-s)}\mathrm{d}s
\end{align*}
and $\kappaca\left(  m,0;\alpha,\tau\right)  =\infty$. For $m=1$, we
obtain
\[
\kappaca\left(  1,z;\alpha,\tau\right)  =\frac{z^{-1}}{\alpha-\tau
}\int_{\tau}^{\alpha}sz^{s}\mathrm{d}s.
\]
For $z=1$, we have $\kappaca\left(  1,1;\alpha,\tau\right)
=\frac{1}{\alpha-\tau}\frac{\alpha^{2}-\tau^{2}}{2}=\frac{\alpha+\tau}{2}%
$.\ For $z\neq1$, using the change of variable $u=z^{s}$, $du=z^{s}\log zds$, we obtain
\begin{align*}
\kappaca\left(  1,z;\alpha,\tau\right)   &  =\frac{z^{-1}}%
{\alpha-\tau}\int_{z^{\tau}}^{z^{\alpha}}sz^{s}ds=\frac{z^{-1}}{\alpha-\tau
}\frac{1}{\log^{2}z}\int_{z^{\tau}}^{z^{\alpha}}\log u\mathrm{d}u\\
&  =\frac{z^{-1}}{\alpha-\tau}\frac{1}{\log^{2}z}\left(  z^{\alpha}\log
z^{\alpha}-z^{\alpha}-\left[  z^{\tau}\log z^{\tau}-z^{\tau}\right]  \right)
\\
&  =\frac{z^{-1}}{\alpha-\tau}\frac{z^{\tau}-z^{\alpha}+(\alpha z^{\alpha
}-\tau z^{\tau})\log z}{\left(  \log z\right)  ^{2}}.
\end{align*}
For $m=2$, we have
\[
\kappaca\left(  2,z;\alpha,\tau\right)  =\frac{z^{-2}}{\alpha-\tau
}\int_{\tau}^{\alpha}s(1-s)z^{s}\mathrm{d}s.
\]
For $z=1$, this gives $\kappaca\left(  2,1;\alpha,\tau\right)
=\frac{1}{\alpha-\tau}\left[  \frac{\alpha^{2}}{2}-\frac{\alpha^{3}}{3}%
-(\frac{\tau^{2}}{2}-\frac{\tau^{3}}{3})\right]  $. For $z\neq1$,
\begin{align*}
\kappaca\left(  2,z;\alpha,\tau\right)   &  =\frac{z^{-2}}%
{\alpha-\tau}\int_{\tau}^{\alpha}s(1-s)z^{s}\mathrm{d}s\\
&  =z^{-1}\kappaca\left(  1,z;\alpha,\tau\right)  -\frac{z^{-2}%
}{\alpha-\tau}\int_{\tau}^{\alpha}s^{2}z^{s}\mathrm{d}s\\
&  =z^{-1}\kappaca\left(  1,z;\alpha,\tau\right)  -\frac{z^{-2}%
}{\alpha-\tau}\int_{z^{\tau}}^{z^{\alpha}}\frac{\log^{2}u}{\log^{3}z}\mathrm{d}u.
\end{align*}
using $\int\log^{2}u\mathrm{d}u=u\log^{2}u-2u\log u+2u+C$ gives the result.

\subsubsection{Proof of \cref{thm:sizebiasedmix}}

The proof follows directly from \cite{Perman1992}, whose construction is recalled in \cref{app:sizebiased}. Indeed, the conditional density $p(w|t)$ in \cref{app:sizebiased} takes the form

\begin{align*}
p(w|t)&=\frac{we^{-(t+\beta)w}\rhoca(w)}{\kappaca(1,t+\beta)}\\
&=\frac{\int_\tau^\alpha s/\Gamma(1-s)w^{-s}e^{-(t+\beta)w}\mathrm{d}s }{\int_\tau^\alpha s' (t+\beta)^{s'-1}\mathrm{d}s' }\\
&=\int_\tau^\alpha p(s|t) p(w|t,s)\mathrm{d}s
\end{align*}
where $p(s|t)=\frac{s(t+\beta)^{s}}{\int_\tau^\alpha s' (t+\beta)^{s'}ds'} \1{s\in(\tau,\alpha)}$ and $p(w|t,s)=w^{-s}e^{-(t+\beta)w}\frac{(t+\beta)^s}{\Gamma(1-s)}$.

\subsubsection{Proof of \cref{thm:totalmasslost}}

We can  follow the same proof strategy as in \citet{Lee2023}. Let $\alpha=1$, $\tau=0$. Without loss of generality, set $c=1$. According to \citet[Proposition 5.1]{Lee2023},  we have
\begin{align*}
\mathbb{E}[R_n|T_{n+1}]&=\int_0^{\infty}we^{-wT_{n+1}}\rhomix(w; 1,0,\beta,1,\eta)\mathrm{d}w\\
&=\kappamix(1,T_{n+1}; 1,0,\beta,1,\eta)\\
&=\eta \kappaca(1,\beta+T_{n+1};1,0)
\end{align*}
with $T_{n+1}  =\psica^{-1}\left(  \frac{\xi_{n+1}}{\eta}+\psica(\beta)\right)  -\beta$ where $\xi_{n+1}$ is a Gamma$(n+1,1)$ random variable. Similarly,
\begin{align*}
\mathbb{V}[R_n|T_{n+1}]&=\int_0^{\infty}w^2e^{-wT_{n+1}}\rhomix(w; 1,0,\beta,1,\eta)\mathrm{d}w\\
&=\kappamix(2,T_{n+1}; 1,0,\beta,1,\eta)\\
&=\eta \kappaca(2,\beta+T_{n+1};1,0).
\end{align*}
Using  \cref{thm:AsymptoticsLevyLaplace} and \cref{Ann:Asymptoticinverse} we get
 $$\psica^{-1}(x)\sim x\log(x)\text{~~as~~}x\rightarrow+\infty.$$
And so $$z_n:=\beta+T_{n+1}\sim\frac{n}{\eta}\log(n)\text{~~as~~}n\rightarrow \infty.$$
In particular, $T_{n+1}\rightarrow +\infty$, so as $n\to\infty$,
$$\mathbb{E}[R_n|T_{n+1}]=\eta \kappamix(1,z_n;1,0)=\frac{\eta}{z_n}\frac{1-\beta-T_{n+1}+(z_n)\log(z_n)}{\log^2(T_{n+1})}\sim\frac{\eta}{\log(z_n)}\sim\frac{\eta}{\log(n)}$$
and 
 \begin{align}
 \mathbb{V}[R_n|T_{n+1}] &=\eta \kappamix(2,z_n;1,0)=\frac{1}{z_n}\mathbb{E}[R_n|T_{n+1}]-\frac{\eta}{z_n^2}\int_1^{z_n}\frac{\log^2(u)\mathrm{d}u}{\log^3(z_n)}\nonumber\\
 &\sim\frac{\eta^2}{n\log^2(n)}-\frac{\eta^2 n\log^3(n)}{n^2\log^5(n)}+\frac{2\eta^2n\log^2(n)}{n^2\log^5(n)}\sim\frac{2\eta^2}{n\log^3(n)}.
 \end{align}
Finally, for all $\epsilon>0$, using Chebyshev's inequality, we obtain
$$\mathbb{P}\left(\left|\frac{R_n}{\mathbb{E}[R_n|T_{n+1}]}-1\right|>\epsilon|T_{n+1}\right)\leq\frac{\mathbb{V}[R_n|T_{n+1}]}{\epsilon^2\mathbb{E}[R_n|T_{n+1}]^2}$$
and as  $$\frac{\mathbb{V}[R_n|T_{n+1}]}{\mathbb{E}[R_n|T_{n+1}]^2}\sim\frac{2}{n\log(n)},$$
then by the Borel-Cantelli lemma, almost surely as $n\rightarrow+\infty$,
$$R_n\sim\mathbb{E}[R_n|T_{n+1}]\sim\frac{\eta}{\log(n)}.$$

\subsubsection{Proof of \cref{thm:mathapprox}}

The total mass $G(\Theta)$ has Laplace transform, for $t\geq 0$,
$$
\bbE[e^{-tG(\Theta)}]=\exp(-\psimix(t;\alpha,\tau,\beta,c,\eta)).
$$
We have
$$
\bbE [e^{-tc\left(\frac{\eta}{n}\right)^{\frac{1}{s_i}} X_i}]=\exp\left(-\frac{\eta}{n}((ct+\beta)^{s_i}-\beta^{s_i})\right).
$$
Hence
$$
\bbE [e^{-t S_n}]=\exp\left(-\frac{\eta}{n}\sum_{i=1}^n((ct+\beta)^{s_i}-\beta^{s_i})\right).
$$
Additionally, as
\[
\lim_{n\to \infty}\frac{1}{n}\sum_{i=1}^{n}t^{s_i} = \psica(t;\alpha,\tau),
\]
it follows that $\bbE [e^{-t S_n}]\to \bbE[e^{-tG(\Theta)}]$ as $n\to\infty$. As $\psimix(t;\alpha,\tau,\beta,c,\eta)$ is continuous, the result follows.

\subsection{Proofs of \cref{sec:graphs}}

Let
$$\bar\rho_{\text{mGG}}(x;\alpha,\tau,\beta,c,\eta) := \int_x^\infty \rhomix(w;\alpha,\tau,\beta,c,\eta)\mathrm{d}w $$
be the tail L\'evy intensity of the mGG CRM. The proofs of \cref{sec:graphs} mostly follow from \citet[Proposition 11, Theorem 2 and Corollary 1]{Caron2023}, which state the asymptotic properties of a CRM-based Caron-Fox graph model described in \cref{sec:CaronFox}, given the behavior of the tail L\'evy intensity at 0. Proposition 11 in \citet{Caron2023} requires that $\bar\rho_{\text{mGG}}(x;1,0,\beta,c,\eta)  \underset{x\rightarrow 0}{\sim} w^{-1}\widetilde\ell(1/x)$ for some slowly varying function $\widetilde\ell$. We prove in \cref{thm:AsymptoticsRhobar} below that this holds with
\begin{align*}
\widetilde\ell(1/x) =\frac{\eta c}{1-\tau}\frac{x^{-1}}{\log^2(1/x)}.
\end{align*}  
The mean
$$
m=\int_0^\infty w \rhomix(w;1,0,\beta,c,\eta)\mathrm{d}w=\left \{\begin{array}{cc}
                                       \eta c\frac{1-\beta+\beta\log(\beta)}{\beta\log^{2}(\beta)}  & \beta\neq 1,\ \\
                                    \frac{\eta c}{2} & \beta=1.
                                       \end{array}\right.
$$ is given by \cref{thm:moment}, and is finite for $\beta>0$. Define 
$$
\ell(x)=2m \widetilde\ell(x)=\frac{C}{\log^2(x)},
$$
where $$
C=\left \{\begin{array}{cc}
                                        2(\eta c)^2\frac{\beta^{-1}-1+\log \beta}{\left(  \log \beta\right)  ^{2}}& \beta\neq 1,\\
                                         (\eta c)^2 & \beta=1.
                                       \end{array}\right .
 $$        
Let                                 
$$
\ell_1(x)=\int_x^\infty y^{-1} \ell(y)\mathrm{d}y=\frac{C}{\log(x)}
$$
and define 
$$
\ell_1^*(y)=\left[\left\{\sqrt{\ell_1(\sqrt{y})}\right\}^{\#}\right]^2.
$$
As shown in \cref{lemma:debruijn} below, 
$$
\ell_1^*(y)=\frac{1}{2C}\log(y).
$$

\subsubsection{Secondary lemmas}

\begin{lemma}
\label{thm:AsymptoticsRhobar}
We have
$$
\bar\rho_{\text{mGG}}(x;1,0,\beta,c,\eta)  \underset{x\rightarrow 0}{\sim} x^{-1}\widetilde\ell(1/x) 
$$                                   
where
\begin{align*}
\widetilde\ell(1/x) =\frac{\eta c}{1-\tau}\frac{x^{-1}}{\log^2(1/x)}.
\end{align*}                                   
\end{lemma}

\begin{proof}
\cref{thm:AsymptoticsLevyLaplace} implies that $\rhomix(w;1,\tau,\beta,c,\eta)=w^{-2}L(1/x)$, where $L(1/x)\underset{x\to 0}{\sim} \frac{1}{(1-\tau)\log^2(1/w)}$. Using \cref{thm:RVKaramata2}, we obtain
$$
\bar\rho_{\text{mGG}}(x;1,\tau,\beta,c,\eta) \underset{x\rightarrow 0}{\sim} x^{-1}\frac{\eta c}{1-\tau}\frac{1}{\log^2(x)}.
$$
\end{proof}

\begin{lemma}
Let $\ell_1(x)=\frac{C}{\log(x)}$ for some constant $C$, and define
$$
\ell_1^*(y)=\left[\left\{\sqrt{\ell_1(\sqrt{y})}\right\}^{\#}\right]^2,
$$
where $\ell^\#$ denotes the de Bruijn conjugate (see \cref{Ann:Asymptoticinverse}) of the slowly varying function $\ell$. Then 
$$
\ell_1^*(y)=\frac{1}{2C}\log(y).
$$
\label{lemma:debruijn}
\end{lemma}
\begin{proof}
The pair $\left(\log(y),\log^{-1}(y)\right)$ is conjugate. Using \cref{thm:conjugaterules} the following pairs are also conjugate
$$\left(\frac{1}{4C}\log(y),4C\log^{-1}(y)\right),$$
$$\left(\sqrt{\frac{1}{4C}\log(y^2)},\sqrt{4C(\log^{-1}(y^2)}\right),$$
$$\left(\sqrt{\frac{1}{C}\log(\sqrt{y})},\sqrt{C\log^{-1}(\sqrt{y})}\right),$$
$$\left(\sqrt{\frac{1}{2C}\log(y)},\sqrt{\ell_1(\sqrt{y})}\right).$$
Hence $\ell_1^*(y)=\frac{1}{2C}\log(y)$.
\end{proof}

\subsubsection{Proof of \cref{thm:degdistrimix}}

Finally, while most of \cref{thm:degdistrimix} follows from Theorem 2 in \cite{Caron2023}, that theorem did not provide an exact asymptotic expression for the proportion of nodes of degree $j$ for $j\geq 2$, which we now derive.  Let  $\tilde{N}_{t,j}=\sum_{k\geq j}N_{t,k}$ denote the number of nodes of degree at least $j$. From \citet[Theorem 2]{Caron2023}, we have, almost surely,
$$\tilde{N}_{t,j}\sim \bbE(\tilde{N}_{t,j}).$$ Then, for $j\geq 2$, almost surely, using \citet[Theorem 1]{Caron2023},
$$
N_{t,j}=\tilde{N}_{t,j}-\tilde{N}_{t,j+1}\sim \bbE(\tilde{N}_{t,j})- \bbE(\tilde{N}_{t,j+1})= \bbE(\tilde{N}_{t,j}-\tilde{N}_{t,j+1})=\bbE(N_{t,j})\sim\frac{t^2}{j(j-1)}\frac{C}{\log^2(t)}. 
$$
The proof is similar for $j=1$. The dominance of nodes of degree 1 in the graph follows directly from $$\log(t)^{-2}=o(\log(t)^{-1}).$$ For \eqref{eq:powerlawdistribb}, we know from Theorem 2 in \cite{Caron2023} that
$$
\tilde{N}_{t,2}\sim \bbE(\tilde{N}_{t,2}) = \sum_{k\geq2}\bbE(N_{t,k}) \sim  \sum_{k\geq2} \frac{t^2}{j(j-1)}\frac{C}{\log^2(t)} = t^2\frac{C}{\log^2(t)}.
$$
Therefore, we obtain
$$\frac{N_{t,j}}{\tilde{N}_{t,2}} \sim  \frac{\frac{t^2}{j(j-1)}\frac{C}{\log^2(t)}}{t^2\frac{C}{\log^2(t)}}=\frac{1}{j(j-1)}.$$ 

\section{Directed multigraph posterior and MCMC Algorithm}
\label{sup:mcmc}
We give here the details of the algorithm introduced in \cref{sec:graphinference}.

\subsection{Directed multigraph posterior}
\label{sup:multigraph}

Formally, the symmetric 0-1 valued atomic measure \cref{eq:undirectedgraph} can be viewed as a transformation of a integer-valued atomic measure (also called directed multigraph) ; this view is useful when it comes to posterior inference. We note the directed multigraph
\begin{equation}
M=\sum_{i=0}^\infty\sum_{j=0}^\infty Q_{ij}\delta_{\theta_i,\theta_j},
\end{equation}
where $Q_{ij}$ counts the number of directed edges from node $i$ to node $j$, with locatiion $\theta_i$ and $\theta_j$. Given a CRM $G\sim\CRM(\rho,\lambda)$ for any Lévy measure $\rho$, $M$ is simply generated from a Poisson process with intensity given by the product measure $G\times G$. To derive the undirected graphs $U$ we simply put $Z_{ij}=Z_{ji}=\min(\tilde{q}_{ij},1)$ where $\tilde{q}_{ji}=\tilde{q}_{ij}=Q_{ij}+Q_{ji}$. This formulation is equivalent to the one introduce in \cref{sec:CaronFox} as explained in  \cite{Caron2017}.

For $t\in\bbR_+$ we note $M_t$  the measure $M$ retrained to the square $[0,t]^2$. Similarly we note $G_t$ the measure $G$ restrained to $[0,t]$ and $G_t^*$ the total masse of $G_t$. The directed multigraph posterior given $(Q_{ij})_{1\leq i,j\leq N_t}$ is characterise in the following theorem 

\begin{theorem}[Theorem 6, \citealp{Caron2017}]
\label{thm:multipost}
For $N_t\geq1$, let $\theta_1<...<\theta_{N_t}$ be the set of support points of the measure $M_t$. Let $w_i=G_t(\{\theta_i\})$ and $w_*=G_t^*-\sum_{i=1}^{N_t}w_i$. We have

$$
\begin{aligned}
& \mathbb{P}\left[\left(w_i \in \mathrm{~d} w_i\right)_{1:N_t}, w_* \in \mathrm{~d} w_* \mid (Q_{ij})_{1\leq i,j\leq N_t}\right] \\
& \propto \exp \left\{-\left(\sum_{i=1}^{N_t} w_i+w_*\right)^2\right\}\left\{\prod_{i=1}^{N_t} w_i^{m_i} \rho\left(\mathrm{d} w_i\right)\right\} W_t^*\left(\mathrm{d} w_*\right)
\end{aligned}
$$
where $m_i=\Sigma_{j=1}^{N_t}\tilde{q}_{ij}>0$ for $i=1, \ldots, N_t$ are the node degrees of the multigraph and $W_t^*$ is the probability distribution of the random variable $G_t^*$, with Laplace transform
$$
\mathbb{E}\left[\exp \left(-x G_t^*\right)\right]=\exp(-\psi(x)).
$$
Additionally, conditionally on observing an empty graph, \emph{i.e.}, $N_t=0$, we have
$$
\mathbb{P}\left[w_* \in \mathrm{~d} w_* \mid N_t=0\right] \propto \exp \left(-w_*^2\right) W_t^*\left(\mathrm{d} w_*\right).
$$
\end{theorem}

\subsection{Our MCMC Algorithm}
\subsubsection{Step 1 : update of the $w$ and $s$}

We propose to use Hamiltonian Monte Carlo to perform the step 1 of our previous algorithm which requires the computation of the gradient of the log-posterior. Following  \cref{thm:multipost}  we can compute the conditional probability.

\begin{align*}
&\mathbb{P}\left[(w_i\in \mathrm{d}w_i)_{0:N_t}|m_{1:N_t}, \phi, w_*\right]\\
\propto & \exp \left[-\left(\sum_{i=1}^{N_\alpha} w_i+w_*\right)^2\right]\left(\prod_{i=1}^{N_\alpha} w_i^{m_i} \rhomix\left(\mathrm{d} w_i\right)\right)\\
\propto & \idotsint_{s_{1:N_t}} \exp\left[-\left(\sum_{i=1}^{N_t}w_i+w_*\right)^2\right]\left(\prod_{i=1}^{N_t}w_i^{m_i}\hat{\rho}(w_i,s_i;\phi)\mathrm{d}w_i\right)\left(\prod_{i=1}^{N_t}\1{s_i\in (0,1)}\mathrm{d}s_i\right)\\
\propto & \idotsint_{s_{1:N_t}}\exp\left[-\left(\sum_{i=1}^{N_t}w_i+w_*\right)^2\right]\left(\prod_{i=1}^{N_t}w_i^{m_i}\frac{\eta s_i c^{s_i}}{\Gamma(1-s_i)}w_i^{-1-s_i}e^{-\beta w_i/c}\mathrm{d}w_i\right)\left(\prod_{i=1}^{N_t}\1{s_i\in (0,1)}\mathrm{d}s_i\right).
\end{align*}
So we can deduce the conditional PDF
\begin{align*}
p\left(w_{1:N_t},s_{1:N_t}|m_{1:N_t}, \phi, w_*\right) & \propto \exp\left[-\left(\sum_{i=1}^{N_t}w_i+w_*\right)^2\right]\left(\prod_{i=1}^{N_t}w_i^{m_i}\frac{\eta s_i c^{s_i}}{\Gamma(1-s_i)}w_i^{-1-s_i}e^{-\beta w_i/c}\1{s_i\in (0,1)}\right)\\
&  =e^{-\left(  \sum w_{i}+w_{\ast}\right)  ^{2}-\left(  \sum_{i}w_{i}\right)
\frac{\beta}{c}}\eta^{N}\left(\prod_{i=1}^{N}w_{i}^{m_{i}-s_{i}-1}\frac
{s_{i}c^{s_{i}}}{\Gamma(1-s_{i})}\1{s_i\in (0,1)}\right).%
\end{align*}
Let $sig(x)=\frac{1}{1+e^{-x}}$ with inverse $sig^{-1}(z)=\log(\frac{z}{1-z}%
)$.\ Note $sig^{\prime}(x)=sig(x)(1-sig(x))$. Do the change of variable
$v_{i}=\log(w_{i})$, and, $z_{i}=sig^{-1}(s_{i})$. Hence
$\mathrm{d}w_{i}=e^{v_{i}}\mathrm{d}v_{i}=w_{i}\mathrm{d}v_{i}$ and $\mathrm{d}s_{i}=sig(z_{i})(1-sig(z_{i}%
))\mathrm{d}z_{i}=s_{i}(1-s_{i})\mathrm{d}z_{i}$. It follows that
\begin{align*}
 p\left(v_{1:N_t},z_{1:N_t}|m_{1:N_t}, \phi, w_*\right) &  \propto e^{-\left(  \sum w_{i}+w_{\ast}\right)
^{2}-\left(  \sum_{i}w_{i}\right)  \frac{\beta}{c}}\eta^{N}\left[  \prod
_{i=1}^{N}w_{i}^{m_{i}-s_{i}-1}\frac{s_{i}c^{s_{i}}}{\Gamma(1-s_{i})}%
w_{i}s_{i}(1-s_{i})\right]  \\
&  =e^{-\left(  \sum w_{i}+w_{\ast}\right)  ^{2}-\left(  \sum_{i}w_{i}\right)
\frac{\beta}{c}}\eta^{N}\left[  \prod_{i=1}^{N}w_{i}^{m_{i}-s_{i}}\frac
{s_{i}^{2}(1-s_{i})}{\Gamma(1-s_{i})}c^{s_{i}}\right]  \\
&  =e^{-\left(  \sum w_{i}+w_{\ast}\right)  ^{2}-\left(  \sum_{i}w_{i}\right)
\frac{\beta}{c}}\eta^{N}\left[  \prod_{i=1}^{N}w_{i}^{m_{i}-s_{i}}\frac
{s_{i}^{2}(1-s_{i})^{2}}{\Gamma(2-s_{i})}c^{s_{i}}\right],
\end{align*}
where the last line comes from $\Gamma(2-s)=(1-s)\Gamma(1-s)$.\ This form is
preferred to avoid numerical instabilities when $s_{i}$ is close to 1. We
obtain (up to a constant independent of $v_{i},z_{i}$)
\begin{align*}
\ln(p\left(v_{1:N_t},z_{1:N_t}|m_{1:N_t}, \phi, w_*\right))& =-\left(  \sum w_{i}+w_{\ast}\right)  ^{2}-\left(
\sum_{i}w_{i}\right)  \frac{\beta}{c}\\
& +\sum_{i=1}^{N}\left(  m_{i}-s_{i}\right)  \log w_{i}+2\log s_{i}%
+2\log(1-s_{i})-\log\Gamma(2-s_{i})+s_{i}\log(c)\\
& :=g(v_{1:N_t},z_{1:N_t}).
\end{align*}
Now we can compute the gradient of the log posterior. We have%
\begin{align}
\label{gradlogpost1}
\frac{\partial g(v_{1:N_t},z_{1:N_t})}{\partial v_{i}} &  =\frac{\partial w_{i}}{\partial v_{i}}%
\frac{\partial g(v_{1:N_t},z_{1:N_t})}{\partial w_{i}}\nonumber\\
&  =w_{i}\times\left[  -2\left(  \sum w_{i}+w_{\ast}\right)  -\frac{\beta}%
{c}+\frac{m_{i}-s_{i}}{w_{i}}\right]  \\
&  =m_{i}-s_{i}-w_{i}\left[  2\left(  \sum w_{i}+w_{\ast}\right)  +\frac
{\beta}{c}\right]\nonumber
\end{align}
and%
\begin{align}
\label{gradlogpost2}
\frac{\partial g(v_{1:N_t},z_{1:N_t})}{\partial z_{i}} &  =\frac{\partial s_{i}}{\partial z_{i}}%
\frac{\partial g(v_{1:N_t},z_{1:N_t})}{\partial s_{i}}\nonumber\\
&  =s_{i}(1-s_{i})\times\left[  -\log w_{i}+\frac{2}{s_{i}}-\frac{2}{1-s_{i}%
}+\psi(2-s_{i})+\log c\right]  \\
&  =s_{i}(1-s_{i})\times\left[  \psi(2-s_{i})-\log w_{i}+\log c\right]
+2-4s_{i},\nonumber%
\end{align}
where $\psi=\frac{\Gamma'}{\Gamma}$ is the bigamma function which is implemented is scipy. Noting that $\psi(x+1)=\psi(x)+\frac{1}{x}$, the above is equivalent to
$$
\frac{\partial g(v_{1:N_t},z_{1:N_t})}{dz_{i}}=s_{i}(1-s_{i})\times\left[
\psi(1-s_{i})-\log w_{i}+\log c\right]  +2-3s_{i}%
$$
but the first expression should be preferred.

In order to perform Step 1, we first update the $w_{1:N_t}$ with $s$ fixed and then the $s_{1:N_t}$ with $w$ fixed, for both update we use an HMC update with momentum variable $\bm{p_w}$ and $\bm{p_s}$. We refer the reader to \cite{Neal2011} and \cite{Betancourt2018} for an overview. We start by explaining how we update the $w$, the update of the latent $s$ is similar.

\paragraph{Update of the $w$.}

Let $L\geq1$ be the number of Leapfrog steps and $\epsilon_w>0$ the step size. The potential energy is given by
$$
U(w,s,m_{1:N_t},\phi, w_*)=\ln(p\left(v_{1:N_t},z_{1:N_t}|m_{1:N_t}, \phi, w_*\right))=g(u,v)
$$
whose gradient $U'_w$ has been given in \ref{gradlogpost1}. The algorithm proceeds by first sampling momentum variable as
\begin{align*}
\bm{p_w} &\sim\mathcal{N}(0,I_{N_t})\\
\end{align*}
The Hamiltonian proposal $q(\tilde{\bm{w}},\tilde{\bm{p_w}},|\bm{w},\bm{p_w})$ (we drop out the $1:N_t$ for lisibility) is obtained by the following Leapfrog algorithm : simulate $L$ steps of discretised Hamiltonian via
 \begin{align*}
 \tilde{\bm{p_w}}^{(0)}&=\bm{p_w}+\frac{\epsilon_w}{2}U'_w(\bm{w},\bm{s},w_*,\phi);   &\tilde{\bm{w}}^{(0)}&=\bm{w}.
 \end{align*}
 And for $l=1,...,L-1$,
  \begin{align*}
  \log(\tilde{\bm{w}}^{(l)})&= \log(\tilde{\bm{w}}^{(l-1)})+\epsilon_w \tilde{\bm{p_w}}^{(l-1)};  & \tilde{\bm{p_w}}^{(l)}&= \tilde{\bm{p_w}}^{(l-1)}+\epsilon_w U'_w(\tilde{\bm{w}}^{(l)},\bm{s},w_*,\phi).\\
 \end{align*}
Finally, set
   \begin{align*}
  \log(\tilde{\bm{w}})&= \log(\tilde{\bm{w}}^{(L-1)})+\epsilon_w \tilde{\bm{p_w}}^{(L-1)};  &  \tilde{\bm{p_w}}&= \tilde{\bm{p_w}}^{(l-1)}+\frac{\epsilon_w}{2} U'_w(\tilde{\bm{w}},\bm{s},w_*,\phi).\\
 \end{align*}
 The proposal  $(\tilde{\bm{w}},\tilde{\bm{p_w}})$ is accepted with probability $\min(1,r_w)$ with
$$
r_w = \exp\left[-\frac{1}{2}\sum_{i=1}^{N_t}(\tilde{p}_{w,i}^2 - p_{w,i}^2)\right]
\frac{p\left(\bm{\tilde{v}} \mid \bm{z},  G, \phi, w_*\right)}{p\left(\bm{v} \mid \bm{z}, G, \phi, w_*\right)}.
$$

\paragraph{Update of the latent $s$.}
Once $w$ has been updated, the update of $s$ follows the same principle: the algorithm proceeds by first sampling momentum variable as
\begin{align*}
\bm{p_s} &\sim\mathcal{N}(0,I_{N_t}).
\end{align*}
 The Hamiltonian proposal $q(\tilde{\bm{s}},\tilde{\bm{p_s}},|\bm{s},\bm{p_s})$ is obtained by $L$ Leapfrog step:
 \begin{align*}
  \tilde{\bm{p_s}}^{(0)}&=\bm{p_s}+\frac{\epsilon_s}{2}U'_s(\bm{w},\bm{s},w_*,\phi);   &\tilde{\bm{s}}^{(0)}&=\bm{s}.
 \end{align*}
 And for $l=1,...,L-1$,
  \begin{align*}
 \log(\tilde{\bm{s}}^{(l)})&= \log(\tilde{\bm{s}}^{(l-1)})+\epsilon_s \tilde{\bm{p_s}}^{(l-1)} & \tilde{\bm{p_s}}^{(l)}&=\tilde{\bm{p_s}}^{(l-1)}+\epsilon_s U'_s(\bm{w},\tilde{\bm{s}}^{(l)},w_*,\phi);
 \end{align*}
and finally set
 \begin{align*}
  \logit(\tilde{\bm{s}})&= \logit(\tilde{\bm{s}}^{(L-1)})+\epsilon_s\tilde{\bm{p_s}}^{(L-1)} &\tilde{\bm{p_s}}=\tilde{\bm{p_s}}^{(l-1)}+\frac{\epsilon_s}{2} U'_s(\bm{w},\tilde{\bm{s}},w_*,\phi).
 \end{align*}
The proposal  $(\tilde{\bm{s}},\tilde{\bm{p_s}})$ is accepted with probability $\min(1,r_s)$ with
$$
r_s = \exp\left[-\frac{1}{2}\sum_{i=1}^{N_t}(\tilde{p}_{s,i}^2 - p_{s,i}^2)\right] \frac{p\left(\bm{\tilde{z}} \mid \bm{v},  G, \phi, w_*\right)}{p\left(\bm{z} \mid \bm{v}, G, \phi, w_*\right)}.
$$

\subsubsection{Step 2 : update of hyperparameters}

Let $G^*_\phi$ be the probability distribution of the random variable of $G(\Theta)$. Recall that we have $\mathbb{E}\left[  e^{-tG(\Theta)}\right]=e^{-\psi(t)}$ for $t\geq0$. Let $g^*_\phi$ denote the probability density of  $G^*_\phi$ with respect of the Lebesgue mesure.  First we can note the interesting property :

\begin{equation}
\label{magictrick}
g^*_{\alpha, \tau, \beta_1+c\beta_2, c, \eta}(x)=g^*_{\alpha, \tau, \beta_1, c, \eta}(x)\frac{\exp(-x\beta_2)}{\exp(-\eta\psimix(\beta_2; \alpha,\tau,\beta_1,c,1)}.
\end{equation}

\begin{proof} For all $t\in\bbR$,
\begin{align*}
&\int_0^{+\infty} e^{-tx}g^*_{\alpha, \tau, \beta_1, c, \eta}(x)\frac{\exp(-x\beta_2)}{\exp(-\eta\psimix(\beta_2; \alpha,\tau,\beta_1,c,1)}\mathrm{d}x\\
=& \exp(-\eta\psimix(\beta_2; \alpha,\tau,\beta_1,c,1)\int_0^{+\infty} e^{-tx}g^*_{\alpha, \tau, \beta_1, c, \eta}(x)\exp(-x\beta_2)\mathrm{d}x\\
=& \exp(-\eta\psimix(\beta_2; \alpha,\tau,\beta_1,c,1)\int_0^{+\infty} e^{-x\left(t+\beta_2\right)}g^*_{\alpha, \tau, \beta_1, c, \eta}(x)\mathrm{d}x\\
=& \exp(-\eta\psimix(\beta_2; \alpha,\tau,\beta_1,c,1)\exp\left(-\psimix\left(t+\beta_2,\alpha,\tau,\beta_1,c,\eta\right)\right)\\
= & \exp\left[\eta\psica(\beta_1+c\beta_2,\alpha,\tau)-\eta\psica(\beta_1,\alpha,\tau)-\eta\psica(\beta_1+c\beta_2+ct,\alpha,\tau)+\eta\psica(\beta_1,\alpha,\tau)\right]\\
= & \exp(-\psimix(t,\alpha,\tau,\beta_1+c\beta_2,c,\eta).
\end{align*}
\end{proof}
Following the same strategy than in \citet[Appendix F.2]{Caron2017}, we propose $(\tilde{\phi},\tilde{w_*})$  from $q(\tilde{\phi},\tilde{w}_* |\phi, w_*)$ and accept with probability $\min(1,r)$ where
\begin{align*}
r&=\frac{\exp\left\{-\left(\sum_{i=1}^{N}w_i+\tilde{w}_*\right)^2\right\}}{\exp\left\{-\left(\sum_{i=1}^{N}w_i+w_*\right)^2\right\}}
\left\{\prod_{i=1}^N\frac{\hat{\rho}(w_i,s_i;\tilde{\phi})}{\hat{\rho}(w_i,s_i;\phi)}\right\}
\frac{g^*_{\tilde{\phi}}(\tilde{w}_*)}{g^*_{\phi}(w_*)}
\frac{p(\tilde{\phi})}{p(\phi)}
\frac{q(\phi,w_*|\tilde{\phi},\tilde{w}_*)}{q(\tilde{\phi},\tilde{w}_*|\phi,w_*)}\\
&=\frac{\exp\left\{-\left(\sum_{i=1}^{N}w_i+\tilde{w}_*\right)^2\right\}}{\exp\left\{-\left(\sum_{i=1}^{N}w_i+w_*\right)^2\right\}}
\left\{\left(\frac{\tilde\eta}{\eta}\right)^N \left(\frac{\tilde c}{c}\right)^{\sum_{i=1}^N s_i} e^{-(\tilde\beta/\tilde c - \beta/c)\left(\sum_{i=1}^N w_i\right) } \right\}
\frac{g^*_{\tilde{\phi}}(\tilde{w}_*)}{g^*_{\phi}(w_*)}
\frac{p(\tilde{\phi})}{p(\phi)}
\frac{q(\phi,w_*|\tilde{\phi},\tilde{w}_*)}{q(\tilde{\phi},\tilde{w}_*|\phi,w_*)}.
\end{align*}
We shall use the proposal
$$q(\tilde{\phi},\tilde{w}_*|\phi,w_*)=q(\tilde{\beta}|\beta)q(\tilde{c}|c)q(\tilde{\eta}|\tilde{\beta},\tilde{c},\eta,w_*)q(\tilde{w}_*|\tilde{\phi},w_*)$$
where
\begin{align*}
q(\tilde{\beta}|\beta)&=\text{lognormal}(\tilde{\beta};\log(\beta),\sigma^2_\beta),\\
q(\tilde{c}|c)&=\text{lognormal}(\tilde{c};\log(c),\sigma^2_c),\\
q(\tilde{w_*}|\tilde{\phi},w_*)&=g^*_{\alpha,\tau,\tilde{\beta}+2\tilde{c}\sum_{i=1}^{N}w_i+2\tilde{c}w_*,\tilde{c},\tilde{\eta}}(\tilde{w}_*).
\end{align*}
For $\eta$, we alternatively use two proposals
\begin{align}
q(\tilde{\eta}|\tilde{\beta},\tilde{c},\eta,w_*) &= \text{Gamma}\left(\tilde{\eta};N,\psi^1_{\tilde{\phi}}\left(2\sum_{i=1}^{N}w_i+2w_*\right)\right)
\label{eq:proposaleta1}
\end{align}
where $\psi^1_\phi(t)=\psimix(t;\alpha,\tau,\beta,c,1)$, or a simple random walk on $\log \eta$
\begin{align}
q(\tilde{\eta}|\eta)&=\text{lognormal}(\tilde{\eta};\log(\eta),\sigma^2_\eta).
\label{eq:proposaleta2}
\end{align}
Note that, under the priors \eqref{eq:priorshyper},
$$
\frac{p(\tilde\beta)q(\beta| \tilde\beta)}{p(\beta)q(\tilde{\beta}|\beta)}=\frac{p(\tilde c)q(c| \tilde c)}{p(c)q(\tilde{c}|c)}=1
$$
and
$$
\frac{g^*_{\tilde{\phi}}(\tilde{w}_*)}{g^*_{\phi}(w_*)}\times \frac{q(w_*|\phi,\tilde w_*)}{q(\tilde{w_*}|\tilde{\phi},w_*)}=\frac{e^{-w_*(2\sum_{i=1}^Nw_i+2\tilde{w}_*)}}{e^{-\tilde w_*(2\sum_{i=1}^Nw_i+2w_*)}}
\times\frac{e^{-\tilde\eta \psi^1_{\tilde{\phi}}(2\sum_{i=1}^Nw_i+2w_*) }}{e^{-\eta \psi^1_{\phi}(2\sum_{i=1}^Nw_i+2\tilde w_*) }}.
$$
Sometimes it is useful to put a gamma prior on $\beta$ and $c-1$:
$$p(\beta)=\text{Gamma}(\beta;a_\beta,b_\beta), ~~~p(c-1)=\text{Gamma}(c;a_c,b_c).$$
In that case we have
$$
\frac{p(\tilde\beta)q(\beta|\tilde\beta)}{p(\beta)q(\tilde{\beta}|\beta)}=\left(\frac{\tilde{\beta}}{\beta}\right)^{a_\beta}\exp(-b_\beta(\tilde{\beta}-\beta)) ,~\text{and}~~\frac{p(\tilde{c})q(c|\tilde{c})}{p(c)q(\tilde{c}|c)}=\left(\frac{\tilde{c}}{c}\right)^{a_c}\exp(-b_c(\tilde{c}-c)).
$$
Note that we can recover the improper prior by taking $a=b=0$ in both case.
It follows that the acceptance ratio is
\begin{align*}
r&=\frac{\exp\left[-\left(\sum_{i=1}^{N}w_i+\tilde{w}_*\right)^2\right]}{\exp\left[-\left(\sum_{i=1}^{N}w_i+w_*\right)^2\right]}
\left\{\left(\frac{\tilde c}{c}\right)^{\sum_{i=1}^N s_i} e^{-(\tilde\beta/\tilde c - \beta/c + 2w_*-2\tilde w_*)\left(\sum_{i=1}^N w_i\right) } \right\}\\
&\times\left(\frac{\tilde{\beta}}{\beta}\right)^{a_\beta}\exp(-b_\beta(\tilde{\beta}-\beta))\times \left(\frac{\tilde{c}}{c}\right)^{a_c}\exp(-b_c(\tilde{c}-c))\\
& \times\underbrace{ \left (\frac{\tilde\eta}{\eta}\right)^{N-1}
\frac{e^{-\tilde\eta \psi^1_{\tilde{\phi}}(2\sum_{i=1}^Nw_i+2w_*) }}{e^{-\eta \psi^1_{\phi}(2\sum_{i=1}^Nw_i+2\tilde w_*) }}
\frac{q(\eta|\beta,c,\tilde \eta,\tilde w_*)}{q(\tilde{\eta}|\tilde{\beta},\tilde{c},\eta,w_*)}}_{r_2}.
\end{align*}
Under the gamma proposal \eqref{eq:proposaleta1} for $\eta$,
$$
r_2=\left(\frac{\psi^1_{\phi}(2\sum_{i=1}^Nw_i+2\tilde{w}_*)}{\psi^1_{\tilde{\phi}}(2\sum_{i=1}^Nw_i+2w_*)}\right)^N.
$$
Under the log-normal proposal \eqref{eq:proposaleta2} for $\eta$,
$$
r_2=\left (\frac{\tilde\eta}{\eta}\right)^{N}
\frac{e^{-\tilde\eta \psi^1_{\tilde{\phi}}(2\sum_{i=1}^Nw_i+2w_*) }}{e^{-\eta \psi^1_{\phi}(2\sum_{i=1}^Nw_i+2\tilde w_*) }} .
$$

\subsubsection{Step 3 : update of the latent counts}

We simply update the latent counts  by using the conditional distribution \cref{latentcountdistrib}. 

\section{Additional experiments}
\label{sup:experim}

\subsection{Synthetic Data}

\subsubsection{Proportion of degree one}

Under our graph model, \cref{{thm:degdistrimix}} states that the proportion of nodes of degree one converges to 1 as the size of the graph goes to infinity. This phenomenon is illustrated in \cref{{fig:degreeone}}.

\begin{figure}[ht]
    \centering
        \includegraphics[scale=0.5]{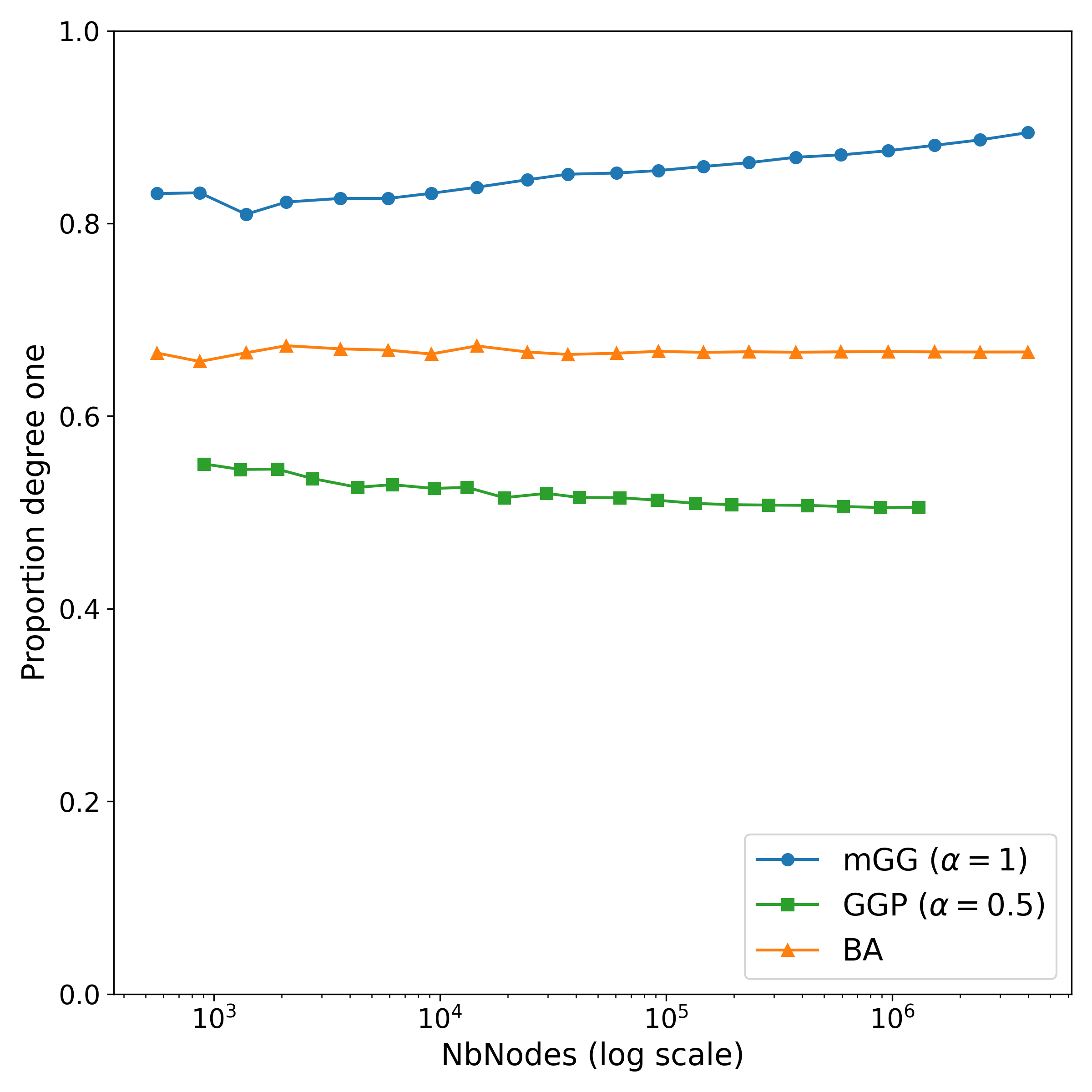}
   \caption{Proportion of nodes of degree one  in the mGG model ({\Large$\bullet$}) with parameters $\alpha=1$, $\tau=0$, $c=1$, and $\beta=1$ for various values of $\eta$ ranging from 50 to 6000, resulting in graphs of different sizes. Comparison with the Generalised Gamma CRM \citep{Caron2017} ({\footnotesize $\blacksquare$}) with parameters $\tau=1$ and $\sigma=0.5$, and with the Barabási–Albert model \citep{Barabasi1999a} ($\blacktriangle$). For every configuration we simulate 20 graph samples and plot the mean of the quantity of interest. These simulations were perfomed using the Python implementation available on Github running on an Apple M2 Chips with 16Gb of RAM.}
   \label{fig:degreeone}
\end{figure}

\subsubsection{Cluster}

\begin{figure}[ht]
    \centering
    \begin{subfigure}{0.45\textwidth}
        \centering
        \includegraphics[width=\textwidth]{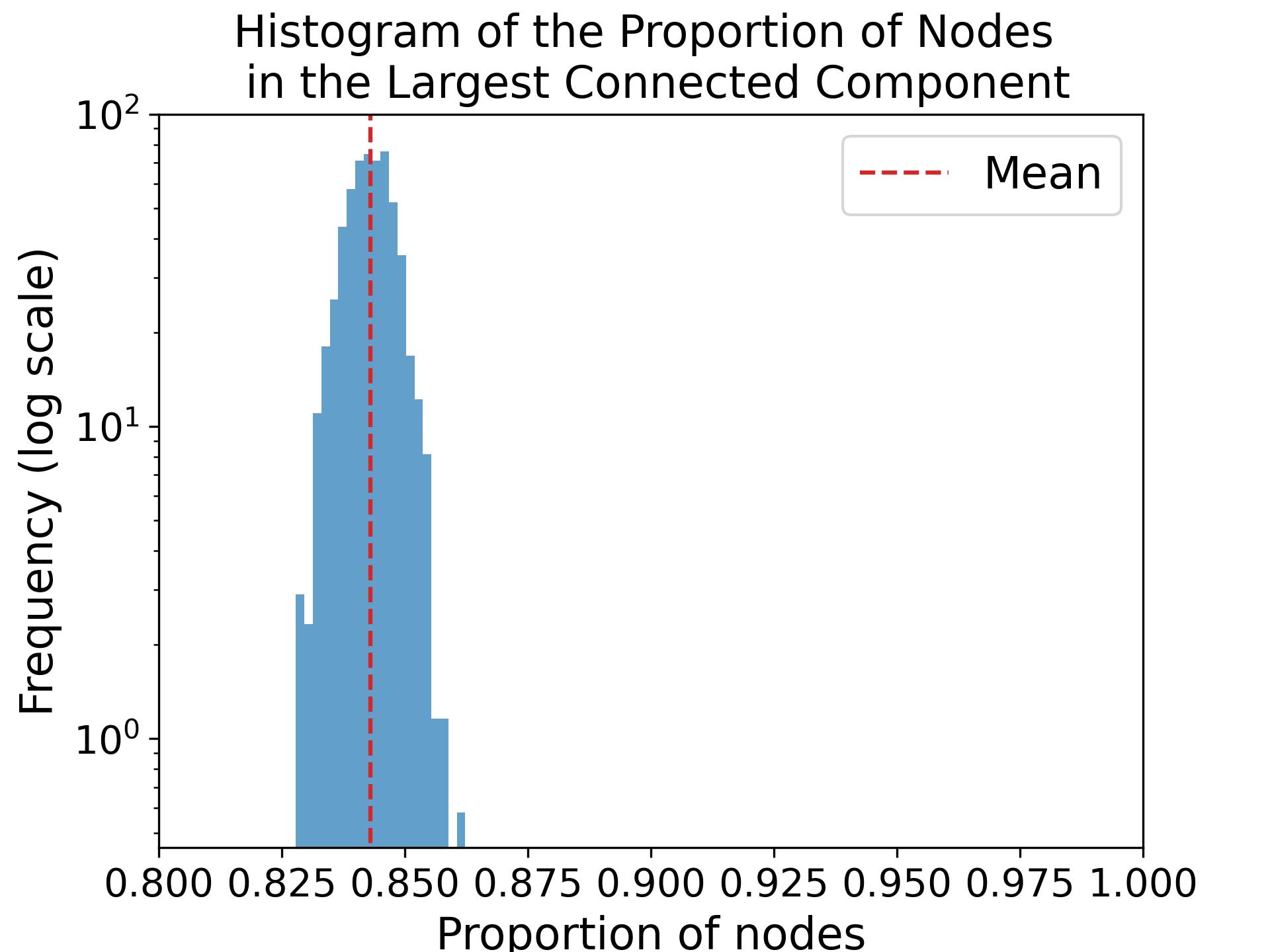}
        \caption{}
    \end{subfigure}
    \begin{subfigure}{0.43\textwidth}
        \centering
        \includegraphics[width=\textwidth]{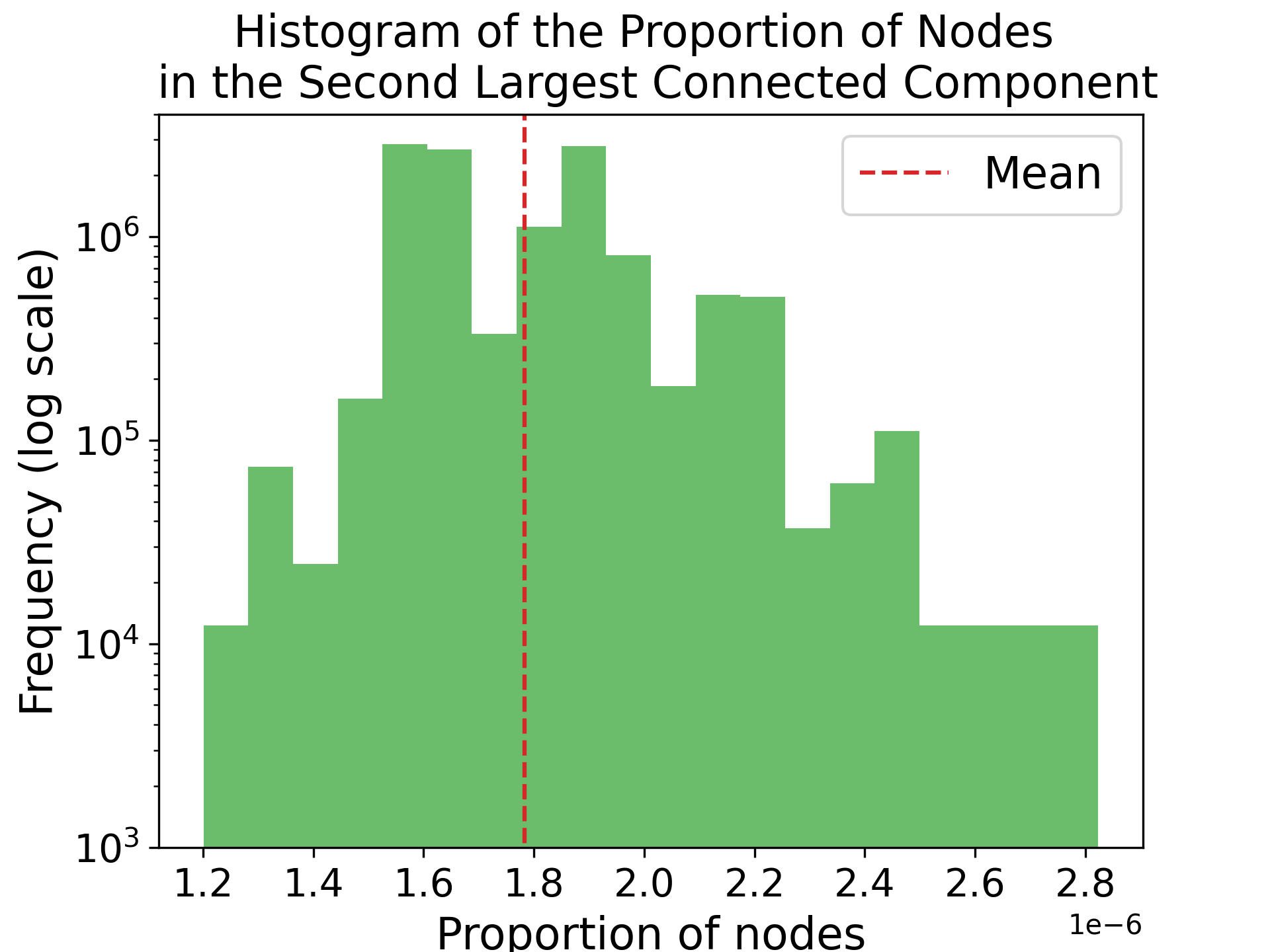}
        \caption{}
    \end{subfigure}
        \caption{Proportion of nodes in the largest connected component (a), second largest connected component (b) for 1000 samples of graph generate with parameters $\alpha=0, \tau=1, \beta=1, c=8$ and $\eta=800$, the typical size of such a graph is 3 574 053 nodes. The mean is 0.842 for the largest connected component and $1.78\times10^{-6}$ for the second largest. These simulations were performed using the Python implementation available on Github running on a standard computer (13th Gen Intel Core i5-1350P). }
\label{fig:propCC}
\end{figure}

\begin{figure}[ht]
\centering
 \includegraphics[scale=0.6]{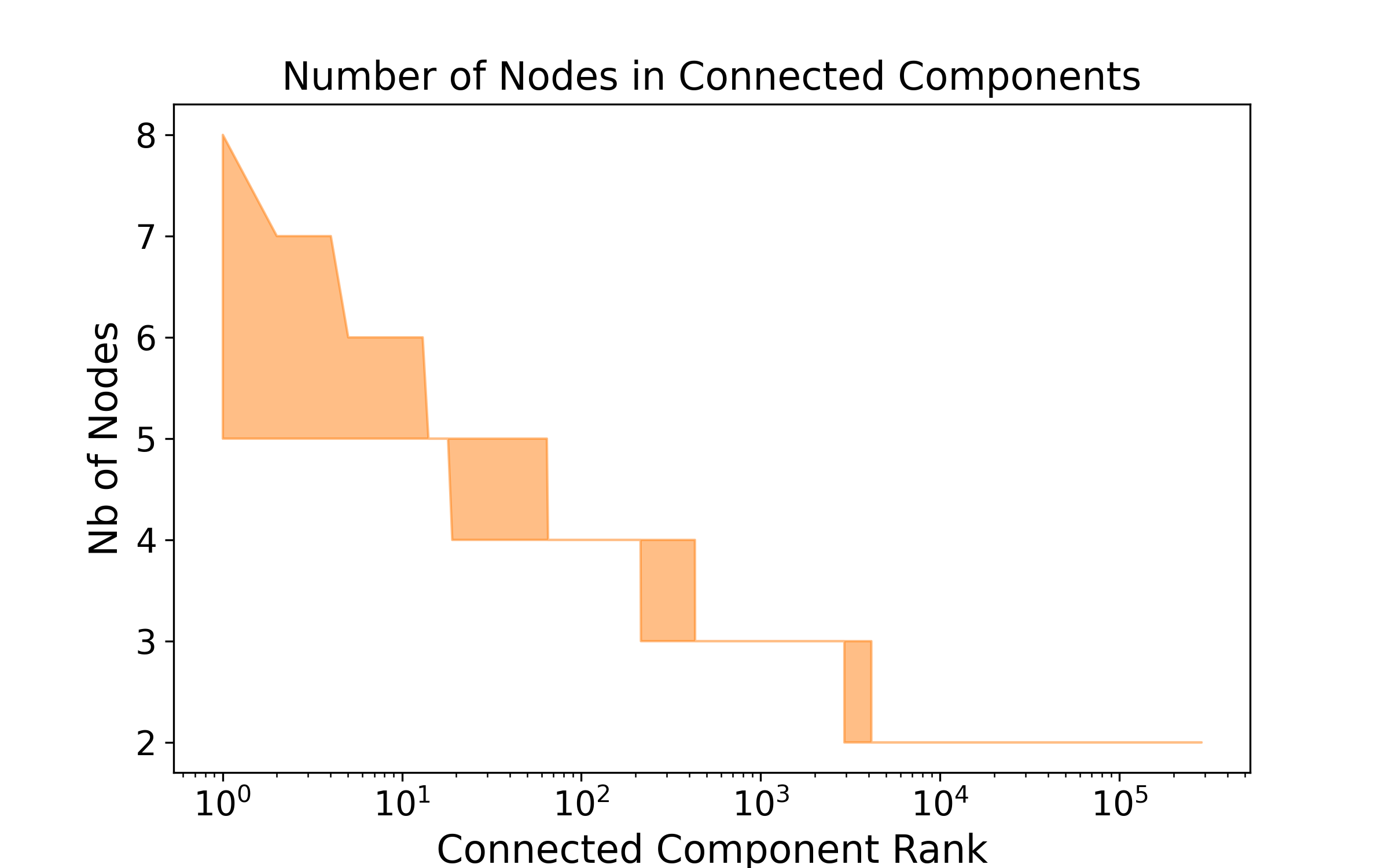}
 \caption{Number of nodes in a connected component with respect of its rank, where the connected components are ranked by decreasing size. The number of nodes in the largest connected component is not displayed. We have superposed the value for the same 1000 graphs than the ones generated for \cref{fig:propCC}.}
 \label{fig:rankCC}
\end{figure}

In order to get a more precise idea of what the cluster structure of the generated networks looks like, we performed additional simulations. These simulations were performed using the Python implementation available on GitHub, running on an Apple M2 chip with 16 GB of RAM. We observe that the generated network consists of one large connected component containing most of the nodes (see \cref{fig:propCC}(a)), along with some very small connected components. As the size of the graph increases (\emph{i.e.}, as $\eta$ grows), the proportion of nodes in the giant connected component also increases. The second-largest component is already of negligible size (see \cref{fig:propCC}(b)). As shown in \cref{fig:rankCC}, apart from the largest component, the other connected components are of small size, with the large majority containing only two nodes.

\FloatBarrier
\subsection{Inference}
\label{sup:inference}

\subsubsection{Details on the inference procedure}
In this section, we provide additional details on the inference procedure from \cref{art:synthetic}. Note that the same procedure is used for real-world data. The weights $w$ and latent $s$ are updated via our HMC procedure with $L=5$ leapfrog steps. As suggested by \citet{Betancourt2018}, the HMC step sizes $\epsilon_s$ and $\epsilon_w$ are automatically adapted during the first half of the iterations to target an acceptance rate of 0.65, following the update rule:
$$\epsilon(t+1)=\epsilon(t)\times\exp(0.005(\text{rate}(t)-0.65)).$$

For the estimated hyperparameters, the standard deviations of the random walk Metropolis-Hastings steps for $\log(\beta), \log(c)$ and $\log(\eta)$ were set to 0.01, 0.01 and 0.02 respectively.  The computing time for running each chain was about 14 hours  using the Python implementation available on Github on a cluster (CPU : Intel Platinum 8628, 2.90GHz).
The latent counts $\tilde{q}_{ij}$ are updated only once every 30 iterations. Every other iteration, the parameter $\eta$ is updated using a simple random walk Metropolis-Hastings step. A thinning procedure is applied to retain only 1,000 iterations per chain.

\subsubsection{Diagnostics of convergence}

To assess the convergence of our procedure, we use the Gelman-Rubin convergence diagnostic \citep{Gelman1992}. The exact  GR values for each hyperparameters are provided in \cref{tab:GR}. Given the large number of parameters, we also report a multivariate version of the GR diagnostic, proposed by \citet{Vats2021} and defined as
\begin{equation}
\label{equa:Rmulti}
R_{multi}=\sqrt{\frac{n_{it}-1}{n_{it}}+\frac{\left(\prod_{i=0}^{n_{par}}X_i\right)^{\frac{1}{n_{par}}}}{n_{it}}},
\end{equation}
where $n_{it}$ is the number of MCMC iteration after the burning phase,  $n_{par}$ is the number of estimated parameters, and $X_i = n_{it} \times R_i^2 - n_{it} + 1$, with $R_i$ denoting the individual GR diagnostic for parameter $i$, where $1 \leq i \leq n_{par}$. In our case, we obtain $R_{multi} = 1.001$, suggesting the convergence of the algorithm.

\begin{table}[ht]
    \centering
    \begin{tabular}{l c}
        Parameter & Value \\
        \hline
        $\beta$ & 1.004 \\
        $c$ & 1.045 \\
        $\eta$ & 1.045 \\
        \hline
    \end{tabular}
    \caption{Values of the Gelman-Rubin convergence diagnostic.}
    \label{tab:GR}
\end{table}

\FloatBarrier
\subsection{Real World Data}
\label{sup:realworld}

\begin{figure}[ht]
    \centering
    \begin{subfigure}{0.45\textwidth}
        \centering
        \includegraphics[width=\textwidth]{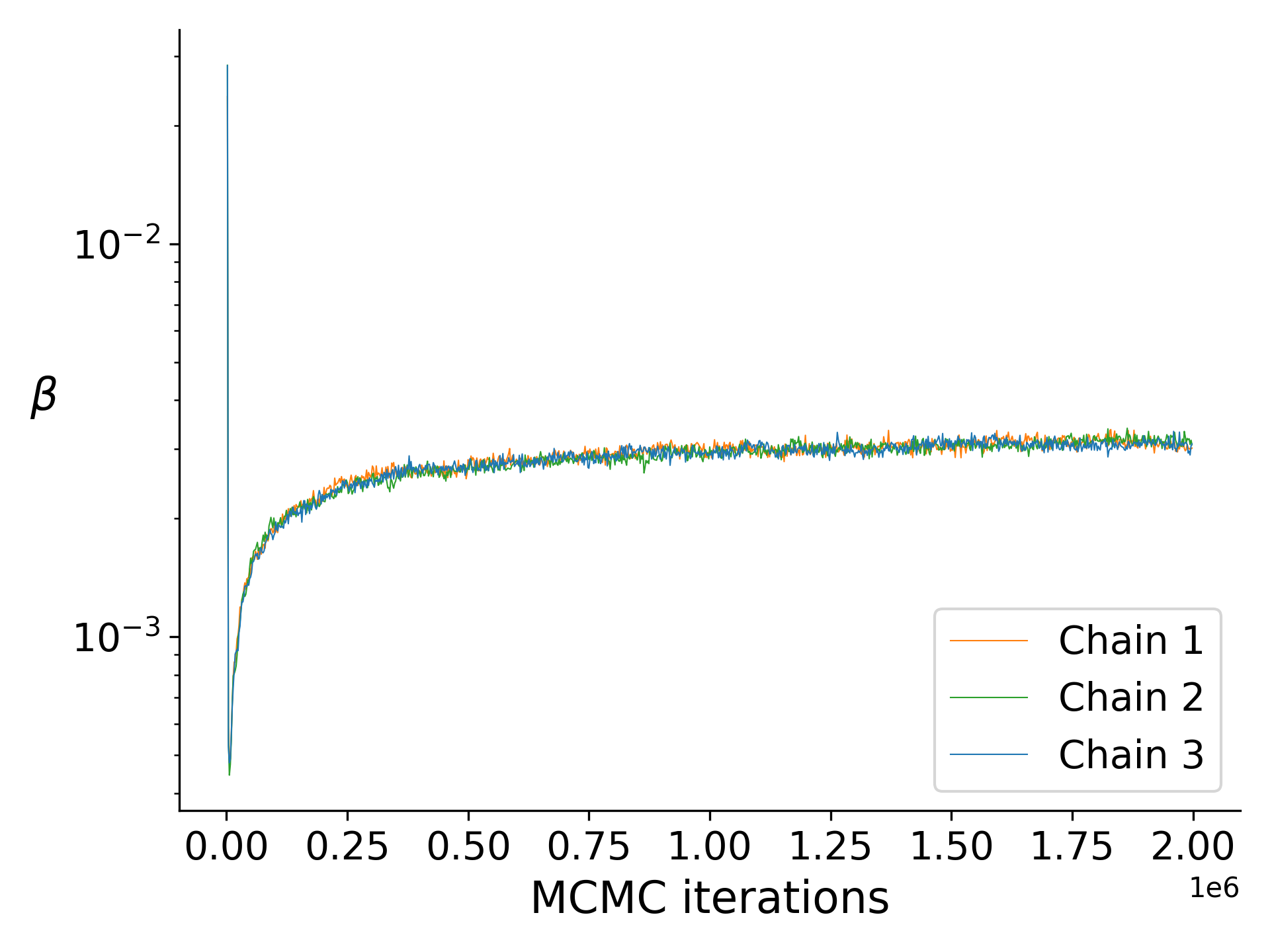}
        \caption{}
    \end{subfigure}
    \begin{subfigure}{0.45\textwidth}
        \centering
        \includegraphics[width=\textwidth]{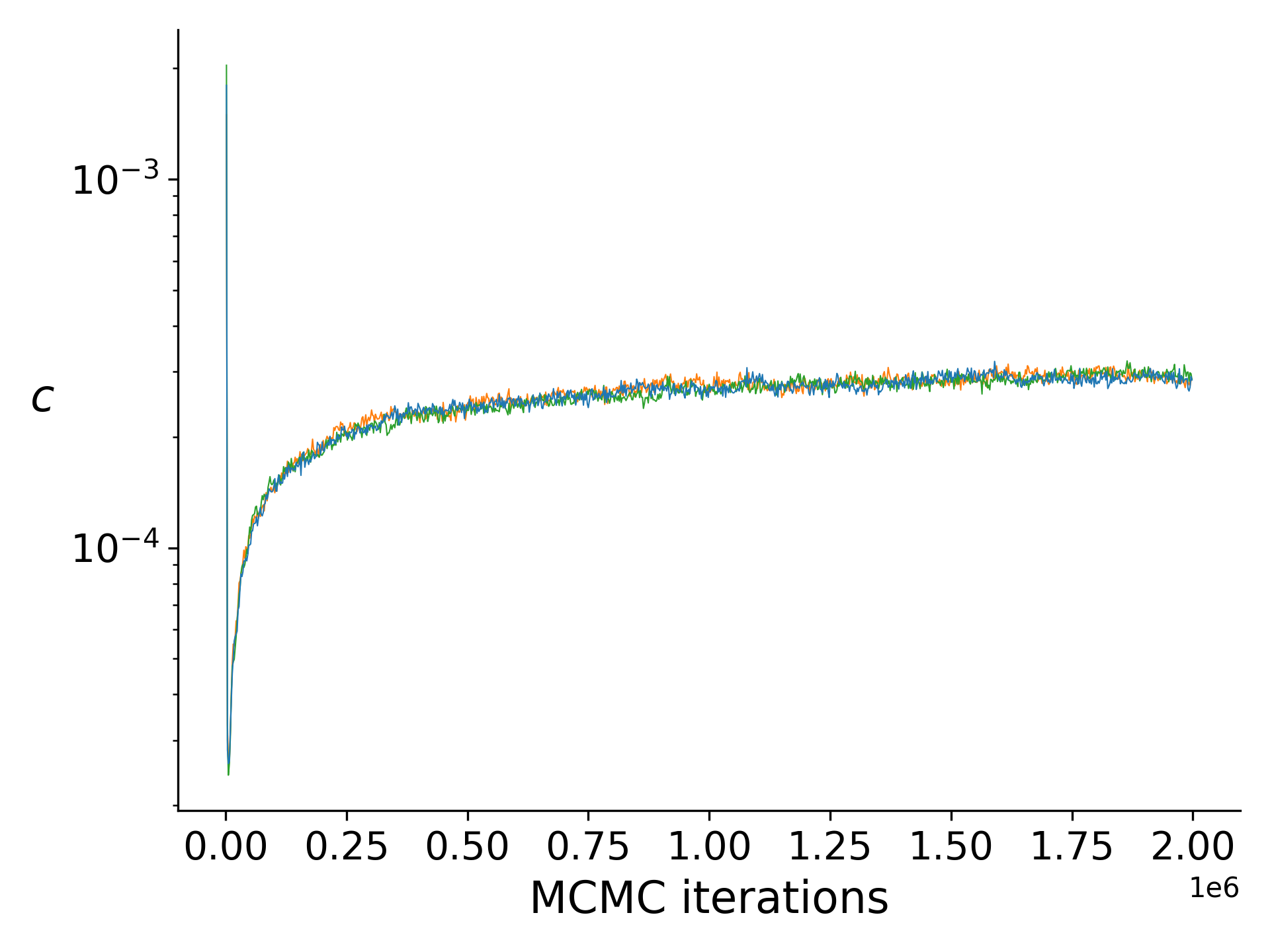}
        \caption{}
    \end{subfigure}
    \begin{subfigure}{0.45\textwidth}
        \centering
        \includegraphics[width=\textwidth]{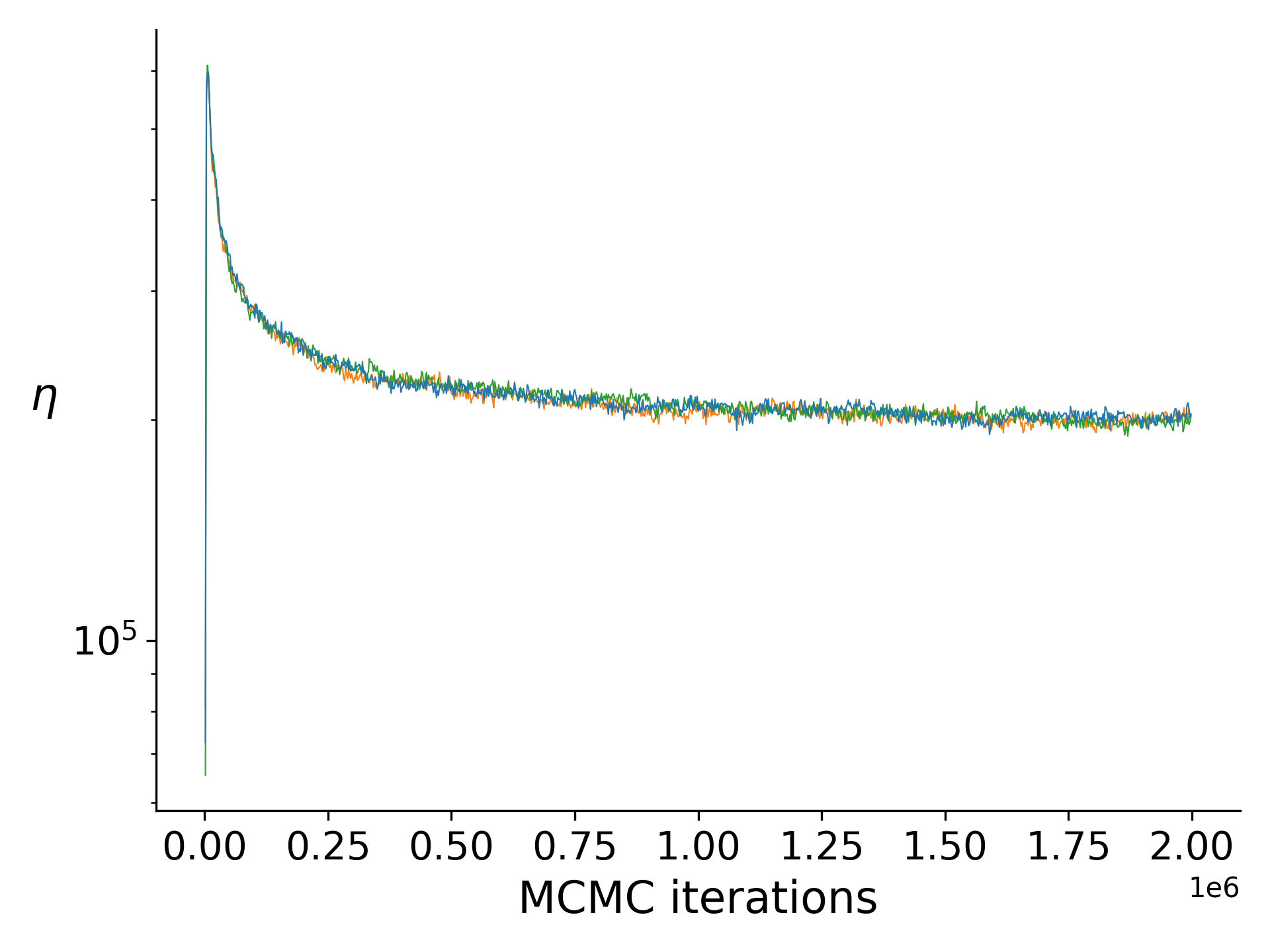}
        \caption{}
    \end{subfigure}
        \begin{subfigure}{0.45\textwidth}
        \centering
        \includegraphics[width=\textwidth]{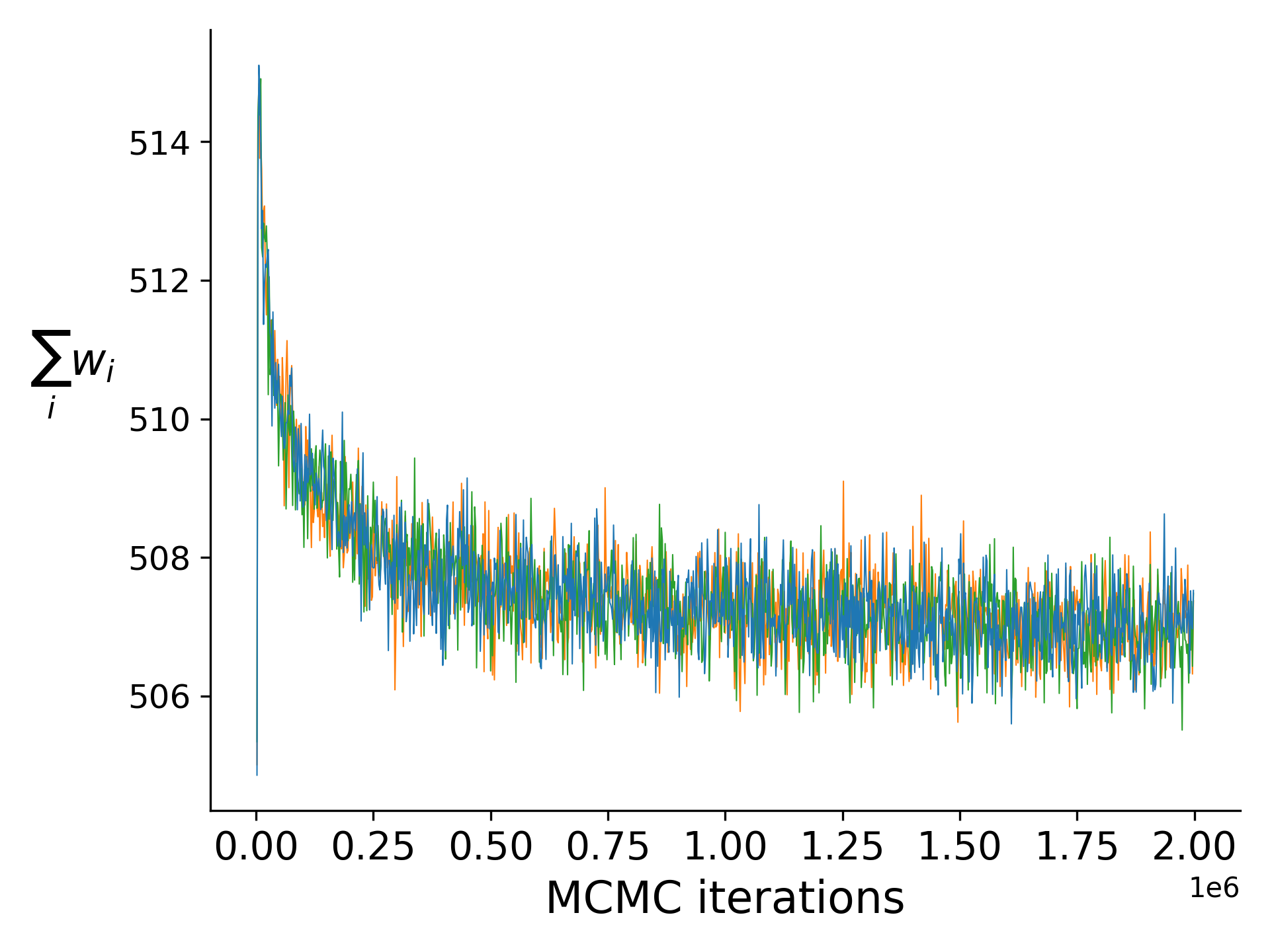}
        \caption{}
    \end{subfigure}
       \caption{MCMC traceplot of parameters (a) $\beta$, (b), c, (c), $\eta$ and (d) the total sum of the weights $w$ for the Flickr subgraph}
    \label{fig:traceplotFlickr}
\end{figure}

\begin{figure}[t]
    \centering
    \begin{subfigure}{0.45\textwidth}
        \centering
        \includegraphics[width=\textwidth]{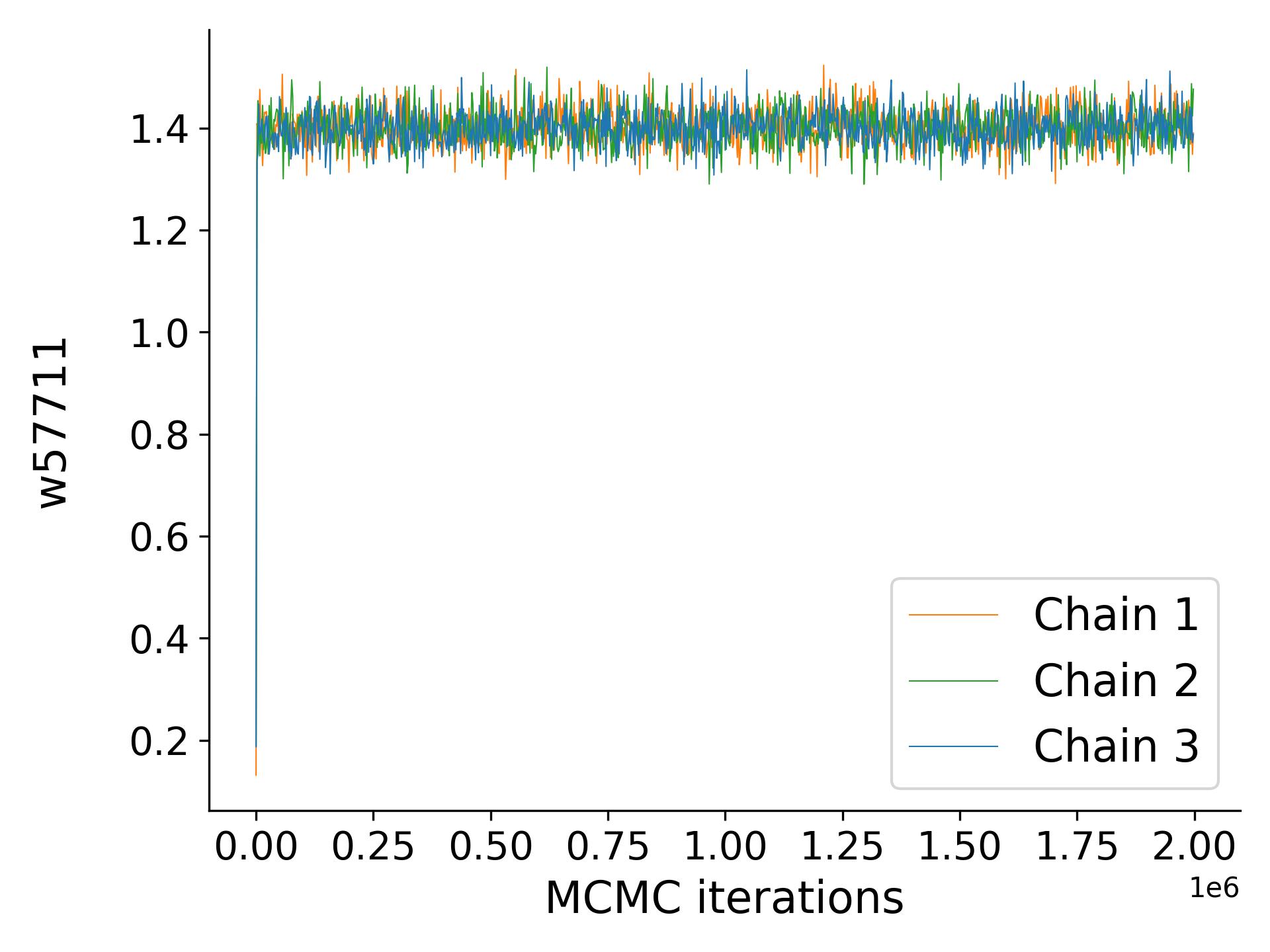}
        \caption{}
    \end{subfigure}
     \begin{subfigure}{0.45\textwidth}
        \centering
        \includegraphics[width=\textwidth]{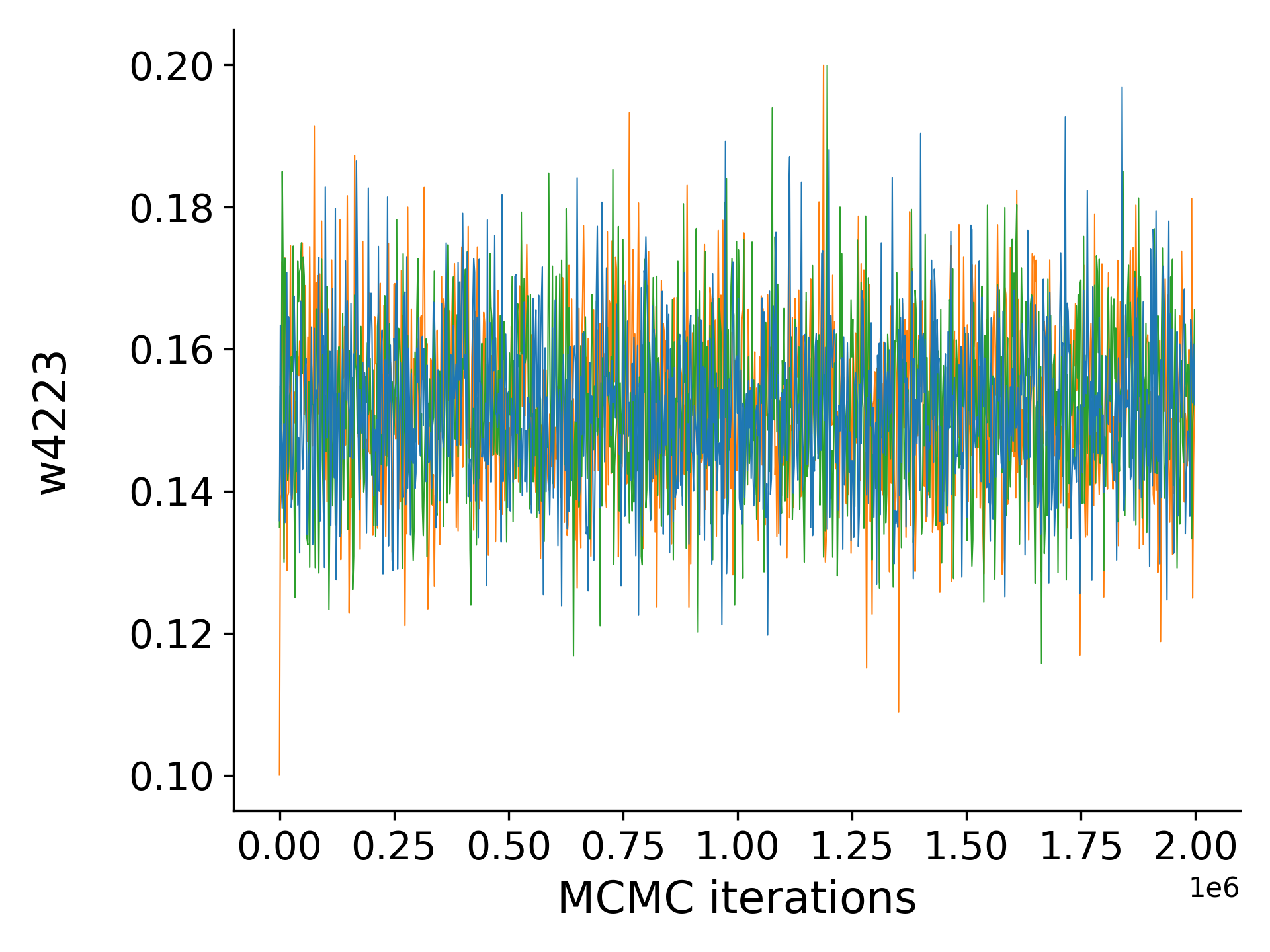}
        \caption{}
    \end{subfigure}

    \begin{subfigure}{0.45\textwidth}
        \centering
        \includegraphics[width=\textwidth]{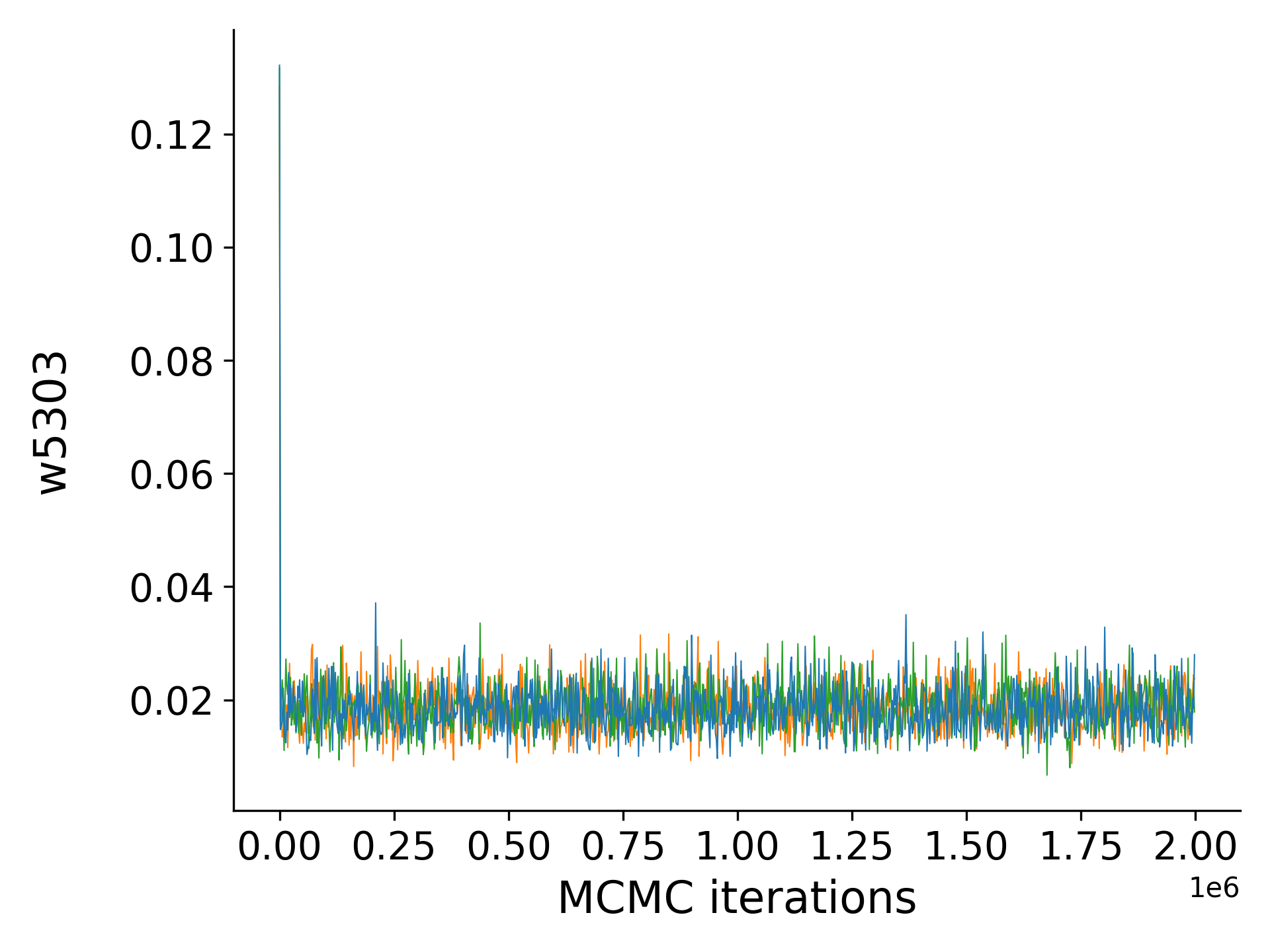}
        \caption{}
    \end{subfigure}
    \begin{subfigure}{0.45\textwidth}
        \centering
        \includegraphics[width=\textwidth]{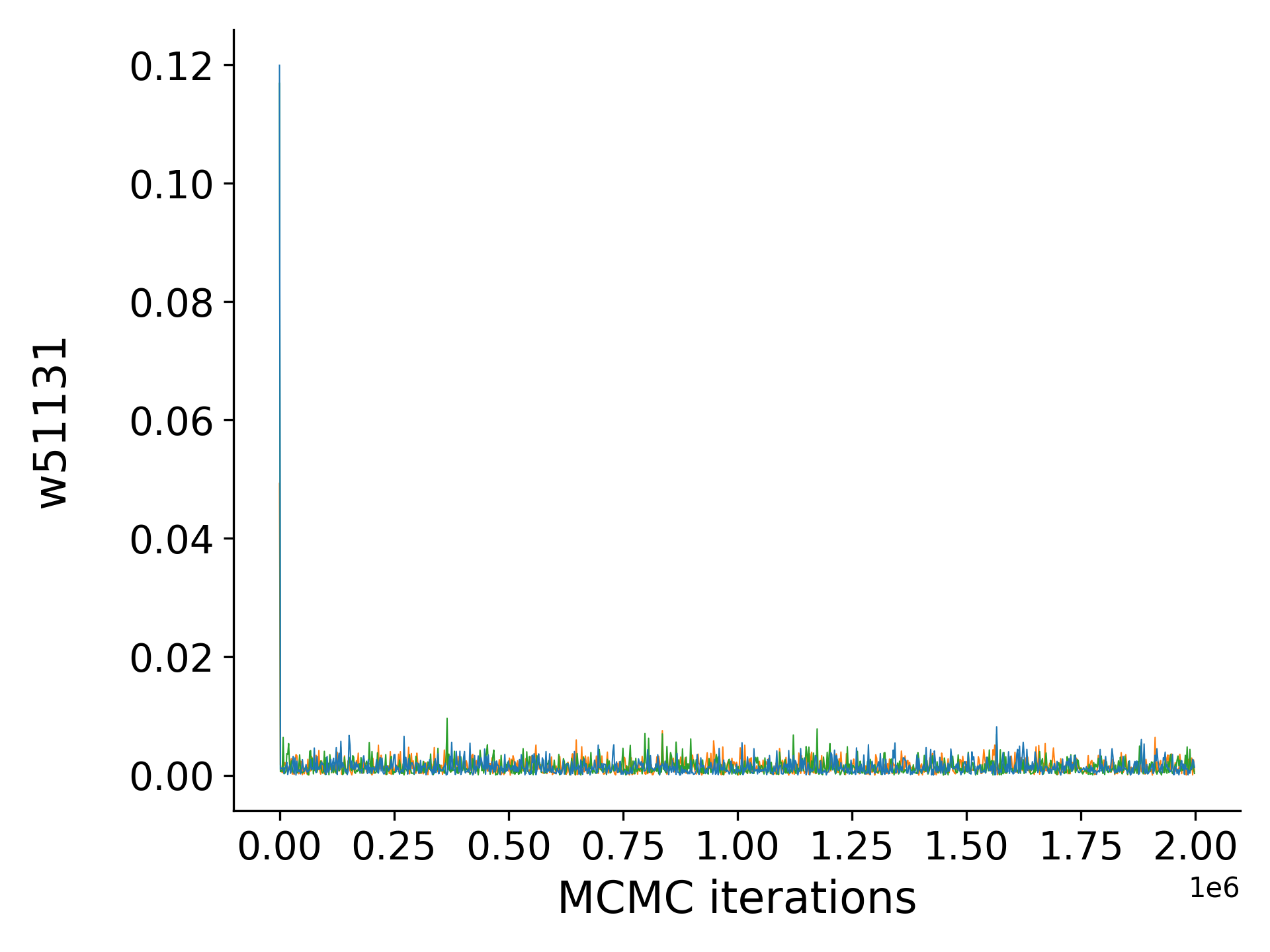}
        \caption{}
    \end{subfigure}
       \caption{MCMC traceplot of four weights for the Flickr subgraph. The degree of the corresponding node are (a) 1445, (b) 166, (c) 21 and (d) 2.}
    \label{fig:traceploweightsFlickr}
\end{figure}

\begin{figure}[t]
    \centering
    \begin{subfigure}{0.45\textwidth}
        \centering
        \includegraphics[width=\textwidth]{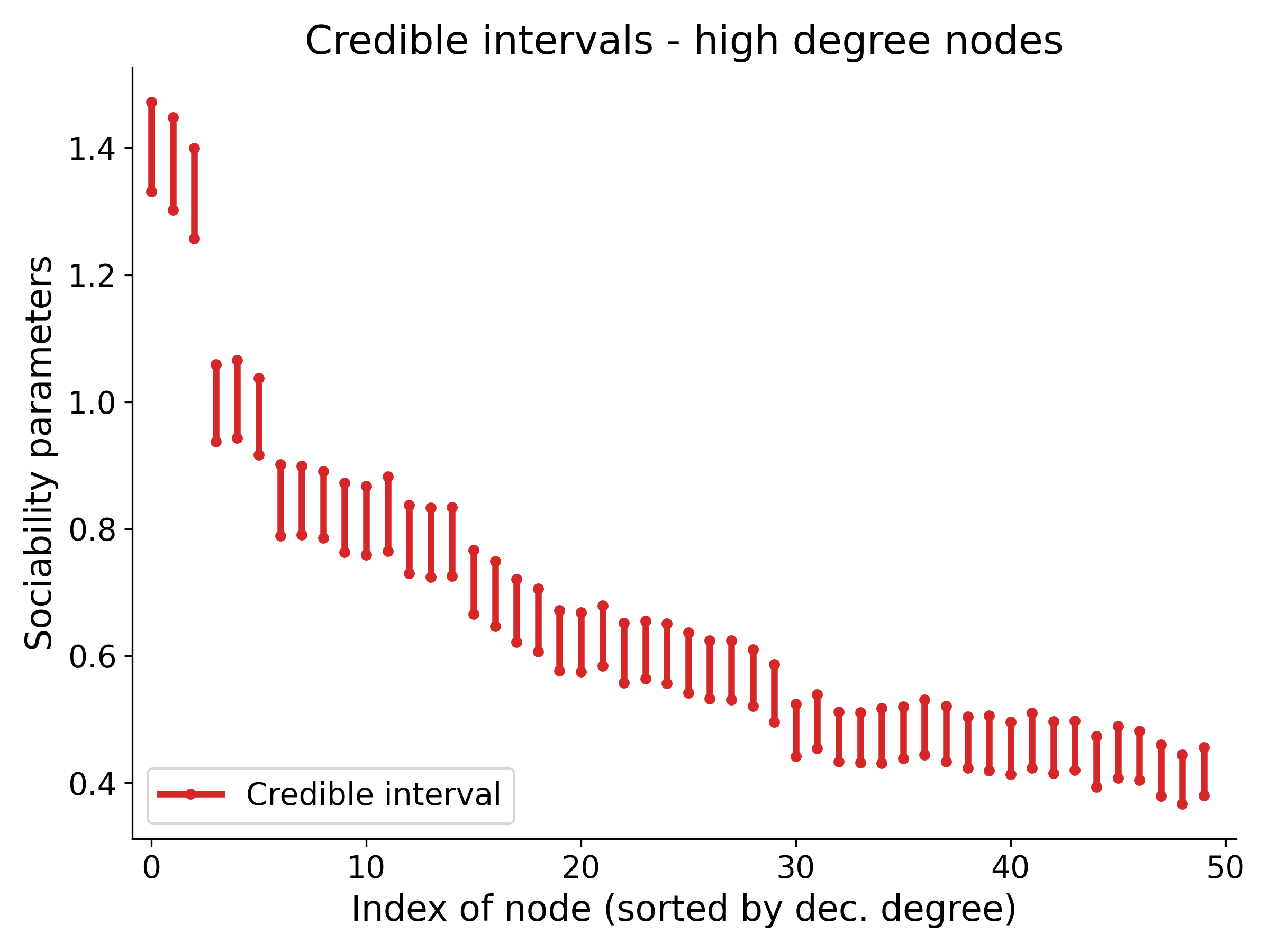}
        \caption{}
    \end{subfigure}
    \begin{subfigure}{0.45\textwidth}
        \centering
        \includegraphics[width=\textwidth]{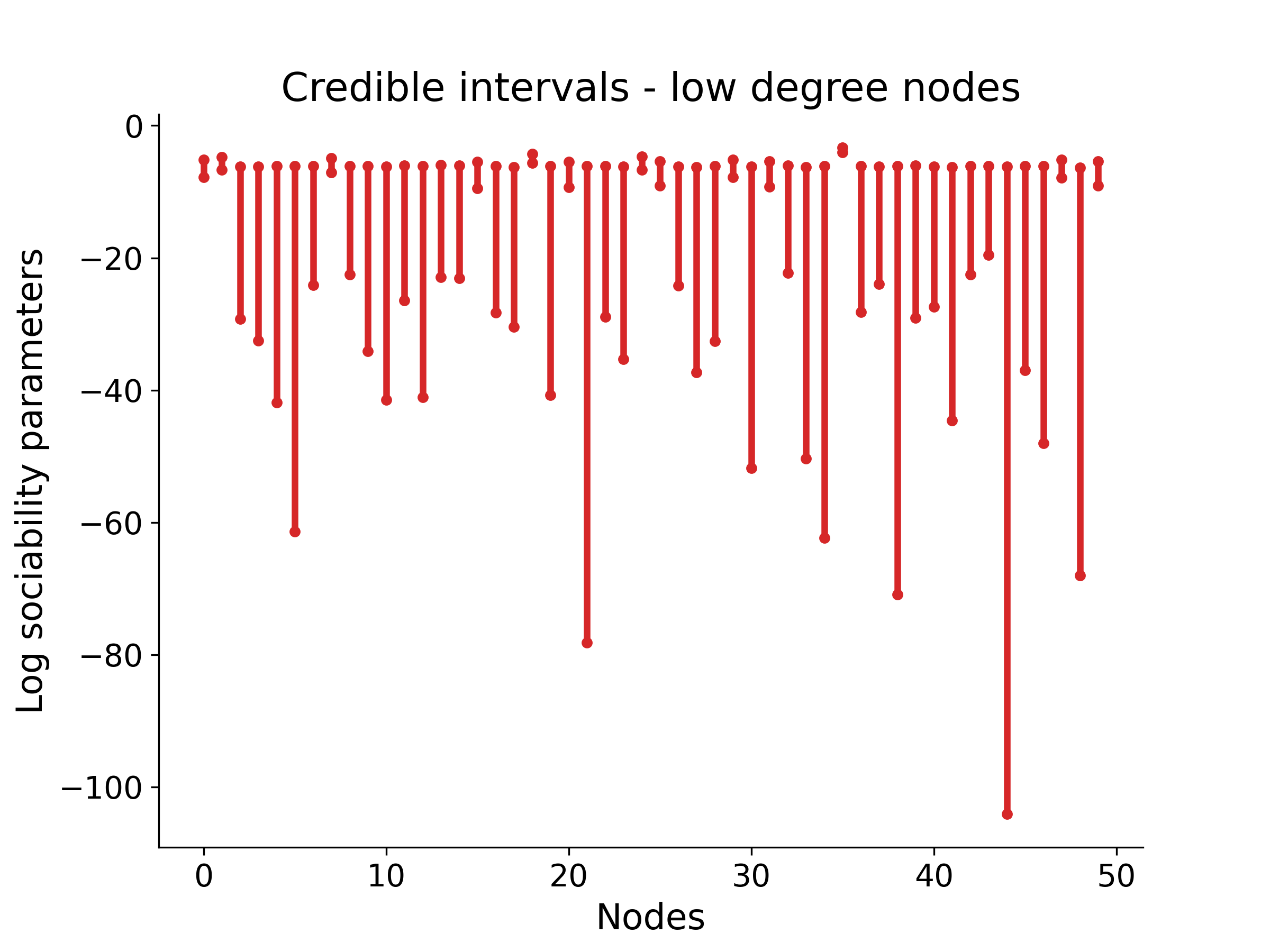}
        \caption{}
    \end{subfigure}
        \caption{95\% posterior intervals of (a) the sociability parameters $w_i$ of the 50 nodes with highest degree and (b) the log-sociability parameters  $\log(w_i)$ of the 50 nodes with lowest degree, for the Flickr subgraph.}
            \label{fig:credibleintervalFlickr}
\end{figure}

\begin{figure}[ht]
    \centering
    \begin{subfigure}{0.45\textwidth}
        \centering
        \includegraphics[width=\textwidth]{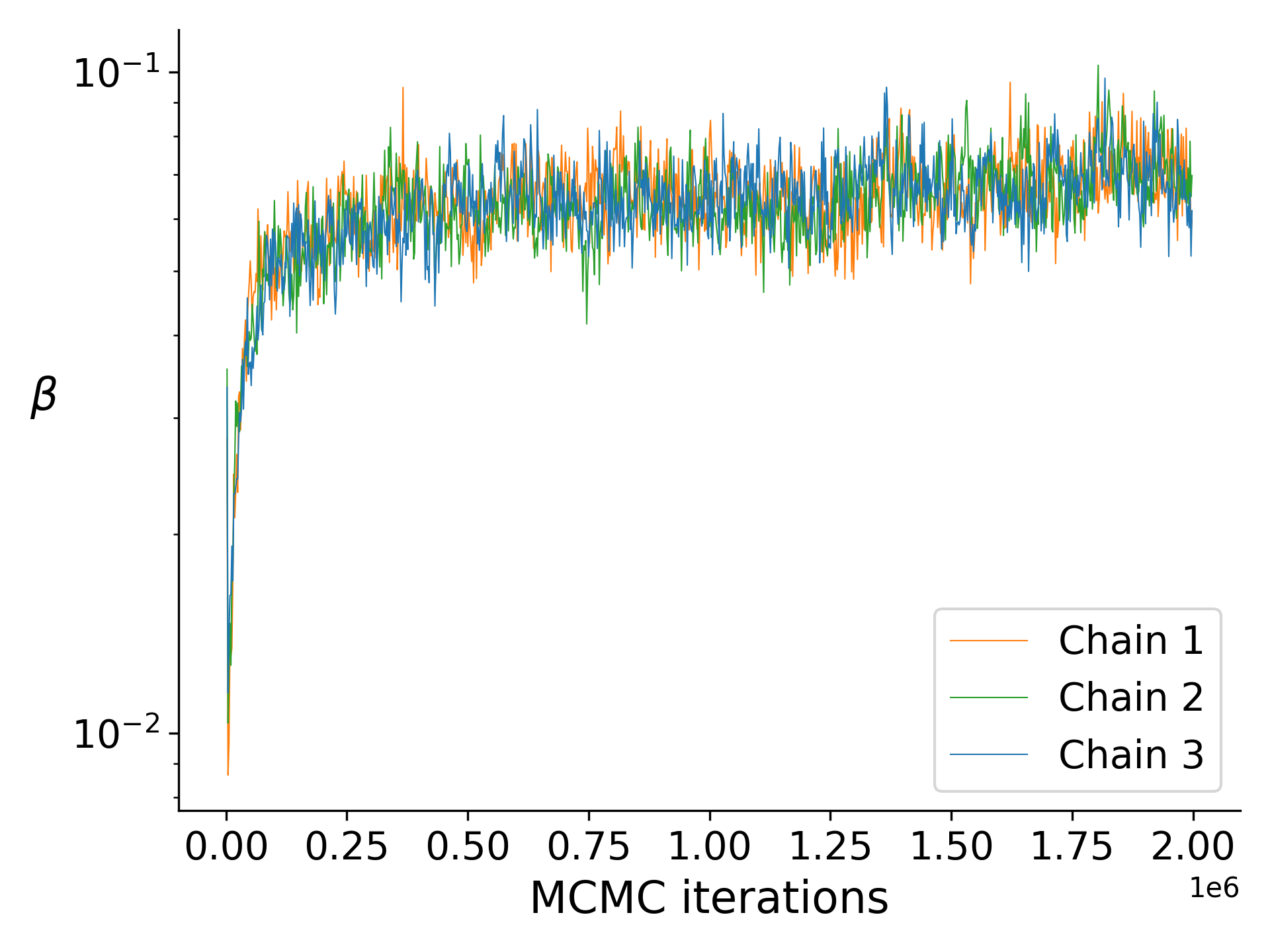}
        \caption{}
    \end{subfigure}
    \begin{subfigure}{0.45\textwidth}
        \centering
        \includegraphics[width=\textwidth]{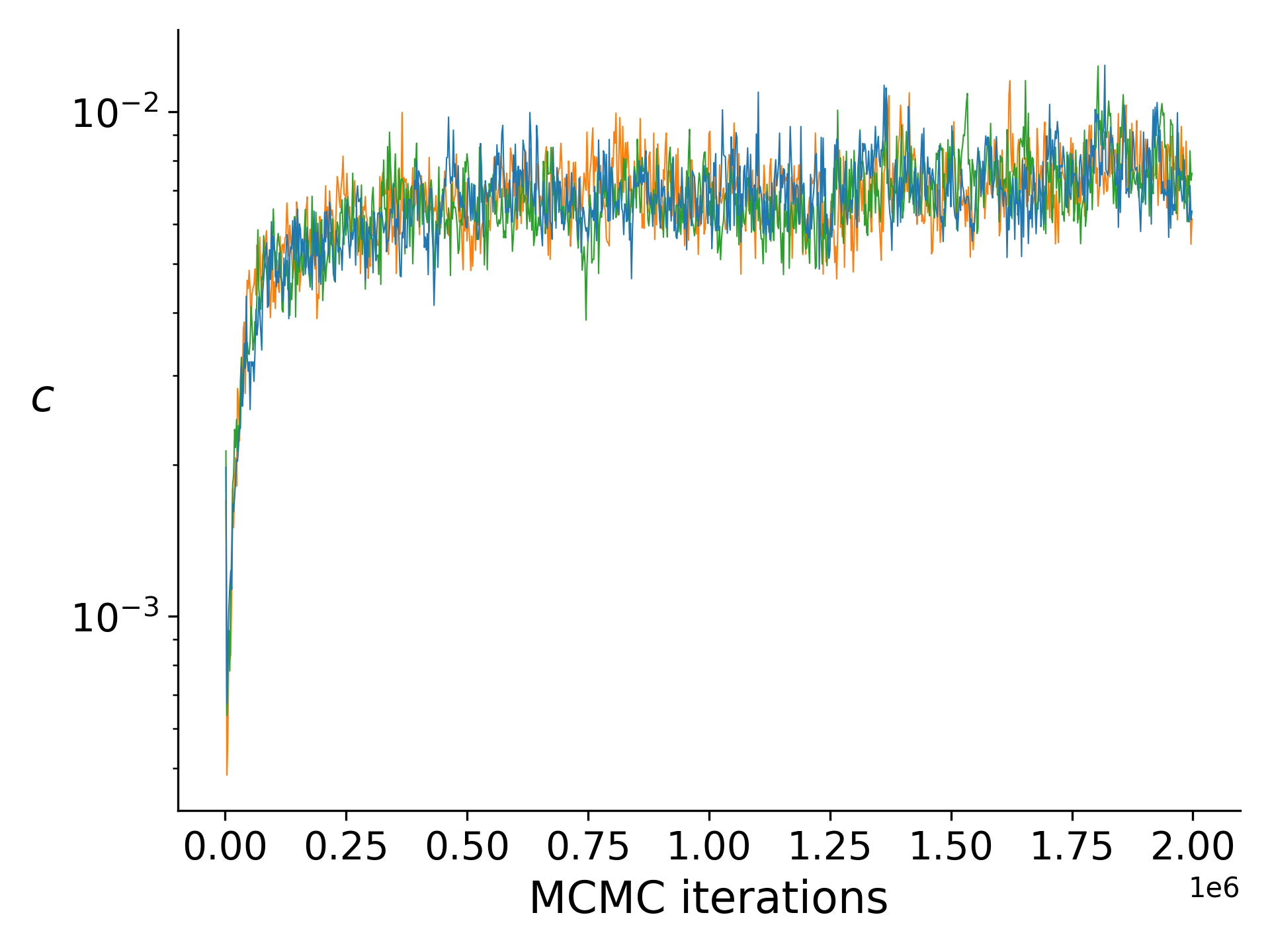}
        \caption{}
    \end{subfigure}
    \begin{subfigure}{0.45\textwidth}
        \centering
        \includegraphics[width=\textwidth]{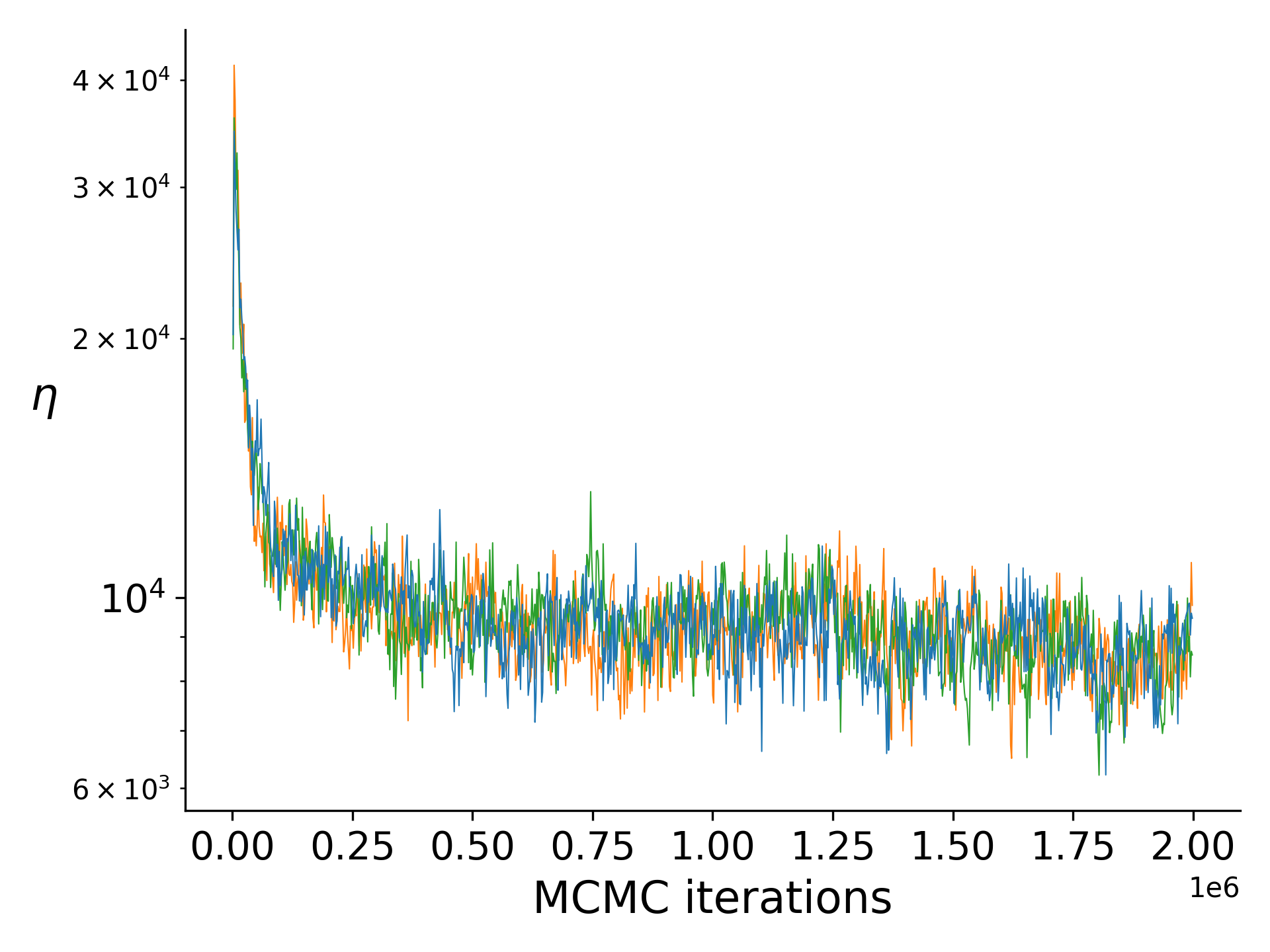}
        \caption{}
    \end{subfigure}
        \begin{subfigure}{0.45\textwidth}
        \centering
        \includegraphics[width=\textwidth]{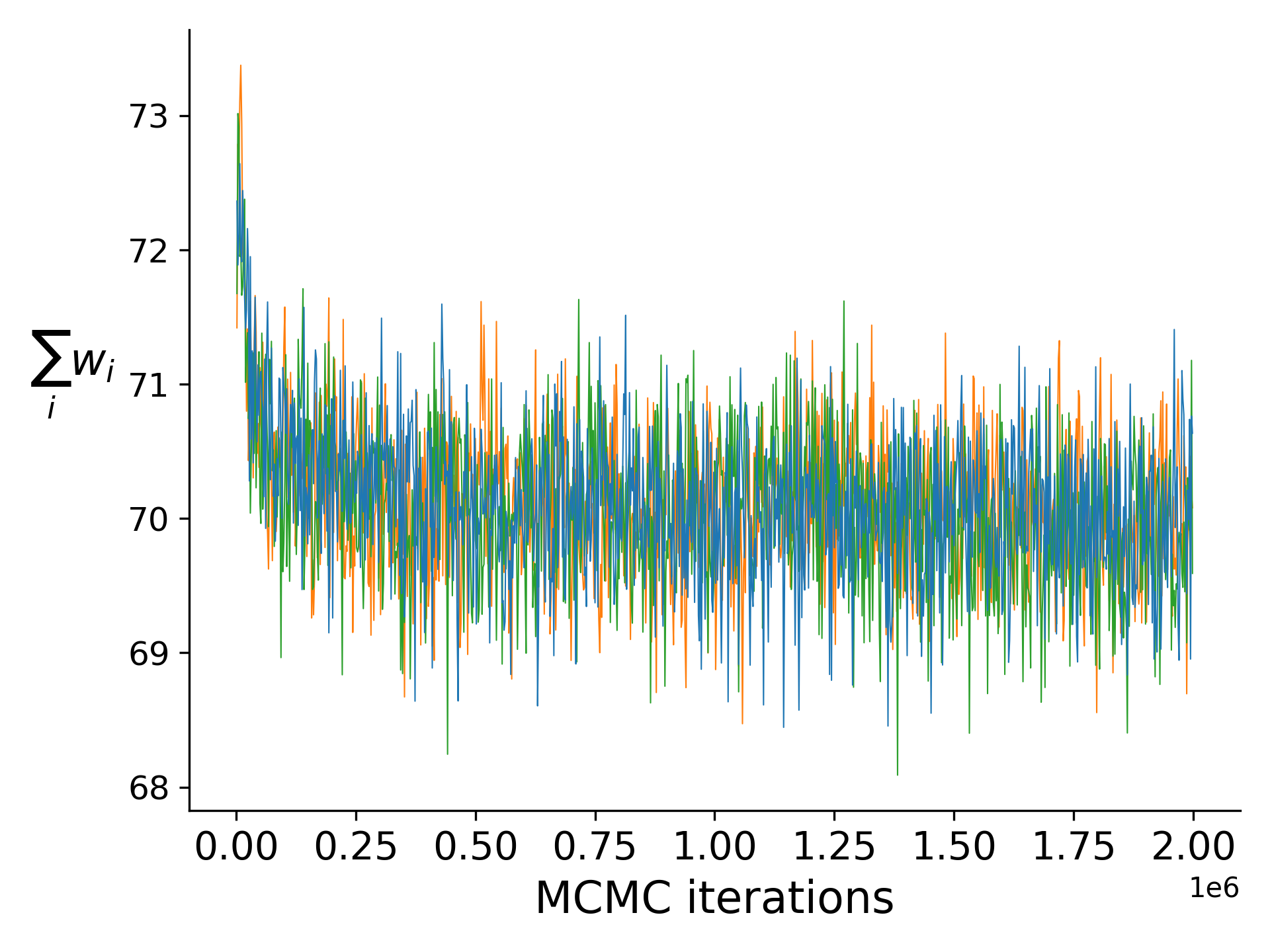}
        \caption{}
    \end{subfigure}
       \caption{MCMC traceplot of parameters (a) $\beta$, (b), c, (c), $\eta$ and (d) the total sum of the weights $w$ for the Douban subgraph}
    \label{fig:traceplotDouban}
\end{figure}

\begin{figure}[t]
    \centering
    \begin{subfigure}{0.45\textwidth}
        \centering
        \includegraphics[width=\textwidth]{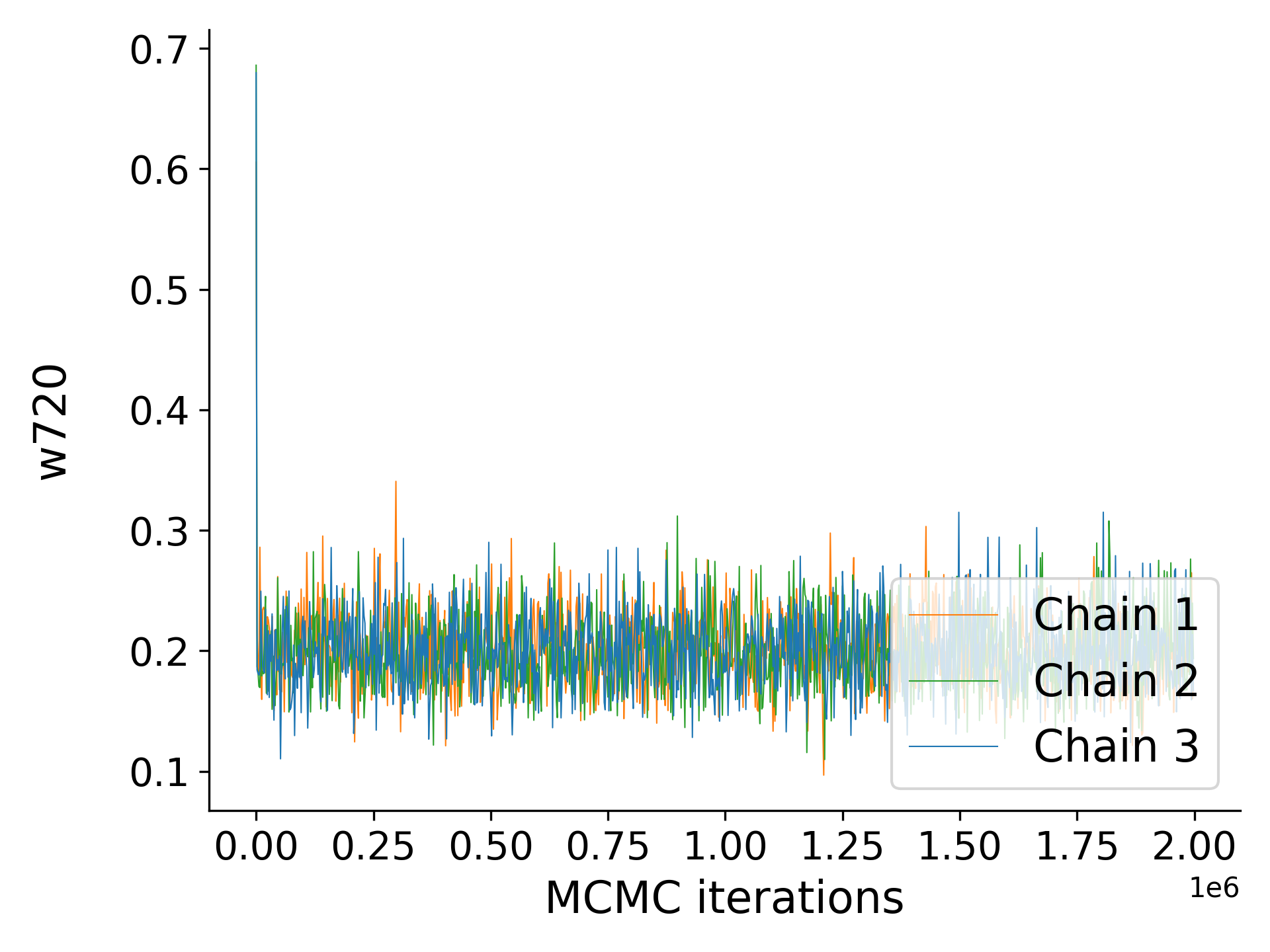}
        \caption{}
    \end{subfigure}
     \begin{subfigure}{0.45\textwidth}
        \centering
        \includegraphics[width=\textwidth]{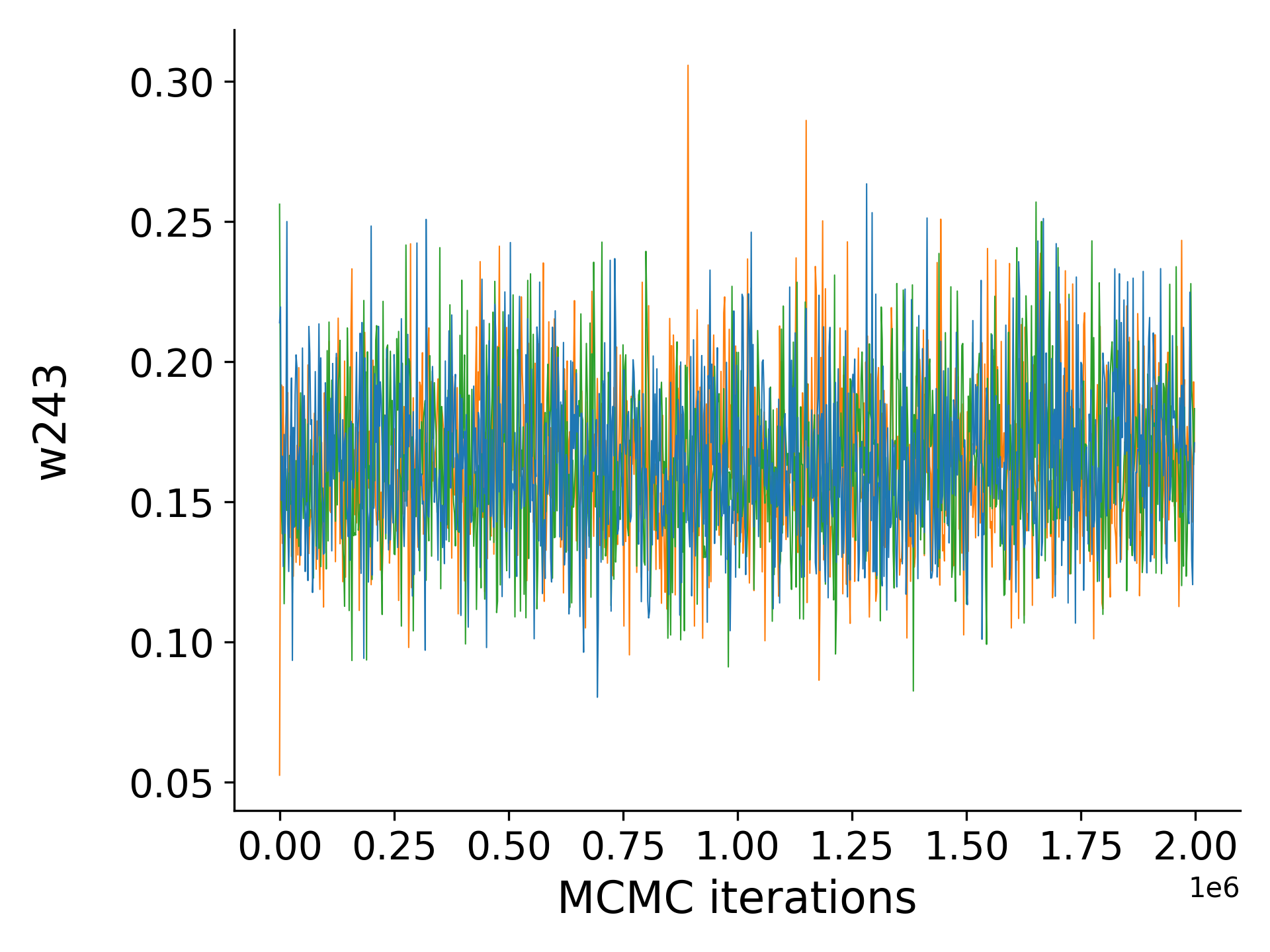}
        \caption{}
    \end{subfigure}

    \begin{subfigure}{0.45\textwidth}
        \centering
        \includegraphics[width=\textwidth]{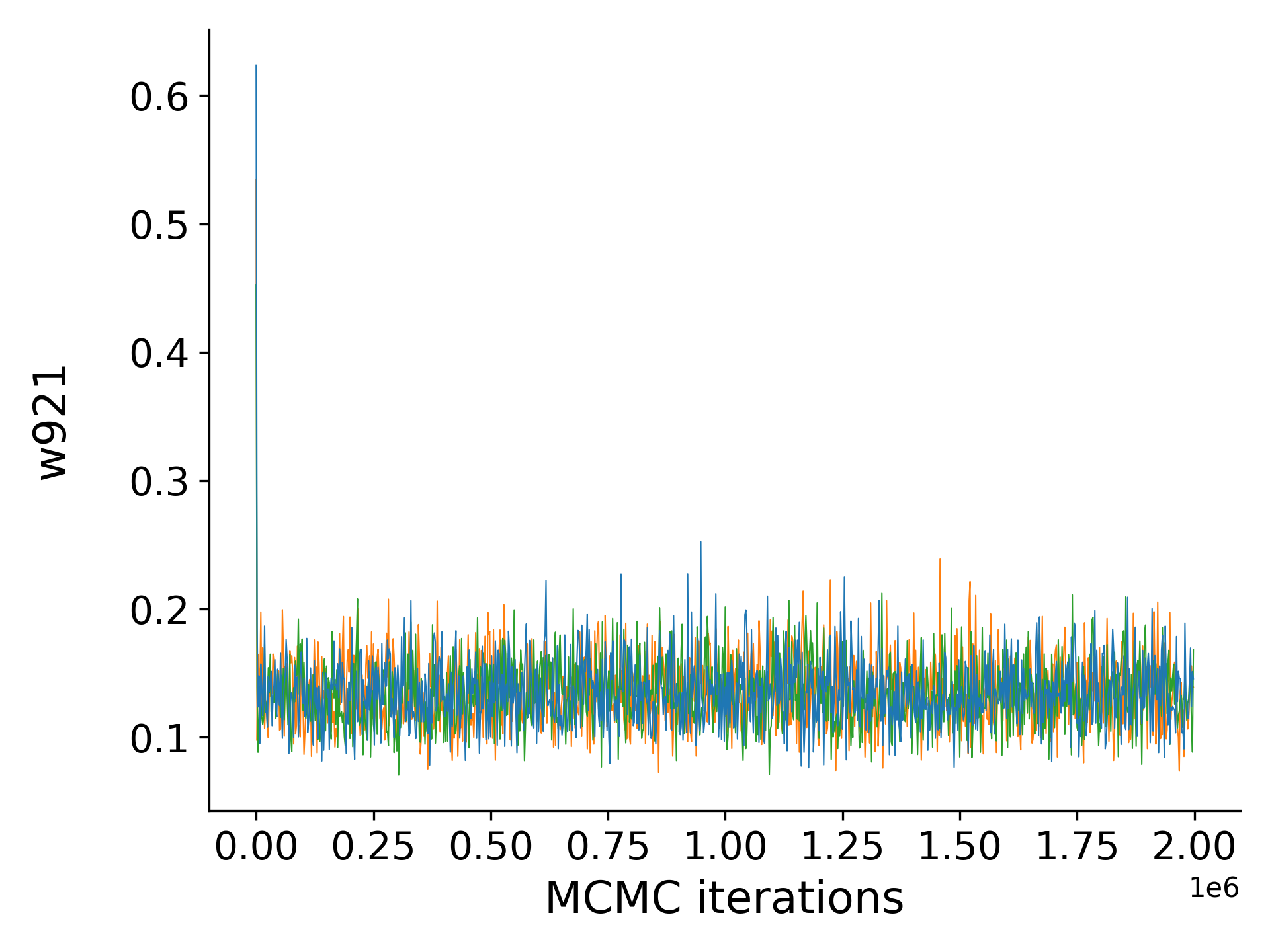}
        \caption{}
    \end{subfigure}
    \begin{subfigure}{0.45\textwidth}
        \centering
        \includegraphics[width=\textwidth]{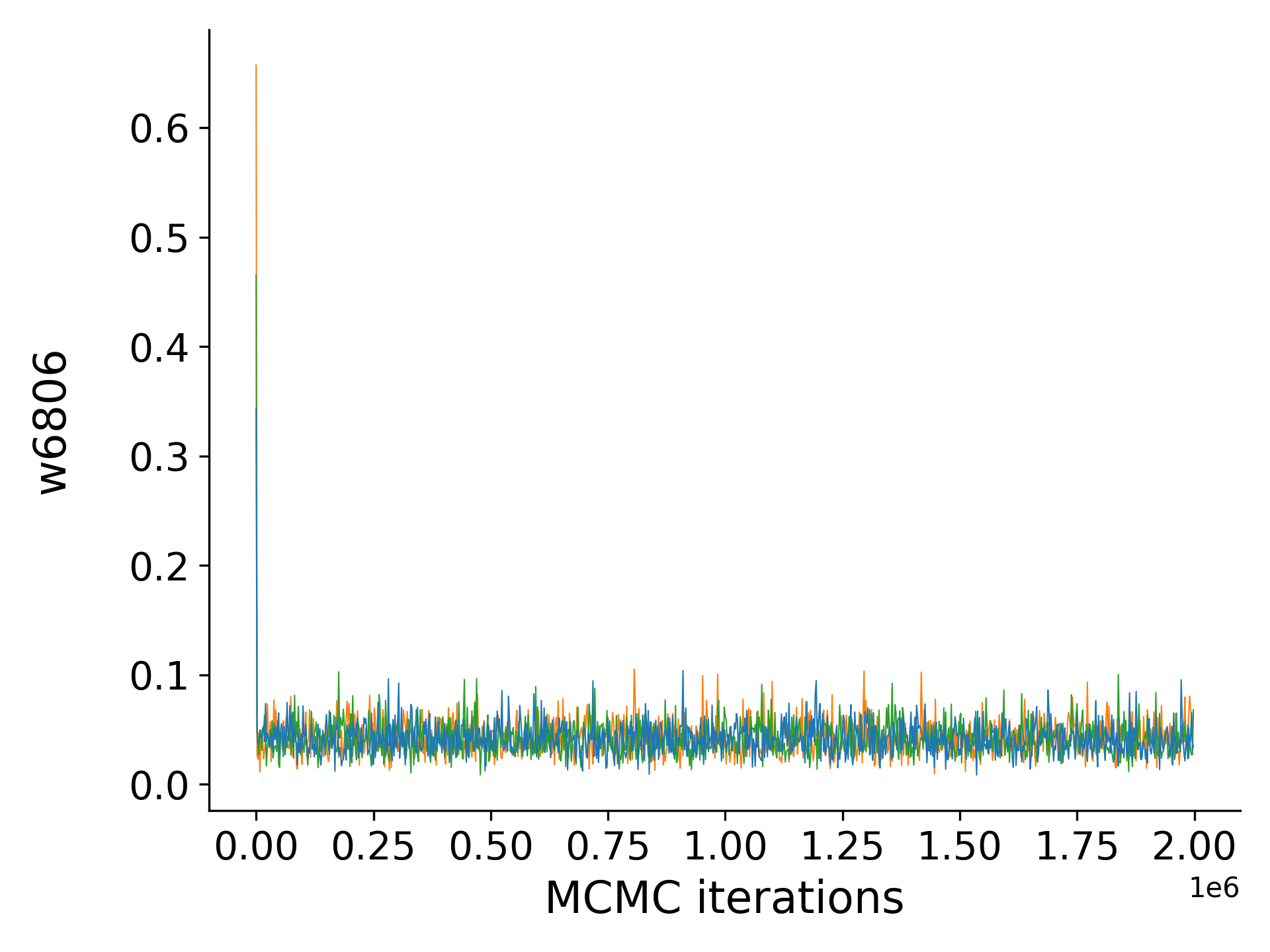}
        \caption{}
    \end{subfigure}
       \caption{MCMC traceplot of four weights for the Douban subgraph. The degree of the corresponding node are (a) 41, (b) 28, (c) 34 and (d) 9.}
    \label{fig:traceploweightsDouban}
\end{figure}

\begin{figure}[t]
    \centering
    \begin{subfigure}{0.45\textwidth}
        \centering
        \includegraphics[width=\textwidth]{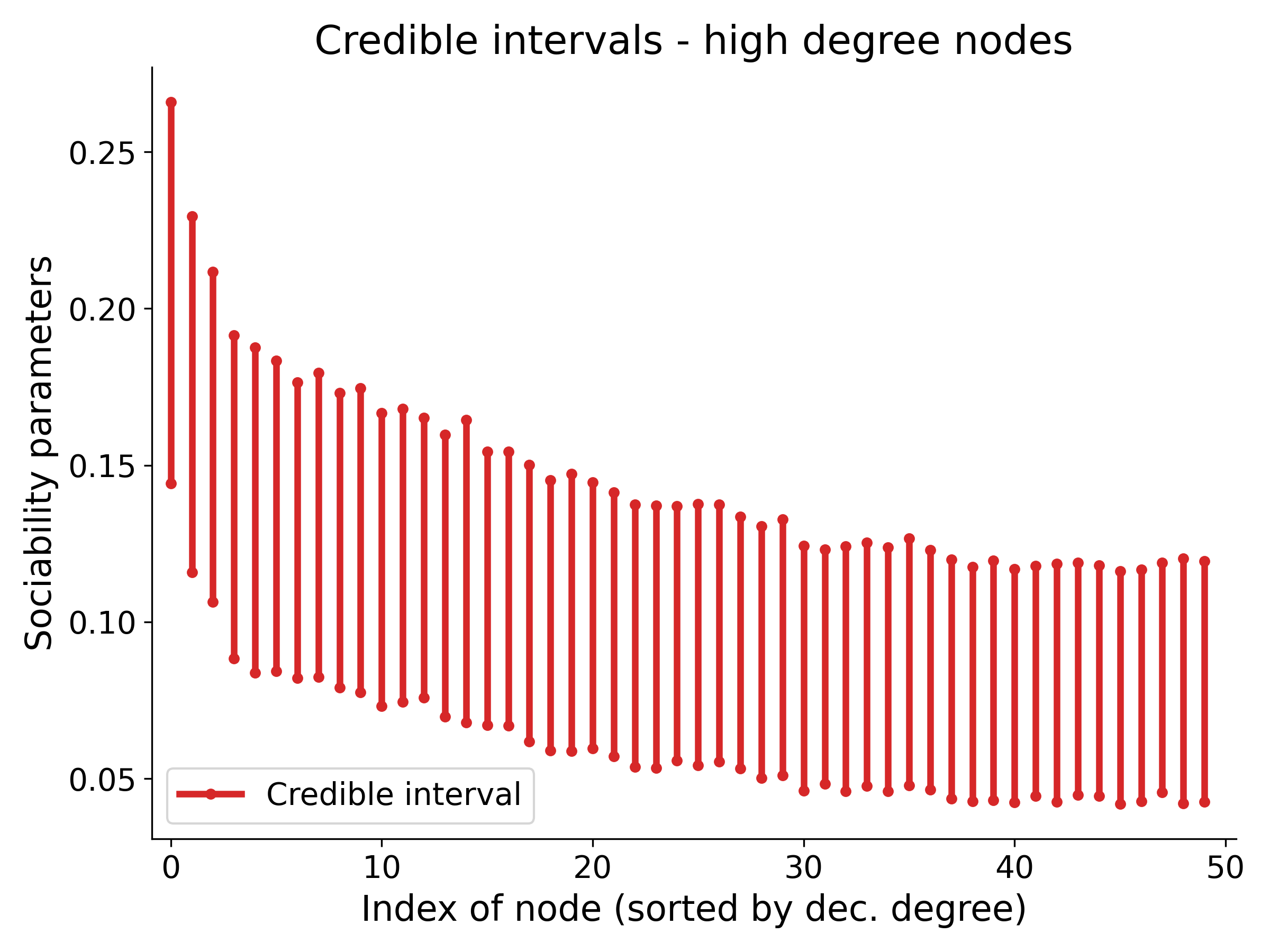}
        \caption{}
    \end{subfigure}
    \begin{subfigure}{0.45\textwidth}
        \centering
        \includegraphics[width=\textwidth]{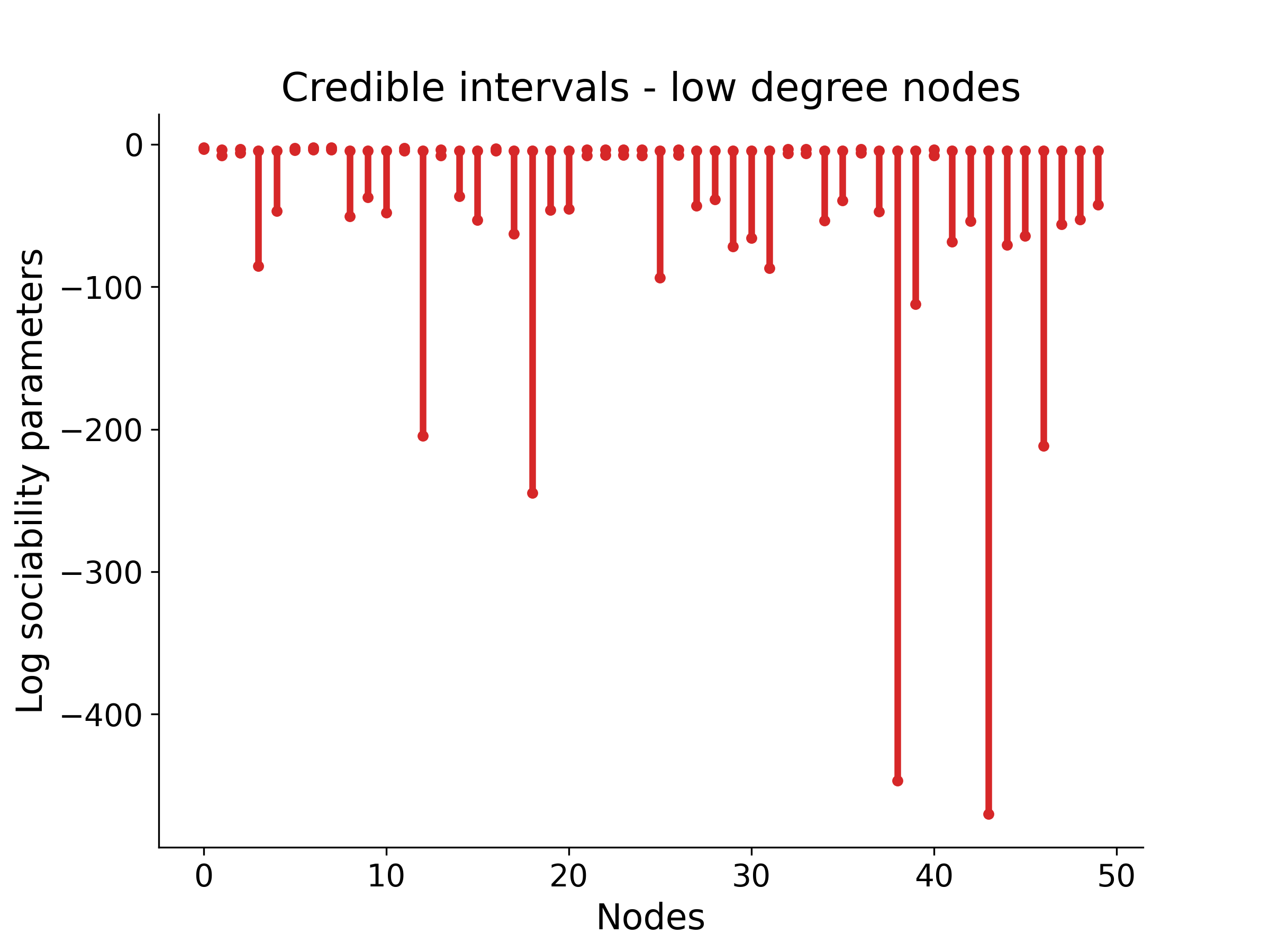}
        \caption{}
    \end{subfigure}
        \caption{95\% posterior intervals of (a) the sociability parameters $w_i$ of the 50 nodes with highest degree and (b) the log-sociability parameters  $\log(w_i)$ of the 50 nodes with lowest degree, for the Douban subgraph.}
            \label{fig:credibleintervalDouban}
\end{figure}

\begin{figure}[ht]
    \centering
    \begin{subfigure}{0.45\textwidth}
        \centering
        \includegraphics[width=\textwidth]{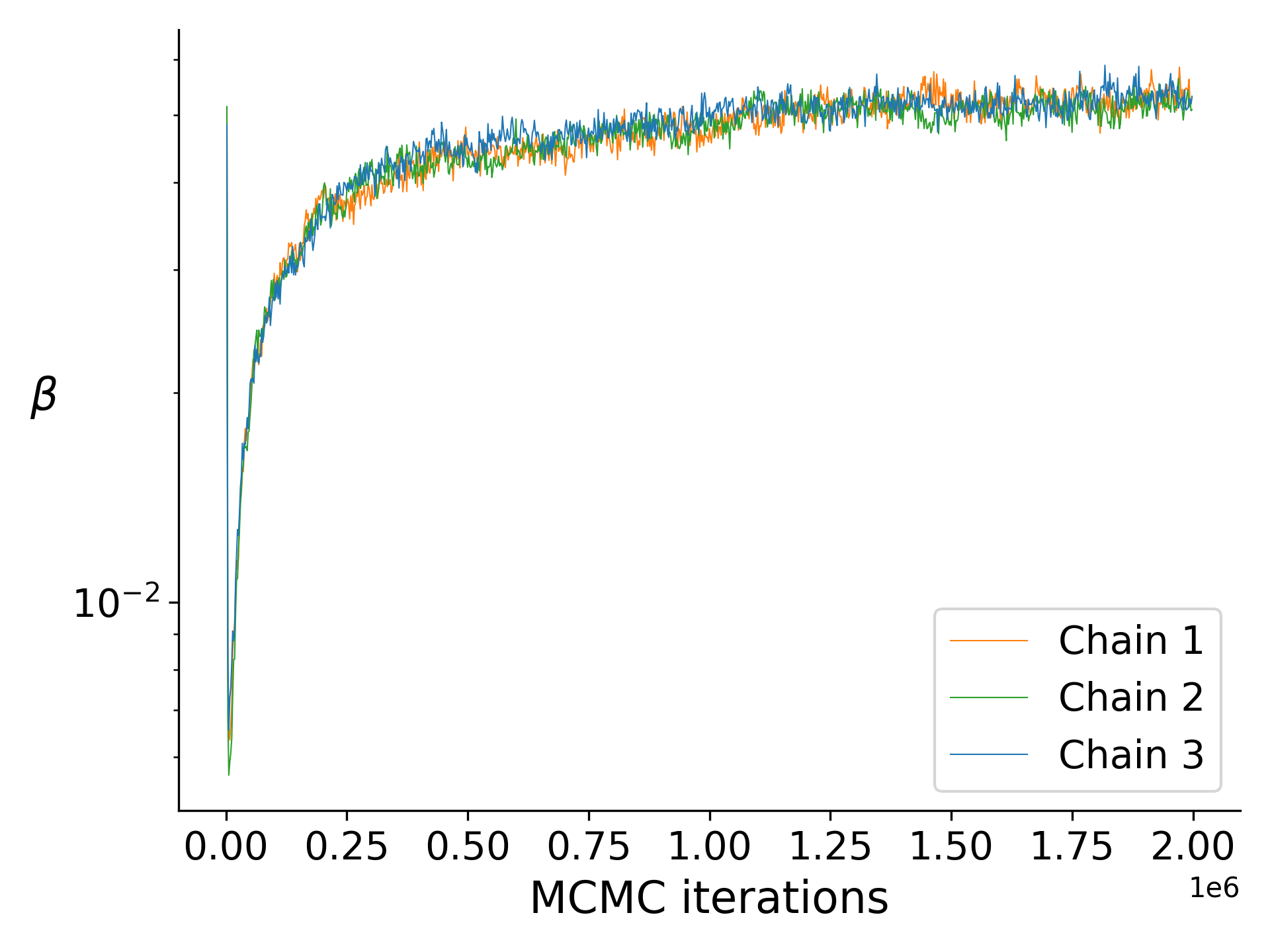}
        \caption{}
    \end{subfigure}
    \begin{subfigure}{0.45\textwidth}
        \centering
        \includegraphics[width=\textwidth]{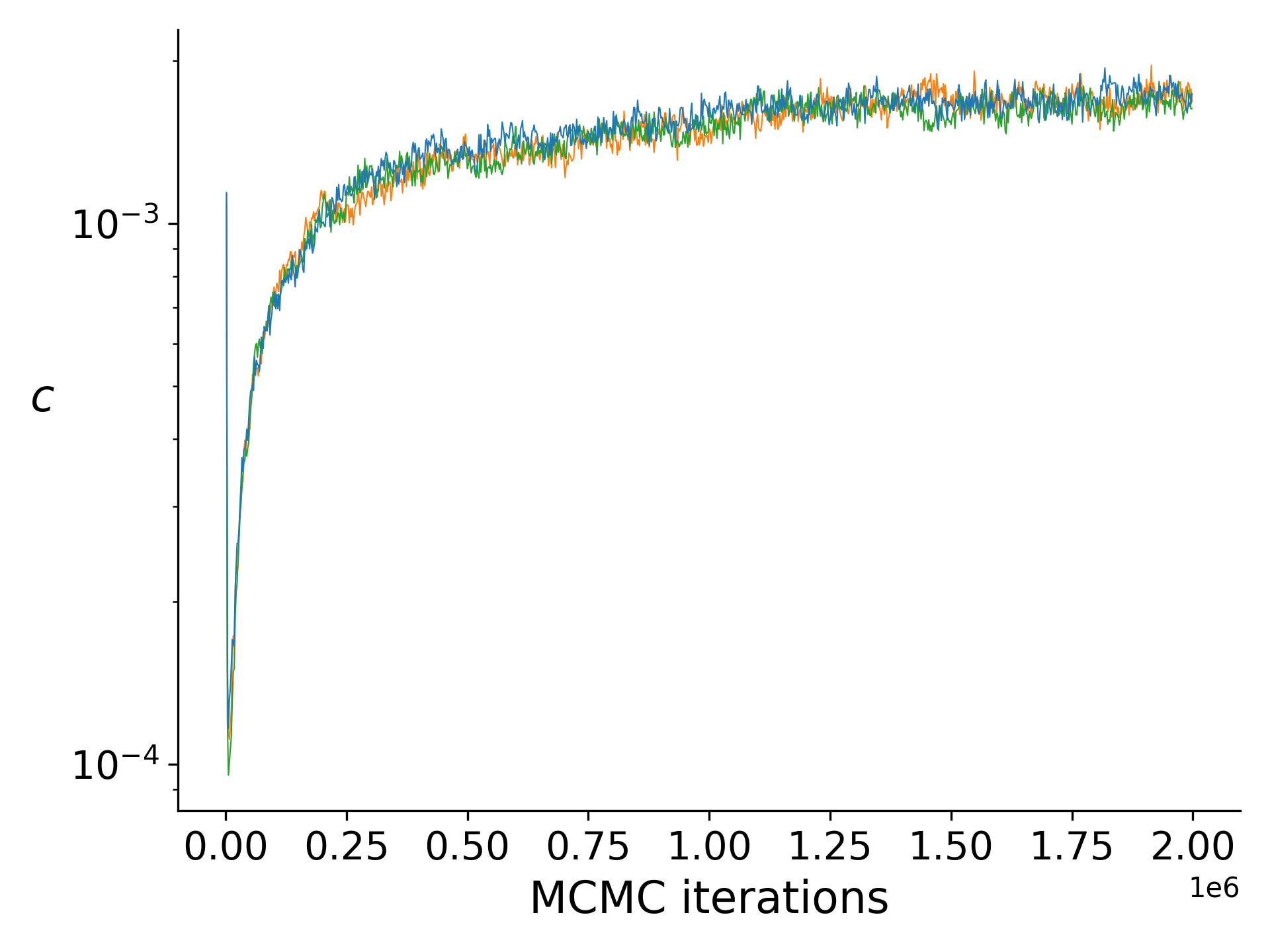}
        \caption{}
    \end{subfigure}
    \begin{subfigure}{0.45\textwidth}
        \centering
        \includegraphics[width=\textwidth]{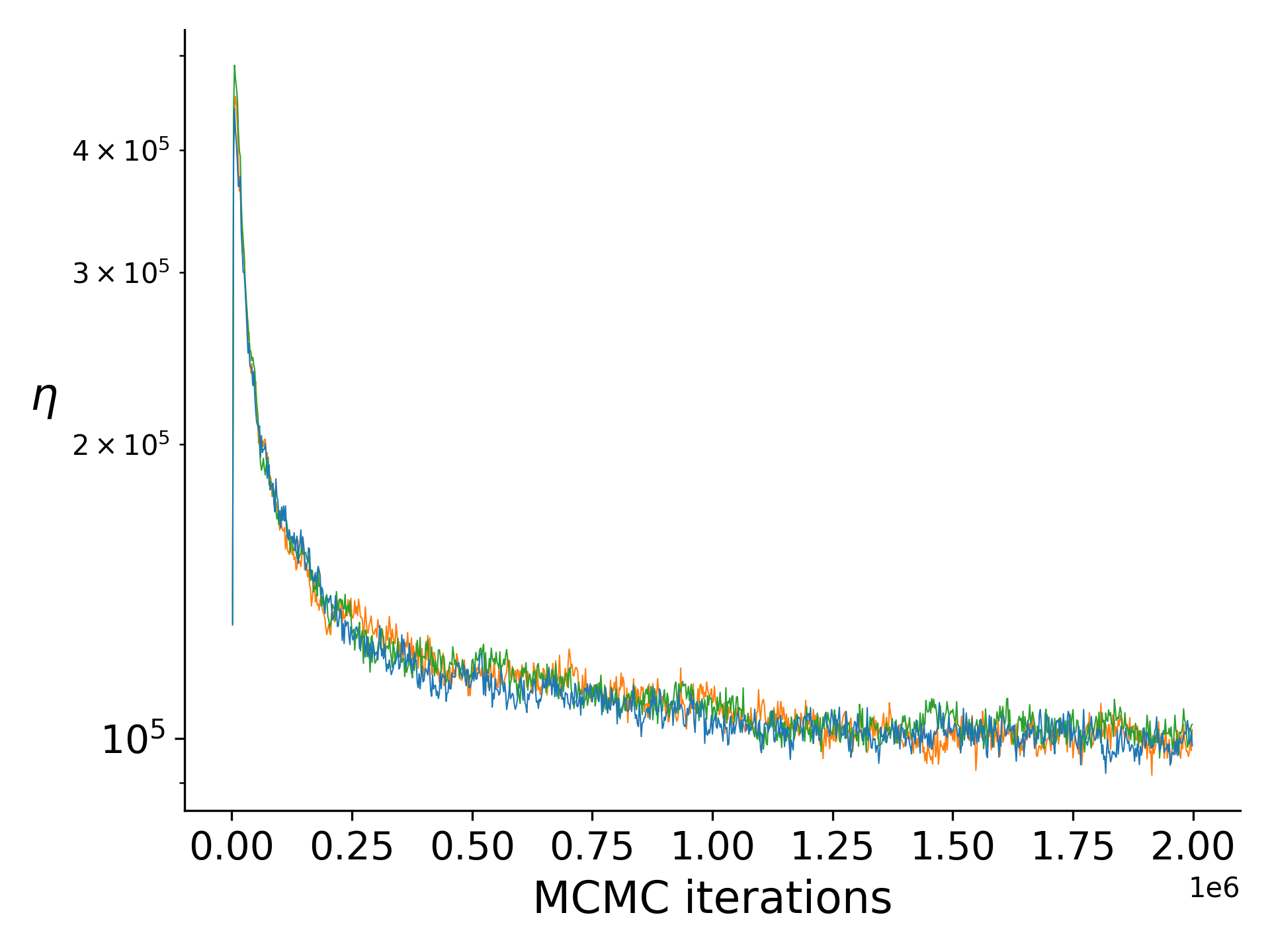}
        \caption{}
    \end{subfigure}
        \begin{subfigure}{0.45\textwidth}
        \centering
        \includegraphics[width=\textwidth]{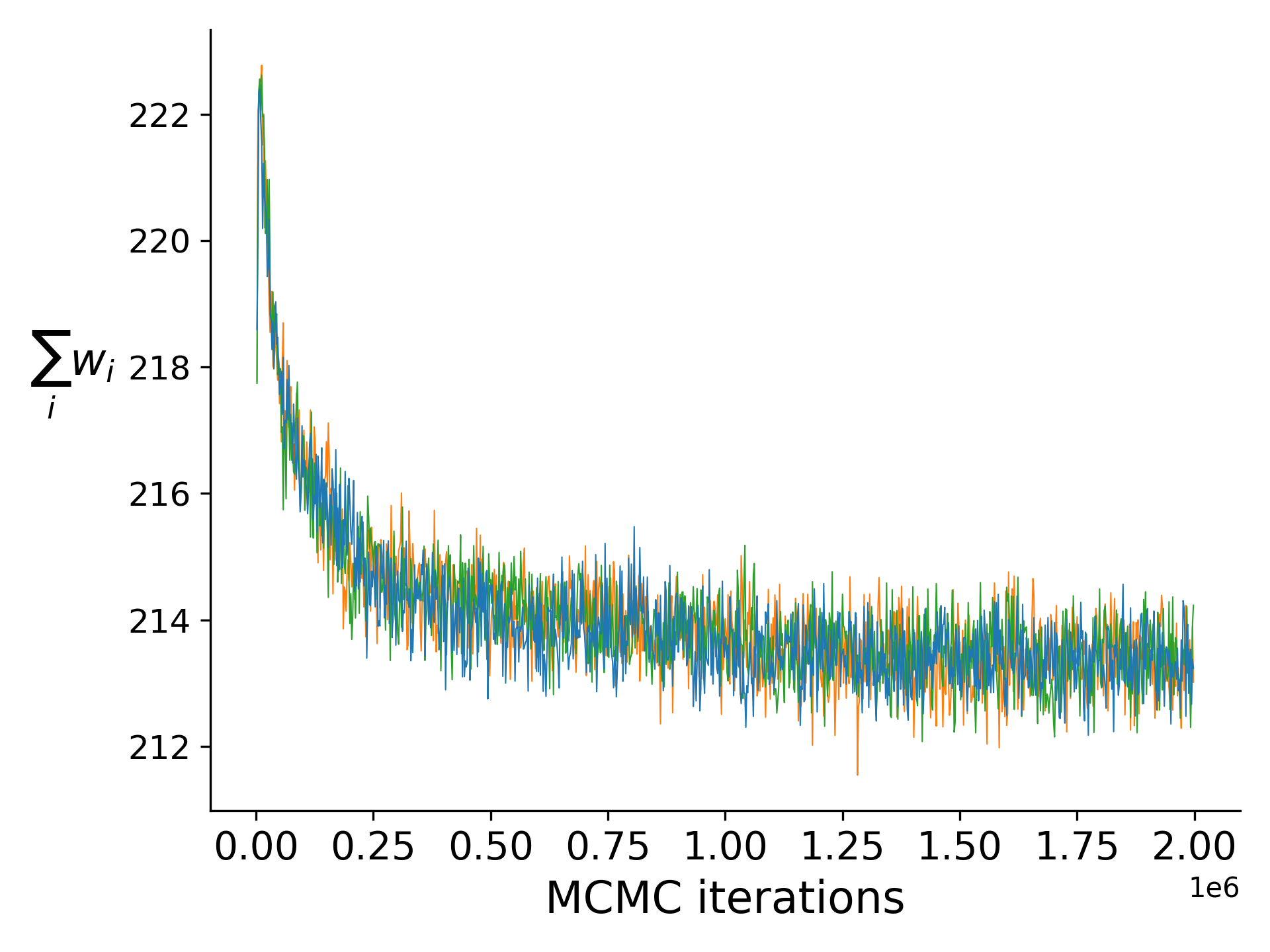}
        \caption{}
    \end{subfigure}
       \caption{MCMC traceplot of parameters (a) $\beta$, (b), c, (c), $\eta$ and (d) the total sum of the weights $w$ for the TwitterCrawl subgraph}
    \label{fig:traceplotTwitterCrawl}
\end{figure}

\begin{figure}[t]
    \centering
    \begin{subfigure}{0.45\textwidth}
        \centering
        \includegraphics[width=\textwidth]{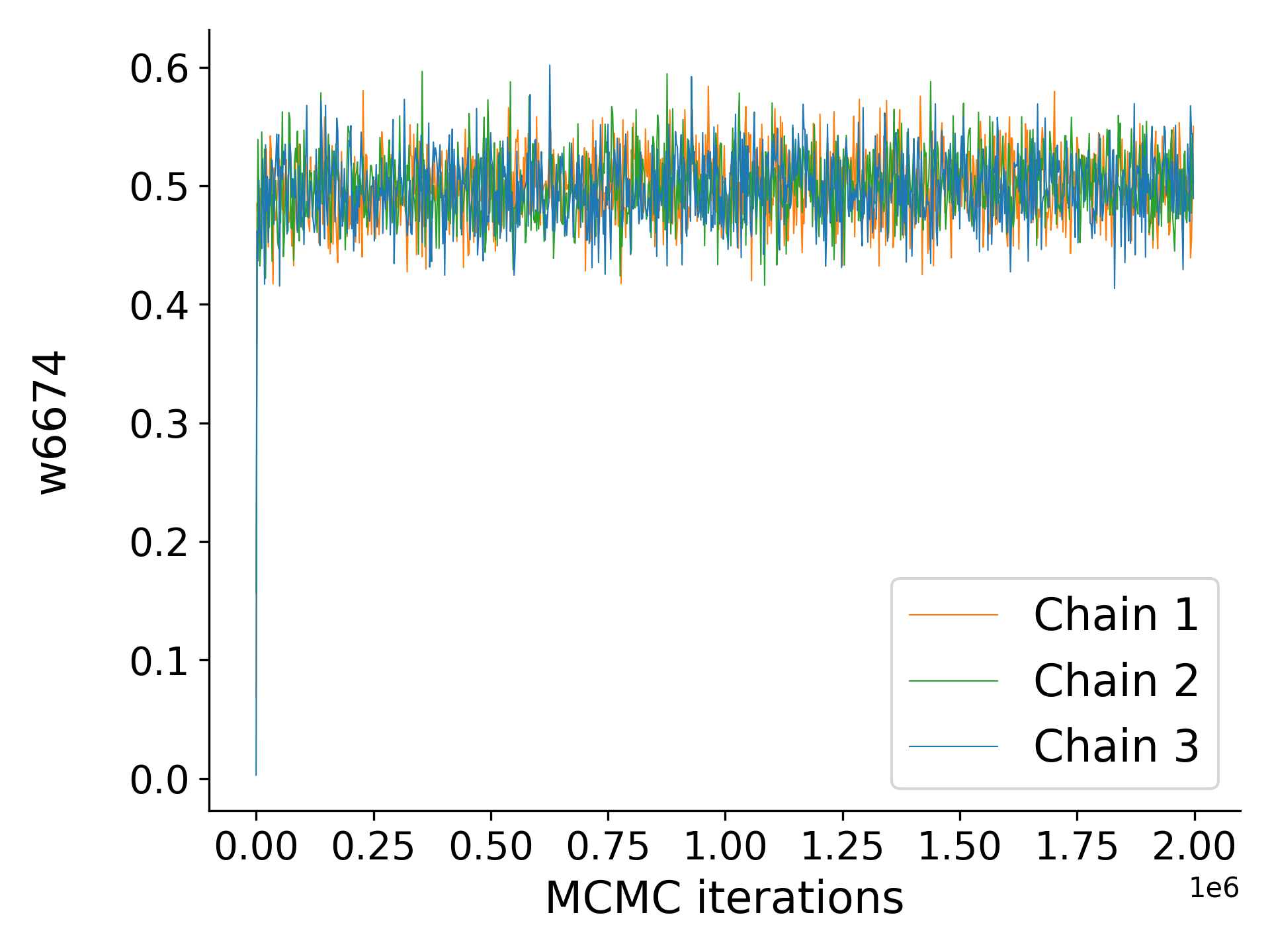}
        \caption{}
    \end{subfigure}
     \begin{subfigure}{0.45\textwidth}
        \centering
        \includegraphics[width=\textwidth]{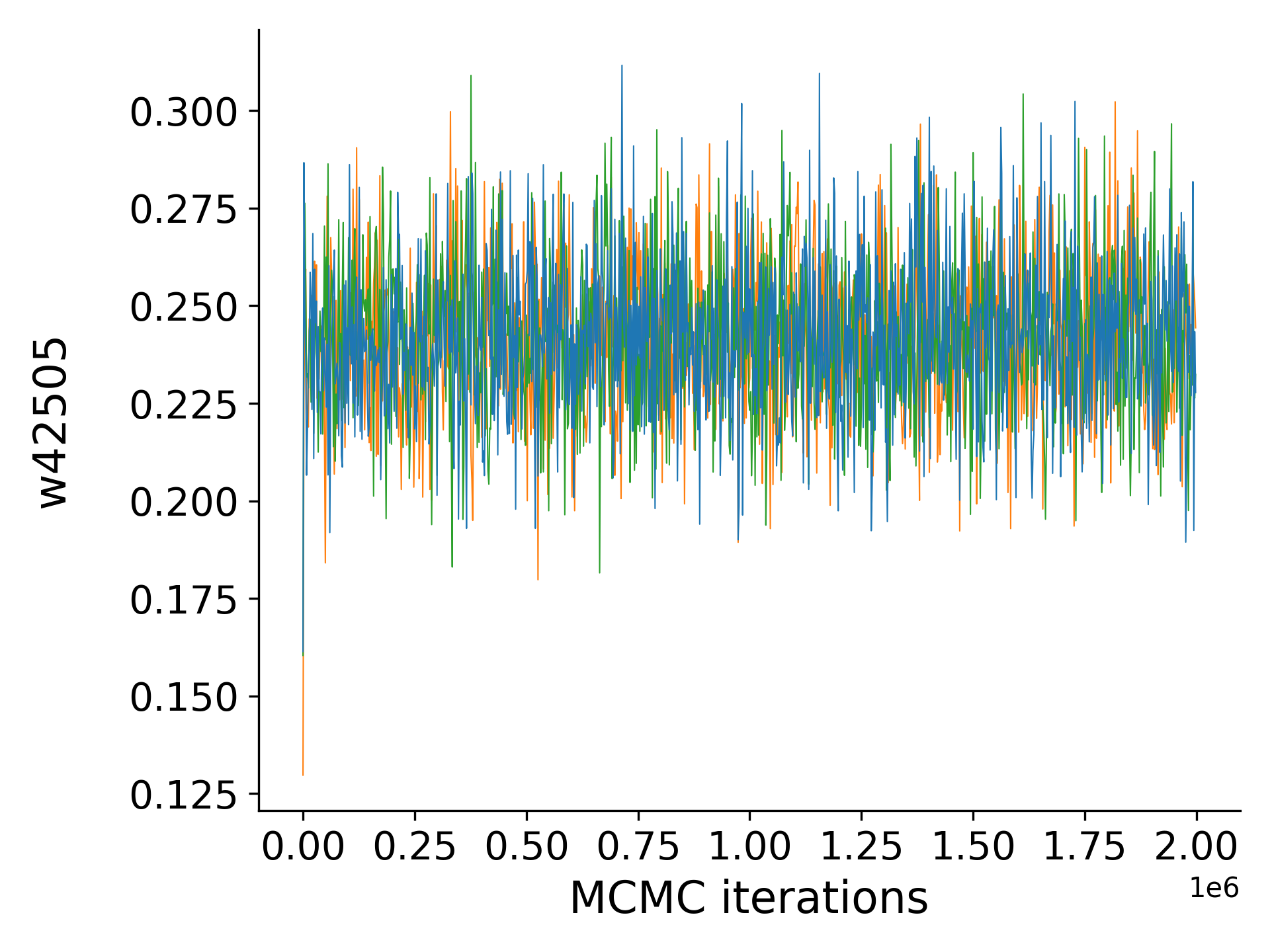}
        \caption{}
    \end{subfigure}

    \begin{subfigure}{0.45\textwidth}
        \centering
        \includegraphics[width=\textwidth]{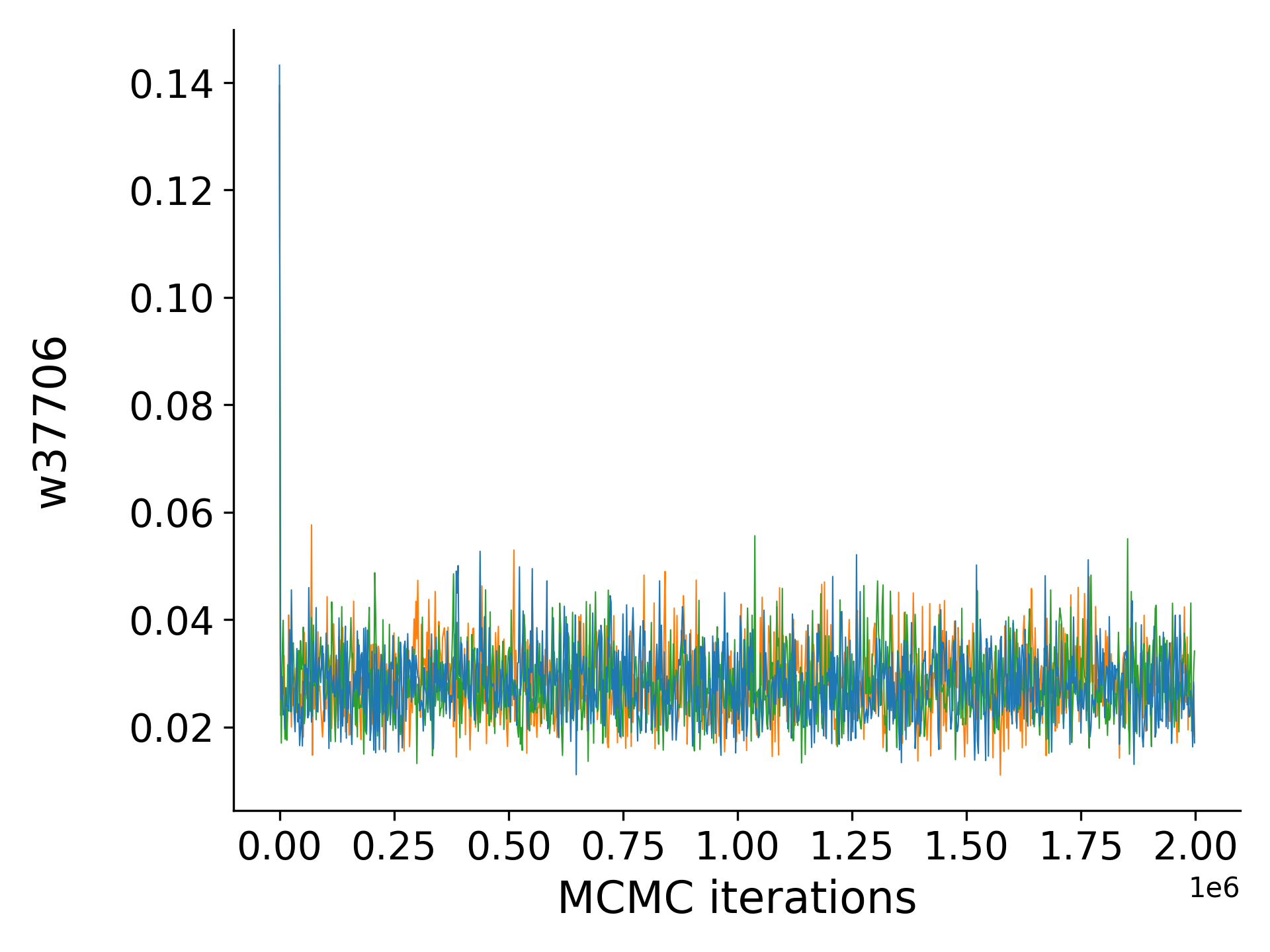}
        \caption{}
    \end{subfigure}
    \begin{subfigure}{0.45\textwidth}
        \centering
        \includegraphics[width=\textwidth]{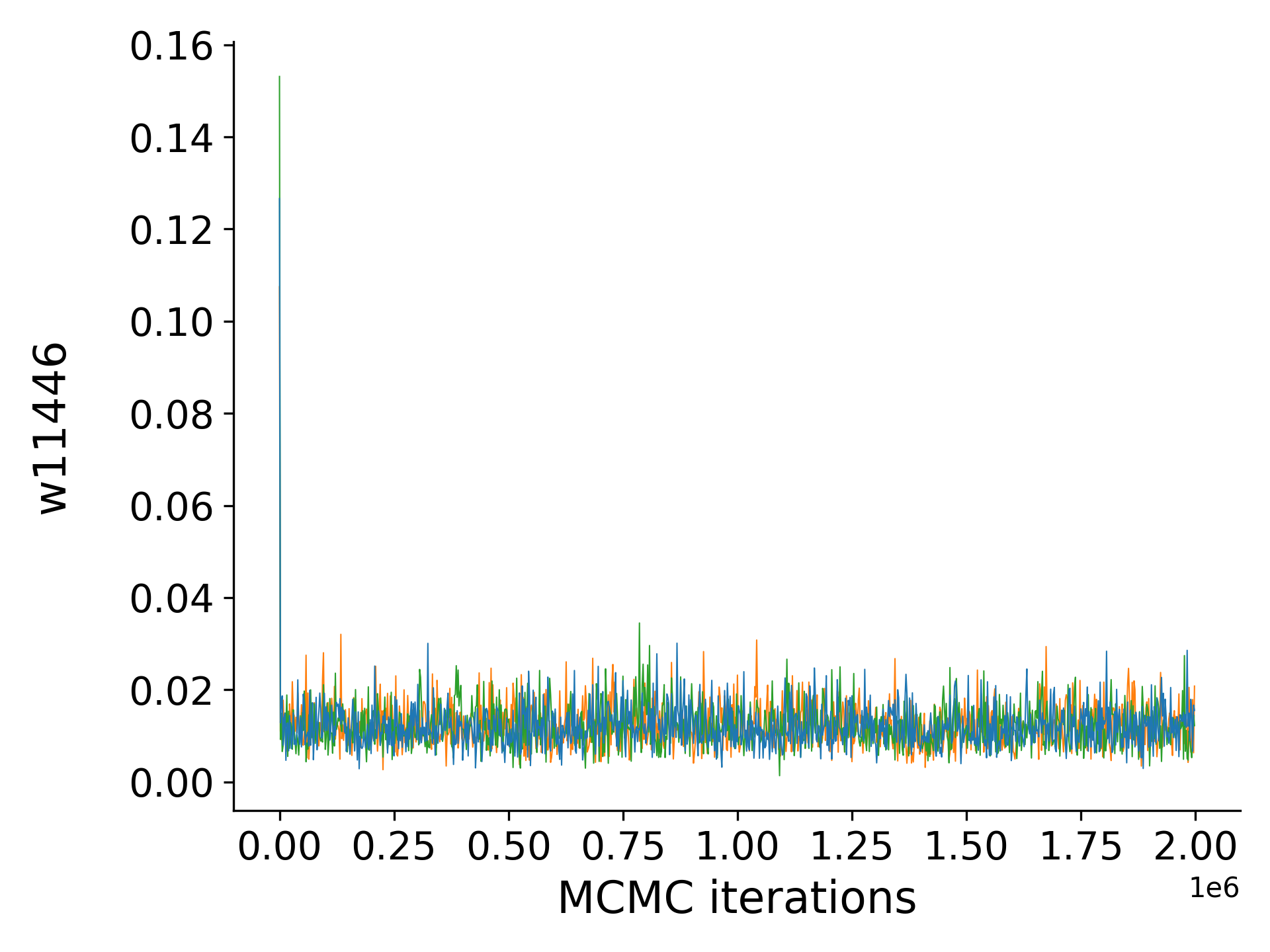}
        \caption{}
    \end{subfigure}
       \caption{MCMC traceplot of four weights for the TwitterCrawl subgraph. The degree of the corresponding node is (a) 309, (b) 151, (c) 18 and (d) 8.}
    \label{fig:traceploweightsTwitterCrawl}
\end{figure}

\begin{figure}[t]
    \centering
    \begin{subfigure}{0.45\textwidth}
        \centering
        \includegraphics[width=\textwidth]{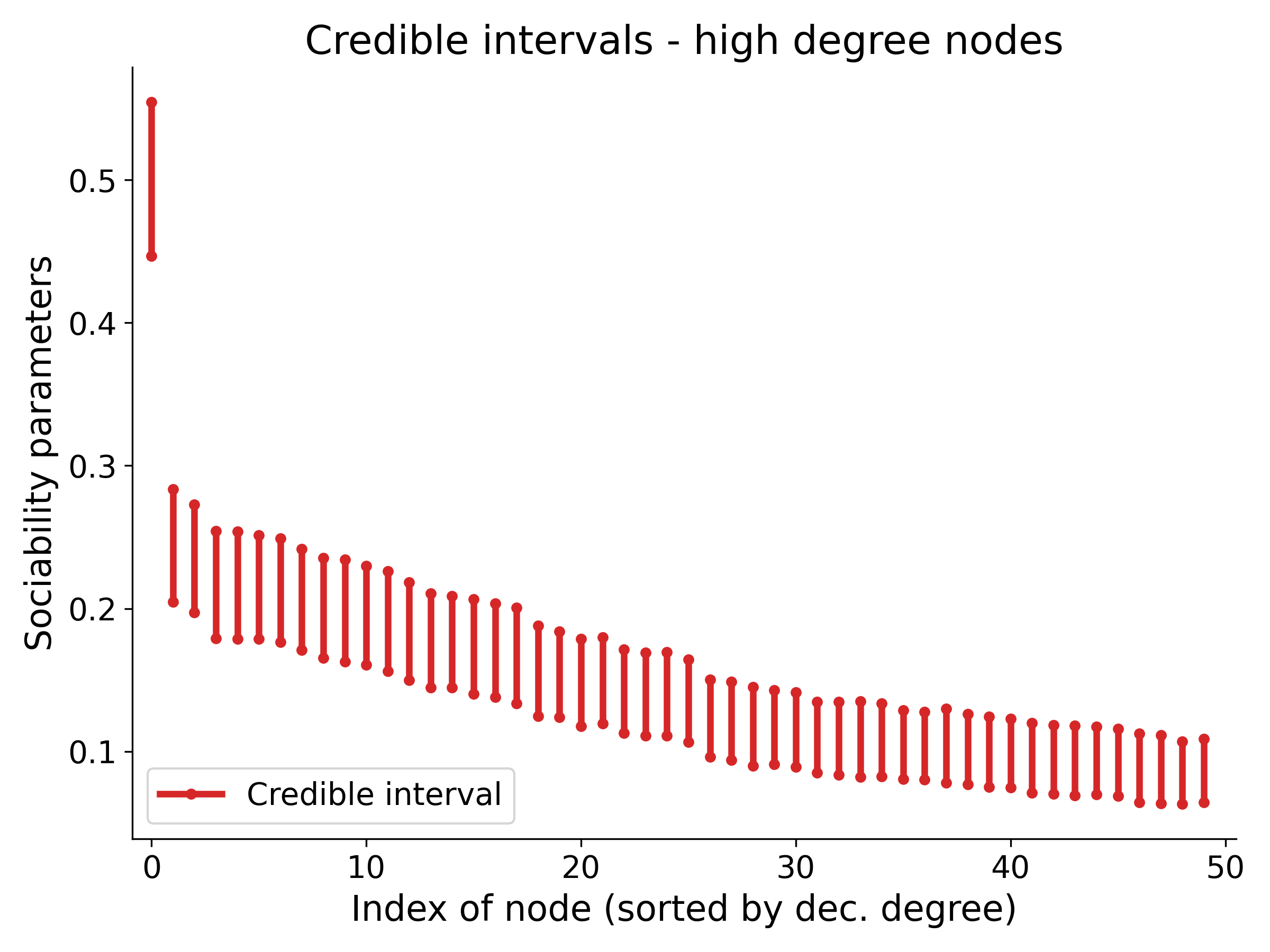}
        \caption{}
    \end{subfigure}
    \begin{subfigure}{0.45\textwidth}
        \centering
        \includegraphics[width=\textwidth]{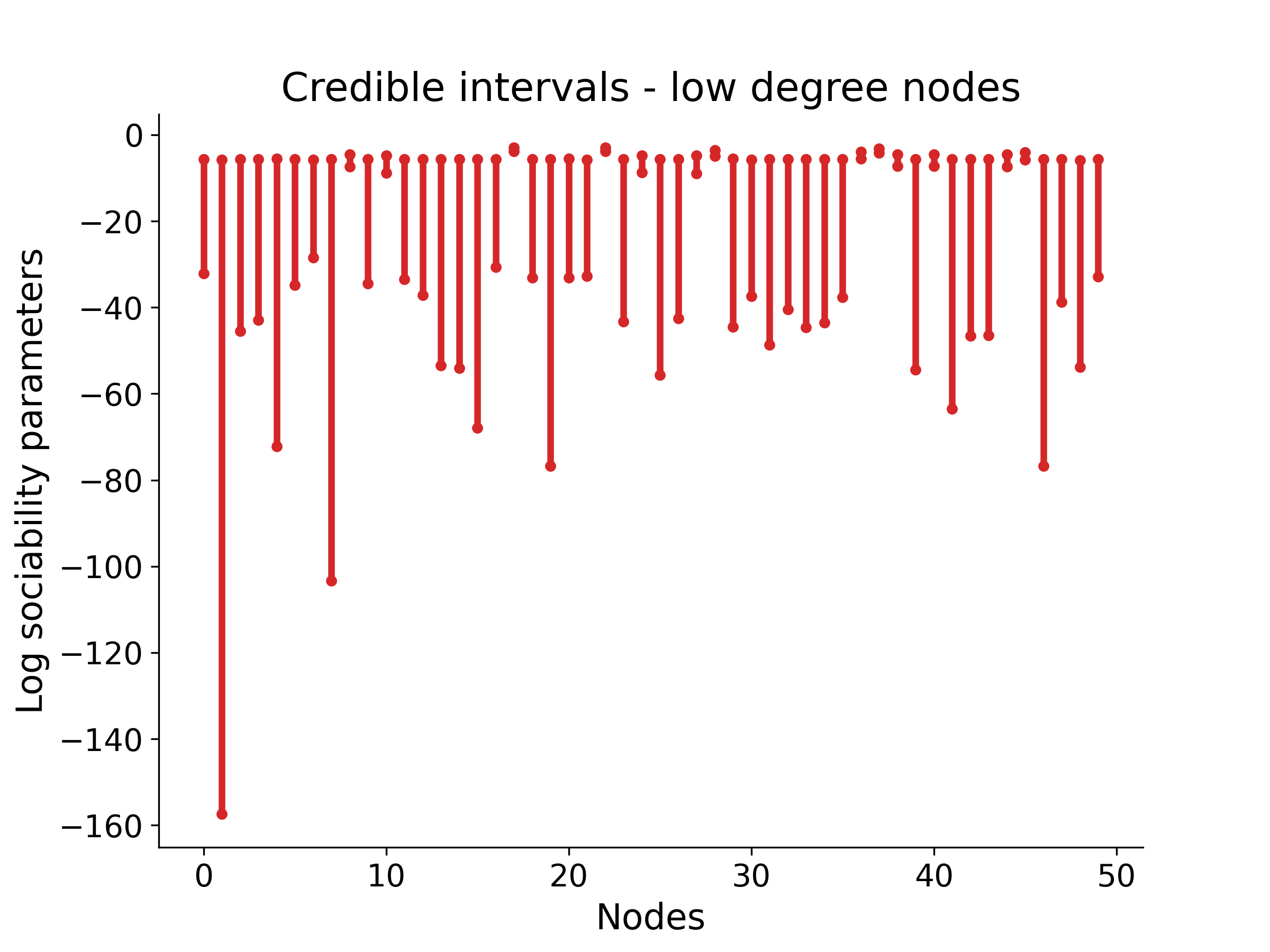}
        \caption{}
    \end{subfigure}
        \caption{95\% posterior intervals of (a) the sociability parameters $w_i$ of the 50 nodes with highest degree and (b) the log-sociability parameters  $\log(w_i)$ of the 50 nodes with lowest degree, for the TwitterCrawl subgraph.}
            \label{fig:credibleintervalTwitterCrawl}
\end{figure}

In this section, we provide additional details on the inference performed on real-world data in \cref{art:real}. To split the network data between a training set and a test set, we use $p$-sampling, as described by \citet{Veitch2019}. Statistics of the sampled graphs used for inference are given in \cref{tab:sizepsample}. The inference was performed on a standard computer (13th Gen Intel Core i5-1350P), each chain taking about 24 hours to be computed. For each graph, we run our inference algorithm for 2 million iterations. Apart from the number of iterations, we follow exactly the same procedure as for the synthetic data. We further illustrate the inference procedure on the three datasets. Unlike in the synthetic case, the true parameter values are unknown. We present the trace plots of the hyperparameters (\cref{fig:traceplotFlickr}, \cref{fig:traceplotDouban}, \cref{fig:traceplotTwitterCrawl}) and selected weights (\cref{fig:traceploweightsFlickr}, \cref{fig:traceploweightsDouban}, \cref{fig:traceploweightsTwitterCrawl}), as well as the credible intervals for the 50 nodes with the highest degree and the 50 nodes with the lowest degree (\cref{fig:credibleintervalFlickr}, \cref{fig:credibleintervalDouban}, \cref{fig:credibleintervalTwitterCrawl}).
We compute the same Gelman-Rubin (GR) convergence diagnostic as for the synthetic data for each of our 3 datasets in \cref{tab:GRreal}. In all three cases, the values suggest the convergence of the algorithm.

\begin{table}[ht]
    \centering
    \begin{tabular}{c c c}
        \begin{subtable}{0.3\textwidth}
            \centering
            \begin{tabular}{l c}
                Parameter & Value \\
                \hline
                $\beta$ & 1.005 \\
                $c$ &  1.009\\
                $\eta$ &  1.010\\
                multi & 1.000\\
                \hline
            \end{tabular}
            \caption{Flickr}
        \end{subtable}
        &
        \begin{subtable}{0.3\textwidth}
            \centering
            \begin{tabular}{l c}
                Parameter & Value \\
                \hline
                $\beta$ & 1.058 \\
                $c$ & 1.070 \\
                $\eta$ & 1.069 \\
                multi &1.003 \\
                \hline
            \end{tabular}
            \caption{Douban}
        \end{subtable}
        &
        \begin{subtable}{0.3\textwidth}
            \centering
            \begin{tabular}{l c}
                Parameter & Value \\
                \hline
                $\beta$ &  1.028\\
                $c$ &  1.044\\
                $\eta$ &  1.046\\
                multi & 1.000\\
                \hline
            \end{tabular}
            \caption{TwitterCrawl}
        \end{subtable}
    \end{tabular}
    \caption{Values of the Gelman-Rubin convergence diagnostic for the 3 datasets.}
    \label{tab:GRreal}
\end{table}

\section{Useful identities and results} 
\label{sup:identity}

\subsection{Gautschi's inequality}
For any $x>0$, any $s\in(0,1),$%
\[
x^{1-s}<\frac{\Gamma(x+1)}{\Gamma(x+s)}<(x+1)^{1-s}.
\]
Taking $x\in(0,1)$ and $s=1-x$, we obtain, for any $x\in(0,1)$%
\[
\frac{1}{x^{1-x}}<\Gamma(x)<\frac{1}{(x+1)^{1-x}}.
\]
For $x>0$,
$$
x\leq 1+ x \log (x).
$$

\subsection{Lambert function}

The LambertW function is defined to be the multivalued inverse of the function $f(y)=y e^{y}$.  When dealing with real number only, the two branches suffice: for real number $x$ and $y$ the equation $f(y)=x$  is given by
\begin{itemize}
\item $y=\LambertW_{0}(x)$ for $x\geq0$,
\item $y=\LambertW_{0}(x)$ and $y=\LambertW_{-1}(x)$ if $-\frac{1}{e}<x<0$,
\end{itemize}
where $\LambertW_{k}$ is the $k$th branch of the Lambert function. The real branches of the LambertW function can't be expressed in terms of elementary function but some numerical methods exist \citep{Loczi2022,Iacono2017,Corless1996}.

The LambertW function satisfies the following identities.

\begin{align}
& \LambertW(x)e^{\LambertW(x)}=x,\\
& \LambertW_{0}(x\log(x))=\log(x) \text{~~for~~} \frac{1}{e}\leq x,\\
& \LambertW_{-1}(x\log(x))=\log(x) \text{~~for~~} 0<x\leq \frac{1}{e},\\
& \LambertW_{0}\left(-\frac{1}{e}\right)=-1.
\end{align}

\subsection{Secondary lemmas}

\begin{proposition}\label{prop:identity1}
For any $z>0$, and $0\leq \tau<\alpha\leq 1$, we have
$$
\int_\tau^\alpha s z^s \mathrm{d}s =\left \{\begin{array}{ll}
                             \frac{\alpha^2-\tau^2}{2} & z=1 \\
                             \frac{z^{\tau}-z^{\alpha}+(\alpha z^{\alpha
}-\tau z^{\tau})\log z}{\left(  \log z\right)  ^{2}} & z\neq 1
                           \end{array}\right .
$$
\end{proposition}
\begin{proof}
For $z=1$, we have $\int_\tau^\alpha s  \mathrm{d}s=\frac{\alpha^2-\tau^2}{2}$. For $z\neq 1$, using the change of variable $u=z^s$,
\begin{align*}
\int_\tau^\alpha s z^s \mathrm{d}s&=\frac{1}{\log^{2}z}\int_{z^{\tau}}^{z^{\alpha}}\log u\mathrm{d}u\\
&  =\frac{1}{\log^{2}z}\left(  z^{\alpha}\log
z^{\alpha}-z^{\alpha}-\left[  z^{\tau}\log z^{\tau}-z^{\tau}\right]  \right)
\\
&  =\frac{z^{\tau}-z^{\alpha}+(\alpha z^{\alpha
}-\tau z^{\tau})\log z}{\left(  \log z\right)  ^{2}}.%
\end{align*}
\end{proof}

\begin{proposition}\label{prop:identity2}\label{prop:inverseCDF}
For $0\leq \tau<x<\alpha\leq 1$ and $z>0$, let
$$
F(x)=\frac{\int_\tau^x sz^s \mathrm{d}s}{\int_\tau^\alpha sz^s \mathrm{d}s}
$$
Then
$$
F(x)=\left \{\begin{array}{ll}
                             \frac{x^2-\tau^2}{\alpha^2-\tau^2} & z=1, \\
                             \frac{z^{\tau}-z^{x}+(x z^{x
}-\tau z^{\tau})\log z}{z^{\tau}-z^{\alpha}+(\alpha z^{\alpha
}-\tau z^{\tau})\log z} & z\neq 1.
                           \end{array}\right .
$$
and
$$
F^{-1}(y)=\left \{\begin{array}{ll}
                             \sqrt{(\alpha^2-\tau^2)y+\tau^2} & z=1, \\
                             \frac{\LambertW_0(c(y)/e) +1}{\log z} & z>1, \\
    \frac{\LambertW_{-1}(c(y)/e) +1}{\log z} & z<1.
\end{array}\right .
$$
where
\begin{equation}
c(y)=(z^{\tau}-z^{\alpha}+(\alpha z^{\alpha
}-\tau z^{\tau})\log z)y -z^\tau +\tau z^\tau\log z,
\label{eq:c}
\end{equation}
and $\LambertW_k$ denotes the $k$th branch of the Lambert W function.

\end{proposition}

\begin{proof}
The expression for $F$ follows from \cref{prop:identity1}. The inverse of $F$ for $z=1$ is straightforward. For $z\neq1$, $y=F(x)$ gives
\begin{align*}
c(y)&=x z^x\log z - z^x \\
&= e \times  (x\log z -1)e^{x\log z -1}.
\end{align*}
where $c(y)=(z^{\tau}-z^{\alpha}+(\alpha z^{\alpha
}-\tau z^{\tau})\log z)y -z^\tau +\tau z^\tau\log z$. Hence we have
$$c(y)/e=f(x\log z -1),$$
where $f(w)=we^w$. We have $-1<c(y)<0$
hence there are two solutions to the equation (in $x$) $c(y)/e=f(x)$, which can be expressed in terms of the $\LambertW_0$ and $\LambertW_{-1}$ functions. If $z>1$, then $x\log z -1 >-1$ for all $x$, hence it is expressed in terms of $\LambertW_0$. Otherwise $x\log z -1 <-1$ hence it is expressed in terms of $\LambertW_{-1}$. We therefore obtain
$$
x=\left \{\begin{array}{cc}
    \frac{\LambertW_0(c(y)/e) +1}{\log z} & z>1, \\
    \frac{\LambertW_{-1}(c(y)/e) +1}{\log z} & z<1.
  \end{array}\right .
$$
\end{proof}

\begin{proposition}\label{prop:inverselaplace}
The inverse of the function $\psica(y;\alpha,0)=\frac{t^{\alpha}-1}{\log t^\alpha}$ is
$$
\psicainv(y;\alpha,0) = \left \{\begin{array}{cc}
\left [ -y\LambertW_{-1}\left(  -\frac{1}{y}e^{-1/y}\right)\right ]^{1/\alpha} & \text{for} t>1,\\
\left [ -y\LambertW_{0}\left(  -\frac{1}{y}e^{-1/y}\right)\right ]^{1/\alpha} & \text{for~} t\in(0,1),\\
1 & \text{for~} t=1.
  \end{array}\right .
$$
where $\LambertW_k$ denotes the $k$th branch of the Lambert W function.
\end{proposition}

\begin{proof}
Consider the case $\tau=0$. We are looking for the inverse of $\psica(t;\alpha,0)$. Let $y=\psica(t;\alpha,0)=\frac{t^{\alpha}-1}{\log t^\alpha}$.
This function is monotone increasing and has a well-defined inverse. We have%
$$
t^\alpha=e^{\frac{t^\alpha-1}{y}}%
$$
and%
$$
-\frac{t^\alpha}{y}e^{-t^\alpha/y}=-\frac{1}{y}e^{-1/y}%
$$
and so $-\frac{1}{y}e^{-1/y}=f(-\frac{t^\alpha}{y})$. Note that $-\frac{1}{y}%
e^{-1/y}\in\lbrack-1/e,0)$, when $y>0$, where the minimum is achieved for
$y=1$. The set of solution is given by $\LambertW_{0}(x)\in\lbrack-1,0)$ and
$\LambertW_{-1}(x)\in\lbrack-1,-\infty)$.

Note that $y=\psica(t;\alpha,0)<t^\alpha$ for $t>1$ as $\log(1+x)\geq \frac{x}{x+1}$ for $x>-1$.
So for $t>1$ we have
\[
t=\left [ -y\LambertW_{-1}\left(  -\frac{1}{y}e^{-1/y}\right)\right ]^{1/\alpha}  =\psicainv(y;\alpha,0).
\]
For $t\in(0,1)$ we have $y>t^\alpha$ as $x>x\ln(x)-1$ if $x\in(0,1)$ so we have
\[
t=\left [ -y\LambertW_{0}\left(  -\frac{1}{y}e^{-1/y}\right)\right ]^{1/\alpha}  =\psicainv(y;\alpha,0).
\]
Finally, for $t=1$ we have $\psica(1;\alpha,0)=1$ so $\psicainv(1;\alpha,0)=1$.
\end{proof}

\section{Background material on regularly varying functions}
\label{app:regularvariation}

Most of the background material in this section originates from the book of \cite{Bingham1987}.

\subsection{Definitions}
\begin{definition}[Slowly varying function]
A function $\ell:(0,\infty)\to(0,\infty)$ is slowly varying at infinity if for all $c>0$,
$$
\frac{\ell(cx)}{\ell(x)}\to 1\text{  as }x\to\infty.
$$
\end{definition}
\noindent Examples of slowly varying functions include $\log^a$, for $a\in\Real$, and functions converging to a constant $c>0$.

\begin{definition}[Regularly varying function]
A function $\ell:(0,\infty)\to(0,\infty)$ is regularly varying at infinity with exponent $\rho\in\Real$ if $f(x)=x^\rho\ell(x)$ for some slowly varying function $\ell$. A function $f$ is regularly varying at
$0$ if $f(1/x)$ is regularly varying at infinity that is, $f(x)=x^{-\rho}\ell(1/x)$ for some $\rho\in\Real$.
\end{definition}

\subsection{Karamata theorems}

The following propositions and corollaries relate to integrals of regularly varying functions.

\begin{proposition}[Karamata theorem] See \citep[Propositions 1.5.8 and 1.5.10]{Bingham1987}. Let $U(t)=t^\rho\ell(t)$ for some locally bounded slowly varying function $\ell$. Then
\begin{itemize}
\item[(i)] If $\rho>-1$
$$
\int_0^x U(t)dt \sim\frac{1}{1+\rho}x^{\rho +1}\ell(x)\text{  as }x\to\infty.
$$
\item[(ii)] If $\rho<-1$
$$
\int_x^\infty U(t)dt\sim -\frac{1}{1+\rho}x^{\rho +1}\ell(x)\text{  as }x\to\infty.
$$
\end{itemize}
\label{thm:RVKaramata1}
\end{proposition}
The following corollaries will be useful.

\begin{corollary}
\label{thm:RVKaramata2}
Let $U(x)=x^\alpha \ell(1/x)$ for some locally bounded slowly varying function $\ell$.
\begin{itemize}
\item[(i)] If $\alpha>-1$
$$
\int_0^x U(t)\mathrm{d}t\sim \frac{1}{\alpha+1}x^{1+\alpha}\ell(1/x)\text{  as }x\to 0.
$$
\item[(ii)] If $\alpha<-1$
$$
\int_x^\infty U(t)\mathrm{d}t\sim -\frac{1}{\alpha+1}x^{1+\alpha}\ell(1/x)\text{  as }x\to 0.
$$

\end{itemize}
If $\ell$ is slowly varying and $\alpha>-1$ then
$\int_{0}^{x}t^{\alpha}\ell(1/t)dt$ converges and%
\[
\frac{x^{1+\alpha}\ell(1/x)}{\int_{0}^{x}t^{\alpha}\ell(1/t)dt}\rightarrow
\alpha+1\text{ as }x\rightarrow0.
\]

\end{corollary}

\begin{proof}
Let $\rho<-1$. We have $\int
_{x}^{\infty}t^{\rho}\ell(t)\mathrm{d}t=\int_{0}^{1/x}t^{-\rho-2}\ell(1/t)\mathrm{d}t$. Writing
$\alpha=-\rho-2>-1$, we obtain, using \cref{thm:RVKaramata1}(ii) %
\[
\frac{x^{-\alpha-1}\ell(x)}{\int_{0}^{1/x}t^{\alpha}\ell(1/t)\mathrm{d}t}%
\rightarrow\rho+1\text{ as }x\rightarrow\infty,
\]
or
\[
\frac{x^{1+\alpha}\ell(1/x)}{\int_{0}^{x}t^{\alpha}\ell(1/t)dt}\rightarrow
\rho+1\text{ as }x\rightarrow0.
\]
The case $\alpha<-1$ follows similarly, using \cref{thm:RVKaramata1}(i).
\end{proof}

\subsection{Asymptotic inverse}
\label{Ann:Asymptoticinverse}

If $f$ is regularly varying then there exists a regularly varying function $g$ such that $f(g(x))\sim g(f(x))\sim x$ as $x\to\infty$. We say that $g$ is the asymptotic inverse of $f$. It is uniquely determined within asymptotic equivalence, and one version of $g$ is the generalised inverse of $f$.

\begin{theorem}{\citep[Proposition 1.5.13]{Bingham1987}}
If $\ell$ varies slowly, there exist a slowly varying function $\ell^{\#}$, unique up to asymptotic equivalent, with
$$\ell(x)\ell^{\#}(x\ell(x))\underset{x\rightarrow\infty}{\longrightarrow} 1 \text{~~and~~} \ell^{\#}(x)\ell(x\ell^{\#}(x))\underset{x\rightarrow\infty}{\rightarrow}1.$$
\end{theorem}
The function $\ell^\#$ is the Bruijn conjugate of $\ell$; $\left(\ell,\ell^{\#}\right)$ is a conjugate pair. Such pairs are of common occurence, for instance  the de\ Bruijn conjugate of $\log^{a}$ is $\log^{-a}$ for any $a\in\mathbb{R}.$
\begin{proposition}{\citep[Proposition 1.5.14]{Bingham1987}}
\label{thm:conjugaterules}
If $\left(\ell,\ell^{\#}\right)$ is a conjugate pair, $A,B,a>0$ each of the following is also a pair
$$\left(\ell(Ax),\ell^{\#}(Ax)\right),$$
$$\left(A\ell(x),A^{-1}\ell^{\#}(x)\right),$$
$$\left([(\ell(x^a)]^{\frac{1}{a}},[\ell^{\#}(x^a)]^{\frac{1}{a}}\right).$$
\end{proposition}
The following is a special case from \citep[Proposition 1.5.15 p.29]{Bingham1987}.
\begin{proposition}{\citep[Proposition 1.5.15]{Bingham1987}}\label{thm:RVInverse}
\noindent Let $a,b>0$. Let $f(x)\sim
x^{ab}\ell^{a}(x^{b})$ where $\ell$ is slowly varying, and let $g$ be the
asymptotic inverse of $f$. Then%
$$g(x)\sim x^{\frac{1}{ab}}[\ell^{\#}]^{\frac{1}{b}}(x^{\frac{1}{a}})\text{ as }x\to\infty.$$
\end{proposition}

\subsection{Tauberian theorem}

The following theorem is a variation of \citet[Theorem 1.7.6 p.46]{Bingham1987}, where the two limits at 0 and infinity are exchanged. The proof is similar. See also \citet[Chapter XIII]{Feller1971}.

\begin{theorem}[Tauberian theorem]\label{thm:RVTauberian}
Assume $U(x)\geq0$, $c\geq0$, $\rho>-1$, $\widehat{U}(s)=s\int_{0}^{\infty
}e^{-sx}U(x)\mathrm{d}x$ convergent for $s>0$, and $\ell$ a slowly varying function.
Then%
\[
U(x)\sim cx^{\rho}\ell(1/x)/\Gamma(1+\rho)\text{ as }x\rightarrow0
\]
implies%
\[
\widehat{U}(s)\sim cs^{-\rho}\ell(s)\text{ as }s\rightarrow\infty.
\]

\end{theorem}

\begin{proof}
Write $V(x)=\int_{0}^{x}U(y)\mathrm{d}y$ (this is finite for any $x$ as $\widehat
{U}(s)$ is convergent for any $s$), then $V$ is non-decreasing and by
\cref{thm:RVKaramata2}
\[
V(x)\sim\frac{c}{\rho+1}x^{\rho+1}\ell(1/x)/\Gamma(1+\rho)\text{ as
}x\rightarrow0.
\]
Then by \citet[Theorem 1.7.1, p.38]{Bingham1987}, this is equivalent to
\[
\widehat{V}(s)=\int_{0}^{\infty}e^{-sx}\mathrm{d}V(x)\sim cs^{-\rho-1}\ell(s)\text{ as
}s\rightarrow\infty.
\]
Finally, note that $\widehat{V}(s)=\frac{\widehat{U}(s)}{s}$. Thus the above equation is
equivalent to%
\[
\widehat{U}(s)\sim cs^{-\rho}\ell(s)\text{ as }s\rightarrow\infty.
\]
\end{proof}

\section{Size-biased representation of a CRM}
\label{app:sizebiased}

\begin{lemma}
Let $\rho_0$ be some L\'{e}vy intensity, with Laplace exponent $\psi_0$ and
inverse Laplace exponent $\psi^{-1}$. Let $\rho_{1}(w)=\frac{\eta}{c}%
\rho_0(\frac{w}{c})e^{-\beta w/c}$, with Laplace exponent $\psi_{1}$. Then
\begin{align*}
\psi_{1}(t)  &  =\eta\left[  \psi_0(ct+\beta)-\psi_0(\beta)\right] \\
\psi_{1}^{-1}(t)  &  =\frac{1}{c}\left[  \psi_0^{-1}\left(  \frac{t}{\eta}%
+\psi_0(\beta)\right)  -\beta\right].
\end{align*}
\end{lemma}

\begin{proof}
\begin{align*}
\psi_{1}(t) &  =\int_{0}^{\infty}(1-e^{-wt})\rho_{1}(w)\mathrm{d}w\\
&  =\frac{\eta}{c}\int_{0}^{\infty}(1-e^{-wt})\rho_0(\frac{w}{c})e^{-\beta
w/c}\mathrm{d}w\\
&  =\eta\int_{0}^{\infty}(1-e^{-uct})\rho_0(u)e^{-\beta u}\mathrm{d}u\\
&  =\eta\left[  \int_{0}^{\infty}(1-e^{-u(\beta+ct)})\rho_0(u)\mathrm{d}u-\int
_{0}^{\infty}(1-e^{-\beta u})\rho_0(u)\mathrm{d}u\right]  \\
&  =\eta\left[  \psi_0(ct+\beta)-\psi_0(\beta)\right].
\end{align*}

\end{proof}
The following proposition follows from \citet{Perman1992}; see also \citet{Campbell2019,Lee2023}.

\begin{proposition}\citep{Perman1992}\label{thm:sizebiased}
Consider a homogeneous CRM  $G\sim\CRM(\rho,H)$ on $\Theta$ with L\'{e}vy intensity $\rho(w)=\frac{\eta}{c}\rho_0
(\frac{w}{c})e^{-\beta w/c}$, where $\rho_0$ is some base L\'evy intensity, with associated Laplace exponent $\psi_0(t)=\int_{0}^{\infty
}(1-e^{-wt})\rho_0(w)\mathrm{d}w$ and tilted moment function $\kappa_0(m,z)=\int_0^\infty w^m e^{-zw}\rho_0(w)\mathrm{d}w$. The Laplace exponent of the CRM $G$ takes the form $$\psi(t)=\eta\left[
\psi_0(ct+\beta)-\psi_0(\beta)\right].$$ 
Let $\psi_0^{-1}$ be the inverse of $\psi_0$. Let $\xi_1,\xi_2,\ldots$ be the points of a unit-rate Poisson process on $(0,\infty)$. The
size-biased representation of the CRM $$G=\sum_{j\geq1}W_{j}\delta_{\theta_{j}%
}$$ is given, for $j\geq 1$, by $\theta_j\sim H$ and %
\begin{align*}
T_{j}  &  =\psi_0^{-1}\left(  \frac{\xi_{j}}{\eta}+\psi_0(\beta)\right)  -\beta\\
W_{j}^{\prime}|\{T_{j}   =t\}&\sim p(w\mid t)=\frac{we^{-(t+\beta)w}\rho_0
(w)}{\kappa_0(1,t+\beta)}\\
W_{j}  &  =cW_{j}^{\prime}.
\end{align*}
\end{proposition}


\begin{thebibliography}{}

\bibitem[Ayed et~al., 2019]{Ayed2019}
Ayed, F., Lee, J., and Caron, F. (2019).
\newblock Beyond the {{Chinese Restaurant}} and {{Pitman-Yor}} processes:
  {{Statistical Models}} with {{Double Power-law Behavior}}.
\newblock In {\em International {{Conference}} on {{Machine Learning}}},
  volume~97, pages 395--404.

\bibitem[Ayed et~al., 2024]{Ayed2024}
Ayed, F., Lee, J., and Caron, F. (2024).
\newblock The {{Normal-Generalised Gamma-Pareto Process}}: {{A Novel Pure-Jump
  L{\'e}vy Process}} with {{Flexible Tail}} and {{Jump-Activity Properties}}.
\newblock {\em Bayesian Analysis}, 19(1):123--152.

\bibitem[Barab{\'a}si and Albert, 1999]{Barabasi1999a}
Barab{\'a}si, A.-L. and Albert, R. (1999).
\newblock Emergence of scaling in random networks.
\newblock {\em Science}, 286(5439):509--512.

\bibitem[Betancourt et~al., 2016]{Zanella2016}
Betancourt, B., Zanella, G., Miller, J., Wallach, H., Zaidi, A., and Steorts,
  R. (2016).
\newblock Flexible models for microclustering with application to entity
  resolution.
\newblock {\em Advances in Neural Information Processing Systems}, 29.

\bibitem[Betancourt et~al., 2022]{Betancourt2022}
Betancourt, B., Zanella, G., and Steorts, R.~C. (2022).
\newblock Random {{Partition Models}} for {{Microclustering Tasks}}.
\newblock {\em Journal of the American Statistical Association},
  117(539):1215--1227.

\bibitem[Borgs et~al., 2018]{Borgs2018}
Borgs, C., Chayes, J.~T., Cohn, H., and Holden, N. (2018).
\newblock Sparse exchangeable graphs and their limits via graphon processes.
\newblock {\em Journal of Machine Learning Research}, 18(210):1--71.

\bibitem[Brix, 1999]{Brix1999}
Brix, A. (1999).
\newblock Generalized gamma measures and shot-noise {{Cox}} processes.
\newblock {\em Advances in Applied Probability}, pages 929--953.

\bibitem[Broderick et~al., 2012]{Broderick2012}
Broderick, T., Jordan, M.~I., and Pitman, J. (2012).
\newblock Beta processes, stick-breaking and power laws.
\newblock {\em Bayesian analysis}, 7(2):439--476.

\bibitem[Cai et~al., 2016]{Cai2016}
Cai, D., Campbell, T., and Broderick, T. (2016).
\newblock Edge-exchangeable graphs and sparsity.
\newblock {\em Advances in Neural Information Processing Systems}, 29.

\bibitem[Camerlenghi et~al., 2019]{Camerlenghi2019}
Camerlenghi, F., Dunson, D.~B., Lijoi, A., Pr{\"u}nster, I., Rodr{\'i}guez, A.,
  et~al. (2019).
\newblock Latent nested nonparametric priors (with discussion).
\newblock {\em Bayesian Analysis}, 14(4):1303--1356.

\bibitem[Caron and Fox, 2017]{Caron2017}
Caron, F. and Fox, E. (2017).
\newblock Sparse graphs using exchangeable random measures.
\newblock {\em Journal of the Royal Statistical Society. Series B (Statistical
  Methodology)}, 79:1295--1366.

\bibitem[Caron et~al., 2023]{Caron2023}
Caron, F., Panero, F., and Rousseau, J. (2023).
\newblock On sparsity, power-law, and clustering properties of graphex
  processes.
\newblock {\em Advances in Applied Probability}, 55(4):1211--1253.

\bibitem[Cont and Tankov, 2004]{Cont2004}
Cont, R. and Tankov, P. (2004).
\newblock {\em Financial Modelling with Jump Processes}, volume~2.
\newblock CRC press.

\bibitem[Crane and Dempsey, 2018]{Crane2018}
Crane, H. and Dempsey, W. (2018).
\newblock Edge exchangeable models for interaction networks.
\newblock {\em Journal of the American Statistical Association},
  113:1311--1326.

\bibitem[Devroye, 2009]{Devroye2009}
Devroye, L. (2009).
\newblock Random variate generation for exponentially and polynomially tilted
  stable distributions.
\newblock {\em ACM Transactions on Modeling and Computer Simulation},
  19(4):1--20.

\bibitem[Di~Benedetto et~al., 2021]{DiBenedetto2021}
Di~Benedetto, G., Caron, F., and Teh, Y.~W. (2021).
\newblock Nonexchangeable random partition models for microclustering.
\newblock {\em The Annals of Statistics}, 49(4):1931--1957.

\bibitem[Favaro et~al., 2013]{Favaro2013}
Favaro, S., Lijoi, A., and Pr{\"u}nster, I. (2013).
\newblock Conditional formulae for {{Gibbs-type}} exchangeable random
  partitions.
\newblock {\em The Annals of Applied Probability}, 23(5):1721--1754.

\bibitem[Gelman and Rubin, 1992]{Gelman1992}
Gelman, A. and Rubin, D.~B. (1992).
\newblock Inference from {{Iterative Simulation Using Multiple Sequences}}.
\newblock {\em Statistical Science}, 7(4):457--472.

\bibitem[Gnedin et~al., 2007]{Gnedin2007}
Gnedin, A., Hansen, B., and Pitman, J. (2007).
\newblock Notes on the occupancy problem with infinitely many boxes: General
  asymptotics and power laws.
\newblock {\em Probab. Surv}, 4(146-171):88.

\bibitem[Griffin and Leisen, 2018]{Griffin2018}
Griffin, J. and Leisen, F. (2018).
\newblock Modelling and {{Computation Using NCoRM Mixtures}} for {{Density
  Regression}}.
\newblock {\em Bayesian Analysis}, 13(3):897--916.

\bibitem[Griffin and Leisen, 2017]{Griffin2017}
Griffin, J.~E. and Leisen, F. (2017).
\newblock Compound random measures and their use in {{Bayesian}}
  non-parametrics.
\newblock {\em Journal of the Royal Statistical Society. Series B (Statistical
  Methodology)}, 79(2):525--545.

\bibitem[Hjort, 1990]{Hjort1990}
Hjort, N. (1990).
\newblock Nonparametric {{Bayes}} estimators based on beta processes in models
  for life history data.
\newblock {\em The Annals of Statistics}, 18(3):1259--1294.

\bibitem[Hofert, 2011]{Hofert2011}
Hofert, M. (2011).
\newblock Sampling exponentially tilted stable distributions.
\newblock {\em ACM Transactions on Modeling and Computer Simulation (TOMACS)},
  22(1):3.

\bibitem[Hougaard, 1986]{Hougaard1986}
Hougaard, P. (1986).
\newblock Survival models for heterogeneous populations derived from stable
  distributions.
\newblock {\em Biometrika}, 73(2):387--396.

\bibitem[James, 2002]{James2002}
James, L.~F. (2002).
\newblock Poisson process partition calculus with applications to exchangeable
  models and {{Bayesian}} nonparametrics.
\newblock {\em arXiv preprint math/0205093}.

\bibitem[James et~al., 2009]{James2009}
James, L.~F., Lijoi, A., and Pr{\"u}nster, I. (2009).
\newblock Posterior analysis for normalized random measures with independent
  increments.
\newblock {\em Scandinavian Journal of Statistics}, 36(1):76--97.

\bibitem[Janson, 2018]{Janson2018}
Janson, S. (2018).
\newblock On {{Edge Exchangeable Random Graphs}}.
\newblock {\em Journal of Statistical Physics}, 173(3):448--484.

\bibitem[Kingman, 1967]{Kingman1967}
Kingman, J. (1967).
\newblock Completely random measures.
\newblock {\em Pacific Journal of Mathematics}, 21(1):59--78.

\bibitem[Kingman, 1993]{Kingman1993}
Kingman, J. (1993).
\newblock {\em Poisson Processes}, volume~3.
\newblock Oxford University Press, USA.

\bibitem[Kunegis, 2013]{konect}
Kunegis, J. (2013).
\newblock {{KONECT}} -- {{The Koblenz Network Collection}}.
\newblock In {\em Proc. {{Int}}. {{Conf}}. on {{World Wide Web Companion}}},
  pages 1343--1350.

\bibitem[Lee et~al., 2023]{Lee2023}
Lee, J., Miscouridou, X., and Caron, F. (2023).
\newblock A unified construction for series representations and finite
  approximations of completely random measures.
\newblock {\em Bernoulli. Official Journal of the Bernoulli Society for
  Mathematical Statistics and Probability}, 29(3):2142--2166.

\bibitem[Lijoi et~al., 2007]{Lijoi2007}
Lijoi, A., Mena, R.~H., and Pr{\"u}nster, I. (2007).
\newblock Controlling the reinforcement in {{Bayesian}} non-parametric mixture
  models.
\newblock {\em Journal of the Royal Statistical Society: Series B (Statistical
  Methodology)}, 69(4):715--740.

\bibitem[Lijoi and Pr{\"u}nster, 2010]{Lijoi2010}
Lijoi, A. and Pr{\"u}nster, I. (2010).
\newblock Models beyond the {{Dirichlet}} process.
\newblock In Hjort, N.~L., Holmes, C., M{\"u}ller, P., and Walker, S.~G.,
  editors, {\em Bayesian Nonparametrics}. Cambridge University Press.

\bibitem[Mislove et~al., 2007]{Mislove2007a}
Mislove, A., Marcon, M., Gummadi, K.~P., Druschel, P., and Bhattacharjee, B.
  (2007).
\newblock Measurement and analysis of online social networks.
\newblock In {\em Proceedings of the 7th {{ACM SIGCOMM}} Conference on
  {{Internet}} Measurement}, pages 29--42, San Diego California USA. ACM.

\bibitem[Naik et~al., 2021]{Naik2021}
Naik, C., Caron, F., and Rousseau, J. (2021).
\newblock Sparse networks with core-periphery structure.
\newblock {\em Electronic Journal of Statistics}, 15(1):1814--1868.

\bibitem[{Nieto-Barajas} et~al., 2004]{Nieto-Barajas2004}
{Nieto-Barajas}, L.~E., Pr{\"u}nster, I., and Walker, S.~G. (2004).
\newblock Normalized random measures driven by increasing additive processes.
\newblock {\em The Annals of Statistics}, 32(6):2343--2360.

\bibitem[Perman et~al., 1992]{Perman1992}
Perman, M., Pitman, J., and Yor, M. (1992).
\newblock Size-biased sampling of {{Poisson}} point processes and excursions.
\newblock {\em Probability Theory and Related Fields}, 92(1):21--39.

\bibitem[Pitman, 2003]{Pitman2003}
Pitman, J. (2003).
\newblock Poisson-kingman partitions.
\newblock In Goldstein, D.~R., editor, {\em Statistics and Science: A
  Festschrift for Terry Speed}, volume Volume 40 of {\em Lecture
  {{Notes}}--{{Monograph}} Series}, pages 1--34. Institute of Mathematical
  Statistics, Beachwood, OH.

\bibitem[Regazzini et~al., 2003]{Regazzini2003}
Regazzini, E., Lijoi, A., and Pr{\"u}nster, I. (2003).
\newblock Distributional results for means of normalized random measures with
  independent increments.
\newblock {\em The Annals of Statistics}, 31(2):560--585.

\bibitem[Rossi and Ahmed, 2015]{nr}
Rossi, R.~A. and Ahmed, N.~K. (2015).
\newblock The {{Network Data Repository}} with {{Interactive Graph Analytics}}
  and {{Visualization}}.
\newblock In {\em {{AAAI}}}.

\bibitem[Rossi et~al., 2014]{Rossi2014}
Rossi, R.~A., Gleich, D.~F., Gebremedhin, A.~H., and Patwary, M. M.~A. (2014).
\newblock Fast maximum clique algorithms for large graphs.
\newblock In {\em Proceedings of the 23rd {{International Conference}} on
  {{World Wide Web}}}, {{WWW}} '14 {{Companion}}, pages 365--366, New York, NY,
  USA. Association for Computing Machinery.

\bibitem[Teh and Gorur, 2009]{Teh2009}
Teh, Y. and Gorur, D. (2009).
\newblock Indian buffet processes with power-law behavior.
\newblock {\em Advances in Neural Information Processing systems}, 22.

\bibitem[Thibaux and Jordan, 2007]{Thibaux2007}
Thibaux, R. and Jordan, M.~I. (2007).
\newblock Hierarchical beta processes and the {{Indian}} buffet process.
\newblock In {\em Proceedings of the 11th International Conference on
  Artificial Intelligence and Statistics ({{AISTATS}}'07)}, pages 564--571.

\bibitem[Todeschini et~al., 2020]{Todeschini2020}
Todeschini, A., Miscouridou, X., and Caron, F. (2020).
\newblock Exchangeable random measures for sparse and modular graphs with
  overlapping communities.
\newblock {\em Journal of the Royal Statistical Society: Series B (Statistical
  Methodology)}, 82(2):487--520.

\bibitem[Van Der~Hofstad, 2024]{VanDerHofstad2024}
Van Der~Hofstad, R. (2024).
\newblock {\em Random {{Graphs}} and {{Complex Networks}}}, volume~2.
\newblock Cambridge University Press, 1 edition.

\bibitem[Vats and Knudson, 2021]{Vats2021}
Vats, D. and Knudson, C. (2021).
\newblock Revisiting the {{Gelman}}--{{Rubin Diagnostic}}.
\newblock {\em Statistical Science}, 36(4):518--529.

\bibitem[Veitch and Roy, 2015]{Veitch2015}
Veitch, V. and Roy, D.~M. (2015).
\newblock The {{Class}} of {{Random Graphs Arising}} from {{Exchangeable Random
  Measures}}.
\newblock {\em arXiv:1512.03099}.

\bibitem[Veitch and Roy, 2019]{Veitch2019}
Veitch, V. and Roy, D.~M. (2019).
\newblock Sampling and estimation for (sparse) exchangeable graphs.
\newblock {\em The Annals of Statistics}, 47(6):3274 -- 3299.

\bibitem[Zafarani and Liu, 2014]{Zafarani2014}
Zafarani, R. and Liu, H. (2014).
\newblock Users {{Joining Multiple Sites}}: {{Distributions}} and {{Patterns}}.
\newblock {\em Proceedings of the International AAAI Conference on Web and
  Social Media}, 8(1):635--638.

\end{thebibliography}

\begin{thebibliography}{}

\bibitem[Barab{\'a}si and Albert, 1999]{Barabasi1999a}
Barab{\'a}si, A.-L. and Albert, R. (1999).
\newblock Emergence of scaling in random networks.
\newblock {\em Science}, 286(5439):509--512.

\bibitem[Betancourt, 2018]{Betancourt2018}
Betancourt, M. (2018).
\newblock A {{Conceptual Introduction}} to {{Hamiltonian Monte Carlo}}.
\newblock {\em arXiv:1701.02434}.

\bibitem[Bingham et~al., 1987]{Bingham1987}
Bingham, N.~H., Goldie, C.~M., and Teugels, J.~L. (1987).
\newblock {\em Regular Variation}, volume~27.
\newblock Cambridge university press.

\bibitem[Campbell et~al., 2019]{Campbell2019}
Campbell, T., Huggins, J.~H., How, J.~P., and Broderick, T. (2019).
\newblock Truncated random measures.
\newblock {\em Bernoulli. Official Journal of the Bernoulli Society for
  Mathematical Statistics and Probability}, 25(2):1256--1288.

\bibitem[Caron and Fox, 2017]{Caron2017}
Caron, F. and Fox, E. (2017).
\newblock Sparse graphs using exchangeable random measures.
\newblock {\em Journal of the Royal Statistical Society. Series B (Statistical
  Methodology)}, 79:1295--1366.

\bibitem[Caron et~al., 2023]{Caron2023}
Caron, F., Panero, F., and Rousseau, J. (2023).
\newblock On sparsity, power-law, and clustering properties of graphex
  processes.
\newblock {\em Advances in Applied Probability}, 55(4):1211--1253.

\bibitem[Corless et~al., 1996]{Corless1996}
Corless, R.~M., Gonnet, G.~H., Hare, D. E.~G., Jeffrey, D.~J., and Knuth, D.~E.
  (1996).
\newblock On the {{LambertW}} function.
\newblock {\em Advances in Computational Mathematics}, 5(1):329--359.

\bibitem[Feller, 1971]{Feller1971}
Feller, W. (1971).
\newblock {\em An Introduction to Probability Theory and Its Applications},
  volume~2.
\newblock John Wiley \& Sons.

\bibitem[Gelman and Rubin, 1992]{Gelman1992}
Gelman, A. and Rubin, D.~B. (1992).
\newblock Inference from {{Iterative Simulation Using Multiple Sequences}}.
\newblock {\em Statistical Science}, 7(4):457--472.

\bibitem[Iacono and Boyd, 2017]{Iacono2017}
Iacono, R. and Boyd, J.~P. (2017).
\newblock New approximations to the principal real-valued branch of the
  {{Lambert W-function}}.
\newblock {\em Advances in Computational Mathematics}, 43(6):1403--1436.

\bibitem[Lee et~al., 2023]{Lee2023}
Lee, J., Miscouridou, X., and Caron, F. (2023).
\newblock A unified construction for series representations and finite
  approximations of completely random measures.
\newblock {\em Bernoulli. Official Journal of the Bernoulli Society for
  Mathematical Statistics and Probability}, 29(3):2142--2166.

\bibitem[L{\'o}czi, 2022]{Loczi2022}
L{\'o}czi, L. (2022).
\newblock Guaranteed- and high-precision evaluation of the {{Lambert W}}.
\newblock {\em Applied Mathematics and Computation}, 433:127406.

\bibitem[Neal, 2011]{Neal2011}
Neal, R.~M. (2011).
\newblock {\em {{MCMC Using Hamiltonian Dynamics}}}, pages 113--162.
\newblock {Chapman and Hall/CRC}, New York, 1 edition.

\bibitem[Perman et~al., 1992]{Perman1992}
Perman, M., Pitman, J., and Yor, M. (1992).
\newblock Size-biased sampling of {{Poisson}} point processes and excursions.
\newblock {\em Probability Theory and Related Fields}, 92(1):21--39.

\bibitem[Vats and Knudson, 2021]{Vats2021}
Vats, D. and Knudson, C. (2021).
\newblock Revisiting the {{Gelman}}--{{Rubin Diagnostic}}.
\newblock {\em Statistical Science}, 36(4):518--529.

\bibitem[Veitch and Roy, 2019]{Veitch2019}
Veitch, V. and Roy, D.~M. (2019).
\newblock Sampling and estimation for (sparse) exchangeable graphs.
\newblock {\em The Annals of Statistics}, 47(6):3274 -- 3299.

\end{thebibliography}
\end{document}